\documentclass[10pt]{article}

\usepackage{geometry}
\usepackage{graphics,amssymb,theorem,latexsym, bbm, amsmath, subcaption, graphicx, multirow}
\usepackage[T1]{fontenc}
\usepackage{lmodern}
\theorembodyfont{\itshape} 
                           \newtheorem{cor}{Corollary}
                           \newtheorem{prop}{Proposition}
                           
\theorembodyfont{\rmfamily} \newtheorem{exa}{\bf Example}

\UseRawInputEncoding
\providecommand{\keywords}[1]
{
  \textbf{\textit{Keywords---}} #1
}

\begin{document}
\newcommand{\ncm}{\newcommand}
\ncm{\bfm}[1]{\mbox{\boldmath $#1$}}
\ncm{\sbfm}[1]{\mbox{\scriptsize\boldmath $#1$}}
\ncm{\scr}[1]{\mbox{\scriptsize #1}}
\ncm{\scrmath}[1]{\mbox{\scriptsize $#1$}}
\ncm{\bfmscr}[1]{\mbox{\scriptsize{\boldmath $#1$}}}

\ncm{\R}{{\mathbb{R}}}
\ncm{\Z}{{\mathbb{Z}}}
\ncm{\T}{{\mathbb{T}}}
\ncm{\Smath}{{\mathbb{S}}}
\ncm{\N}{{\mathbb{N}}}
\ncm{\C}{{\mathbb{C}}}
\ncm{\A}{{\mathbb{A}}}
\ncm{\amath}{\bfm{a}}
\ncm{\V}{{\mathbb{V}}}
\ncm{\Hap}{{\mathbb{H}}}
\ncm{\MMD}{{\mathbb{MD}}}

\ncm{\cA}{{\cal A}}
\ncm{\cB}{{\cal B}}
\ncm{\cC}{{\cal C}}
\ncm{\calF}{{\cal F}}
\ncm{\cD}{{\cal D}}
\ncm{\cG}{{\cal G}}
\ncm{\cL}{{\cal L}}
\ncm{\cN}{{\cal N}}
\ncm{\cI}{{\cal I}}
\ncm{\cJ}{{\cal J}}
\ncm{\cH}{{\cal H}}
\ncm{\cV}{{\cal V}}
\ncm{\cW}{{\cal W}}
\ncm{\cT}{{\cal T}}
\ncm{\cX}{{\cal X}}
\ncm{\cQ}{{\cal Q}}
\ncm{\cR}{{\cal R}}
\ncm{\cS}{{\cal S}}
\ncm{\cM}{{\cal M}}
\ncm{\cU}{{\cal U}}
\ncm{\cP}{{\cal P}}
\ncm{\cZ}{{\cal Z}}
\ncm{\cO}{{\cal O}}
\ncm{\cPzer}{{\cal P}_0}
\ncm{\cPone}{{\cal P}_1}
\ncm{\cPk}{{\cP_{\mbox{\scr{known}}}}}
\ncm{\cF}{{\cal F}}
\ncm{\cE}{{\cal E}}
\ncm{\cMD}{{\cal MD}}
\ncm{\tcV}{\tilde{\cal V}}
\ncm{\cCobs}{{\cal C}_{\scr{obs}}}

\ncm{\Om}{\Omega}
\ncm{\om}{\omega}
\ncm{\va}{\varepsilon}
\ncm{\vam}{\varepsilon_{\scr{max}}}
\ncm{\de}{\delta}
\ncm{\De}{\Delta}
\ncm{\ga}{\gamma}
\ncm{\Ga}{\Gamma}
\ncm{\la}{\lambda}
\ncm{\ka}{\kappa}
\ncm{\si}{\sigma}
\ncm{\Si}{\Sigma}
\ncm{\La}{\Lambda}
\ncm{\eps}{\epsilon}

\ncm{\bY}{\bfm{Y}}
\ncm{\bA}{\bfm{A}}
\ncm{\bB}{\bfm{B}}
\ncm{\bC}{\bfm{C}}
\ncm{\bD}{\bfm{D}}
\ncm{\bF}{\bfm{F}}
\ncm{\bI}{\bfm{I}}
\ncm{\bZ}{\bfm{Z}}
\ncm{\bG}{\bfm{G}}
\ncm{\bH}{\bfm{H}}
\ncm{\bL}{\bfm{L}}
\ncm{\bP}{\bfm{P}}
\ncm{\bQ}{\bfm{Q}}
\ncm{\bR}{\bfm{R}}
\ncm{\bS}{\bfm{S}}
\ncm{\bT}{\bfm{T}}
\ncm{\bU}{\bfm{U}}
\ncm{\bV}{\bfm{V}}
\ncm{\bM}{\bfm{M}}
\ncm{\bN}{\bfm{N}}
\ncm{\bW}{\bfm{W}}
\ncm{\bX}{\bfm{X}}
\ncm{\bu}{\bfm{u}}
\ncm{\bv}{\bfm{v}}
\ncm{\bw}{\bfm{w}}
\ncm{\bwpr}{\bfm{w}^\prime}
\ncm{\bhp}{\bfm{h}^\prime}
\ncm{\bc}{\bfm{c}}
\ncm{\bd}{\bfm{d}}
\ncm{\bm}{\bfm{m}}
\ncm{\bh}{\bfm{h}}
\ncm{\bn}{\bfm{n}}
\ncm{\bb}{\bfm{b}}
\ncm{\bg}{\bfm{g}}
\ncm{\be}{\bfm{e}}
\ncm{\bl}{\bfm{l}}
\ncm{\bp}{\bfm{p}}
\ncm{\bq}{\bfm{q}}
\ncm{\br}{\bfm{r}}
\ncm{\bs}{\bfm{s}}
\ncm{\bx}{\bfm{x}}
\ncm{\by}{\bfm{y}}
\ncm{\bz}{\bfm{z}}
\ncm{\balp}{\bfm{\alpha}}
\ncm{\bbe}{\bfm{\beta}}
\ncm{\bxi}{\bfm{\xi}}
\ncm{\bth}{\bfm{\theta}}
\ncm{\bom}{\bfm{\om}}
\ncm{\bmu}{\bfm{\mu}}
\ncm{\bde}{\bfm{\de}}
\ncm{\bva}{\bfm{\va}}
\ncm{\beps}{\bfm{\eps}}
\ncm{\bga}{\bfm{\ga}}
\ncm{\bka}{\bfm{\ka}}
\ncm{\bla}{\bfm{\la}}
\ncm{\bpi}{\bfm{\pi}}
\ncm{\bpsi}{\bfm{\psi}}
\ncm{\brho}{\bfm{\rho}}
\ncm{\boldeta}{\bfm{\eta}}
\ncm{\bphi}{\bfm{\phi}}
\ncm{\bLa}{\bfm{\Lambda}}
\ncm{\bPi}{\bfm{\Pi}}
\ncm{\bSi}{\bfm{\Si}}
\ncm{\bone}{\bfm{1}}

\ncm{\sbb}{\sbfm{b}}
\ncm{\sbc}{\sbfm{c}}
\ncm{\sbd}{\sbfm{d}}
\ncm{\sbC}{\sbfm{C}}
\ncm{\sbM}{\sbfm{M}}
\ncm{\sbX}{\sbfm{X}}
\ncm{\sbw}{\sbfm{w}}
\ncm{\sbx}{\sbfm{x}}
\ncm{\sby}{\sbfm{y}}
\ncm{\subv}{\scrmath{\ubv}}
\ncm{\subw}{\scrmath{\ubw}}
\ncm{\subx}{\scrmath{\ubx}}
\ncm{\suby}{\scrmath{\uby}}
\ncm{\subX}{\scrmath{\ubX}}

\ncm{\hbe}{\hat{\beta}}
\ncm{\heta}{\hat{\eta}}
\ncm{\hth}{\hat{\theta}}
\ncm{\hLa}{\hat{\Lambda}}
\ncm{\hSi}{\hat{\Sigma}}
\ncm{\hbth}{\hat{\bth}}
\ncm{\hh}{\hat{h}}
\ncm{\hs}{\hat{s}}
\ncm{\hB}{\hat{B}}
\ncm{\hD}{\hat{D}}
\ncm{\hF}{\hat{F}}
\ncm{\hI}{\hat{I}}
\ncm{\hN}{\hat{N}}
\ncm{\hP}{\hat{P}}
\ncm{\hQ}{\hat{Q}}
\ncm{\htau}{\hat{\tau}}
\ncm{\hga}{\hat{\gamma}}
\ncm{\hla}{\hat{\lambda}}
\ncm{\hmu}{\hat{\mu}}
\ncm{\hpi}{\hat{\pi}}
\ncm{\hpsi}{\hat{\psi}}
\ncm{\hbB}{\hat{\bB}}
\ncm{\hbbe}{\hat{\bbe}}
\ncm{\hbga}{\hat{\bga}}
\ncm{\hbpsi}{\hat{\bpsi}}

\ncm{\mast}{m^\ast}
\ncm{\cast}{c^\ast}
\ncm{\fast}{f^\ast}
\ncm{\siast}{\si^\ast}
\ncm{\psiast}{\psi^\ast}
\ncm{\tsiast}{\tilde{\si}^\ast}
\ncm{\alfast}{\alpha^\ast}
\ncm{\tkaast}{\tilde{\kappa}^\ast}

\ncm{\ap}{a^\prime}
\ncm{\hp}{h^\prime}
\ncm{\ip}{i^\prime}
\ncm{\jp}{j^\prime}
\ncm{\kp}{k^\prime}
\ncm{\np}{n^\prime}
\ncm{\npr}{n^\prime}
\ncm{\qp}{q^\prime}
\ncm{\spr}{s^\prime}
\ncm{\up}{u^\prime}
\ncm{\vp}{v^\prime}
\ncm{\wpr}{w^\prime}
\ncm{\xp}{x^\prime}
\ncm{\yp}{y^\prime}
\ncm{\zp}{z^\prime}
\ncm{\Cp}{C^\prime}
\ncm{\Gp}{G^\prime}
\ncm{\Ip}{I^\prime}
\ncm{\Mp}{M^\prime}
\ncm{\Np}{N^\prime}
\ncm{\Npr}{N^\prime}
\ncm{\Tp}{T^\prime}
\ncm{\gap}{\ga^\prime}
\ncm{\phpr}{\phi^\prime}

\ncm{\wbis}{w^{\prime\prime}}

\ncm{\tih}{\tilde{h}}
\ncm{\tZ}{\tilde{Z}}
\ncm{\tA}{\tilde{A}}
\ncm{\tD}{\tilde{D}}
\ncm{\tF}{\tilde{F}}
\ncm{\tI}{\tilde{I}}
\ncm{\tN}{\tilde{N}}
\ncm{\tQ}{\tilde{Q}}
\ncm{\tY}{\tilde{Y}}
\ncm{\tmu}{\tilde{\mu}}
\ncm{\tOm}{\tilde{\Omega}}
\ncm{\tnu}{\tilde{\nu}}
\ncm{\tsi}{\tilde{\sigma}}
\ncm{\tal}{\tilde{\alpha}}
\ncm{\tbeta}{\tilde{\beta}}
\ncm{\tde}{\tilde{\delta}}
\ncm{\txi}{\tilde{\xi}}
\ncm{\tmathV}{\tilde{\V}}
\ncm{\tV}{\tilde{V}}
\ncm{\tr}{\tilde{r}}
\ncm{\tu}{\tilde{u}}
\ncm{\tw}{\tilde{w}}
\ncm{\twpr}{\tilde{w}^\prime}
\ncm{\tb}{\tilde{b}}
\ncm{\td}{\tilde{d}}
\ncm{\tp}{\tilde{p}}
\ncm{\tf}{\tilde{f}}
\ncm{\tn}{\tilde{n}}
\ncm{\tS}{\tilde{S}}
\ncm{\tL}{\tilde{L}}
\ncm{\tl}{\tilde{l}}
\ncm{\tP}{\tilde{P}}
\ncm{\tSmath}{\tilde{\mathbb{S}}}
\ncm{\tT}{\tilde{T}}
\ncm{\tK}{\tilde{K}}
\ncm{\tka}{\tilde{\ka}}
\ncm{\tom}{\tilde{\om}}
\ncm{\tva}{\tilde{\va}}
\ncm{\tla}{\tilde{\la}}
\ncm{\tpi}{\tilde{\pi}}
\ncm{\trho}{\tilde{\rho}}
\ncm{\tbom}{\tilde{\bfm{\om}}}
\ncm{\tbxi}{\tilde{\bfm{\xi}}}
\ncm{\tbrho}{\tilde{\bfm{\rho}}}
\ncm{\tbg}{\tilde{\bg}}
\ncm{\tbb}{\tilde{\bb}}
\ncm{\tbr}{\tilde{\br}}
\ncm{\tbf}{\tilde{\bfm{f}}}
\ncm{\tbD}{\tilde{\bfm{D}}}
\ncm{\tbH}{\tilde{\bfm{H}}}
\ncm{\tbone}{\tilde{\bfm{1}}}
\ncm{\tbe}{\tilde{\bfm{e}}}
\ncm{\tbbe}{\tilde{\bfm{\beta}}}

\ncm{\bae}{\bar{e}}
\ncm{\baf}{\bar{f}}
\ncm{\bah}{\bar{h}}
\ncm{\bal}{\bar{l}}
\ncm{\bam}{\bar{m}}
\ncm{\ban}{\bar{n}}
\ncm{\bap}{\bar{p}}
\ncm{\bav}{\bar{v}}
\ncm{\baw}{\bar{w}}
\ncm{\baF}{\bar{F}}
\ncm{\baZ}{\bar{Z}}
\ncm{\baY}{\bar{Y}}
\ncm{\baS}{\bar{S}}
\ncm{\baH}{\bar{H}}
\ncm{\baA}{\bar{A}}
\ncm{\baD}{\bar{D}}
\ncm{\baC}{\bar{C}}
\ncm{\baN}{\bar{N}}
\ncm{\baQ}{\bar{Q}}
\ncm{\baW}{\bar{W}}
\ncm{\bacW}{\bar{\cW}}
\ncm{\bacV}{\bar{\cV}}
\ncm{\bacR}{\bar{{\cal R}}}
\ncm{\bacP}{\bar{{\cal P}}}
\ncm{\babe}{\bar{\beta}}
\ncm{\baka}{\bar{\kappa}}
\ncm{\bamu}{\bar{\mu}}
\ncm{\banu}{\bar{\nu}}
\ncm{\bade}{\bar{\de}}
\ncm{\bala}{\bar{\la}}
\ncm{\baga}{\bar{\ga}}
\ncm{\barho}{\bar{\rho}}
\ncm{\babf}{\bar{\bfm{f}}}
\ncm{\babD}{\bar{\bfm{D}}}
\ncm{\babA}{\bar{\bfm{A}}}
\ncm{\babQ}{\bar{\bfm{Q}}}
\ncm{\babW}{\bar{\bfm{W}}}
\ncm{\babh}{\bar{\bfm{h}}}
\ncm{\babr}{\bar{\bfm{r}}}
\ncm{\babde}{\bar{\bfm{\de}}}
\ncm{\babrho}{\bar{\bfm{\rho}}}
\ncm{\babone}{\bar{\bfm{1}}}

\ncm{\chnu}{\check{\nu}}

\ncm{\uC}{\underline{C}}
\ncm{\ucX}{\underline{\cX}}
\ncm{\ubx}{\underline{\bx}}
\ncm{\ubv}{\underline{\bv}}
\ncm{\ubw}{\underline{\bw}}
\ncm{\ubX}{\underline{\bX}}
\ncm{\uby}{\underline{\by}}
\ncm{\ubY}{\underline{\bY}}

\ncm{\Lin}{\, \stackrel{\cal L} \in}
\ncm{\Leq}{\, \stackrel{\cal L} =}
\ncm{\Lto}{\, \stackrel{\cal L} \longrightarrow}
\ncm{\pto}{\, \stackrel{p} \longrightarrow}
\ncm{\asto}{\, \stackrel{\rm a.s.} \longrightarrow}
\ncm{\Cov}{\mbox{Cov}}
\ncm{\Var}{\mbox{Var}}
\ncm{\sameord}{\stackrel{\cup}{{\scriptstyle \cap}}}

\ncm{\ith}{i^{\scr{th}}}
\ncm{\jth}{j^{\scr{th}}}
\ncm{\kth}{k^{\scr{th}}}
\ncm{\lth}{l^{\scr{th}}}
\ncm{\Bin}{\mbox{Bin}}
\ncm{\Exp}{\mbox{Exp}}
\ncm{\Hyp}{\mbox{Hyp}}
\ncm{\mm}{\mbox{mm}}
\ncm{\base}{\scr{bl}}
\ncm{\PD}{\mbox{PD}}
\ncm{\sgn}{\mbox{sgn}}
\ncm{\Ctot}{\bar{C}}
\ncm{\Ctottiny}{C_{\mbox{\tiny tot}}}
\ncm{\bzero}{\bfm{0}}
\ncm{\fappr}{\hat{f}}
\ncm{\bappr}{\hat{b}}
\ncm{\laappr}{\hat{\la}}
\ncm{\muappr}{\hat{\mu}}
\ncm{\pappr}{\hat{p}}
\ncm{\piappr}{\hat{\pi}}
\ncm{\kaappr}{\hat{\ka}}
\ncm{\Siappr}{\hat{\Si}}
\ncm{\bSiappr}{\hat{\bSi}}
\ncm{\demax}{\de_{\scr{max}}}
\ncm{\mumin}{\mu_{\scr{min}}}
\ncm{\hmumin}{\hmu_{\scr{min}}}
\ncm{\Ias}{I_{\scr{as}}}
\ncm{\Ibas}{I_B}
\ncm{\cBall}{\cB_{\scr{all}}}
\ncm{\Inonas}{I_{\scr{nas}}}
\ncm{\Ilong}{I_{\scr{long}}}
\ncm{\Ishort}{I_{\scr{short}}}

\ncm{\beq}{\begin{equation}}
\ncm{\eeq}{\end{equation}}
\ncm{\beqr}{\begin{eqnarray}}
\ncm{\eeqr}{\end{eqnarray}}
\ncm{\beqrn}{\begin{eqnarray*}}
\ncm{\eeqrn}{\end{eqnarray*}}
\ncm\rthm[1]{\ref{#1}}
\ncm\lb[1]{\label{#1}}
\ncm\re[1]{(\ref{#1})}
\ncm{\slut}{
  {\unskip\nobreak\hfill\penalty100\hskip1em\vadjust{}\nobreak
  \hfill\mbox{$\Box$}\parfillskip=0pt\finalhyphendemerits=0}}

\parindent=0mm
\newcommand*\samethanks[1][\value{footnote}]{\footnotemark[#1]}

\newcommand{\MAP}{\hat{\theta}^{\text{(MAP)}}}
\newcommand{\Bayes}{\hat{\theta}^{\text{(Bayes)}}}
\ncm{\cDnew}{{\cal D}^{\text{new}}}
\ncm{\cDobs}{{\cal D}^{\text{obs}}}
\ncm{\cDsim}{{\cal D}^{\text{sim}}}
\ncm{\Ipr}{I^\prime}
\ncm{\hatom}{\hat{\omega}}
\ncm{\Prob}{\mathbb{P}}
\ncm{\E}{\mathbb{E}}
\ncm{\I}{\mathbbm{1}}
\newcommand{\mY}{\mathbf{Y}}
\newcommand{\mX}{\mathbf{X}}
\newcommand{\mB}{\mathbf{B}}
\newcommand{\mE}{\mathbf{E}}
\newcommand{\mZ}{\mathbf{Z}}
\newcommand{\mU}{\mathbf{U}}
\newcommand{\mC}{\mathbf{C}}
\newcommand{\mD}{\mathbf{D}}
\newcommand{\mV}{\mathbf{V}}
\newcommand{\mS}{\mathbf{S}}
\newcommand{\mI}{\mathbf{I}}
\newcommand{\bbeta}{\boldsymbol{\beta}}
\newcommand{\Sigmajuv}{\bSigma^\texttt{juv}}
\newcommand{\Sigmaad}{\bSigma^\texttt{ad}}
\newcommand{\sfI}{\mathsf{I}}
\newcommand{\classI}{\hat{\mathrm{I}}}
\newcommand{\rmI}{\mathrm{I}}
\newcommand{\sfN}{\mathsf{N}}
\newcommand{\sfS}{\mathsf{S}}
\newcommand{\sfK}{\mathsf{K}}
\newcommand{\Rtot}{R^\text{tot}}
\newcommand{\rmd}{\,\mathrm{d}}
\newcommand{\bSigma}{\mathbf{\Sigma}}
\newcommand{\wl}{\textit{wing length }}
\newcommand{\nl}{\textit{notch length }}
\newcommand{\Po}{\text{Po}}
\newcommand{\Geo}{\text{Geo}}
\newcommand{\Be}{\text{Be}}
\newcommand{\Wei}{\text{Wei}}
\newcommand{\U}{\text{U}}
\newcommand{\D}{\mathbf{D}}
\newcommand{\B}{\mathbf{B}}
\newcommand{\J}{\mathbf{J}}
\newcommand{\bt}{\mathbf{t}}
\newcommand{\Sigmat}{\mathbf{\Sigma}}
\newcommand{\Lambdat}{\mathbf{\Lambda}}
\newcommand{\rb}{\right\}}
\newcommand{\lh}{\left[}
\newcommand{\rh}{\right]}
\newcommand{\lp}{\left(}
\newcommand{\rp}{\right)}
\newcommand{\Laplace}{\text{L}}
\newcommand{\convp}{\overset{p}{\longrightarrow}}
\newcommand{\convr}{\overset{r}{\longrightarrow}}
\newcommand{\convd}{\overset{d}{\longrightarrow}}
\newcommand{\tntoinf}{\text{när $n\to\infty$}}
\newcommand{\ntoinf}{$n\to\infty$}
\newcommand{\sumiton}{\sum_{i=1}^n}
\newcommand{\Tloc}{T_{\text{loc}}}
\newcommand{\Tscale}{T_{\text{scale}}}
\newcommand{\Tskew}{T_{\text{skew}}}
\newcommand{\Tkurt}{T_{\text{kurt}}}
\newcommand{\logit}{\text{logit}}

\title{On the Use of $L$-functionals in Regression Models}

\author{Ola H\"{o}ssjer\thanks{Department of Mathematics, Stockholm University, 106 91 Stockholm, Sweden}
\\
M{\aa}ns Karlsson\thanks{Department of Mathematics, Stockholm University, 106 91 Stockholm, Sweden}}
\date{}
\maketitle

\begin{abstract}
In this paper we survey and unify a large class or $L$-functionals of the conditional distribution of the response variable in regression models. This includes robust measures of location, scale, skewness, and heavytailedness of the response, conditionally on covariates. We generalize the concepts of $L$-moments (Sittinen, 1969), $L$-skewness, and $L$-kurtosis (Hosking, 1990) and introduce order numbers for a large class of $L$-functionals through orthogonal series expansions of quantile functions. In particular, we motivate why location, scale, skewness, and heavytailedness have order numbers 1, 2, (3,2), and (4,2) respectively and describe how a family of $L$-functionals, with different order numbers, is constructed from Legendre, Hermite, Laguerre or other types of polynomials. Our framework is applied to models where the relationship between quantiles of the response and the covariates follow a transformed linear model, with a link function that determines the appropriate class of $L$-functionals. In this setting, the distribution of the response is treated parametrically or nonparametrically, and the response variable is either censored/truncated or not. We also provide a unified asymptotic theory of estimates of $L$-functionals, and illustrate our approach by analyzing the arrival time distribution of migrating birds. In this context a novel version of the coefficient of determination is introduced, which makes use of the abovementioned orthogonal series expansion.     
\end{abstract}


\section{Introduction}

Linear combinations of order statistics represent a wide class of estimators of location and scale parameters for samples of independent and identically distributed (i.i.d.) observations. These estimators are well known to combine high efficiency and robustness (Bickel and Lehmann, 1975). Their asymptotic properties are conveniently represented in terms of $L$-functionals of the empirical distribution formed by the sample, as summarized in Chapter 8 of Serfling (1980). 
\par\medskip
Many authors have proposed extensions of $L$-statistics for regression models. In a pioneering article Koenker and Bassett (1978) introduced regression quantiles. These nonparametric estimators make few assumptions on the conditional distribution of the response variable. They naturally extend order statistics to linear regression models, and have been applied to many fields of science (Koenker and Hallock, 2001). Koenker and Portnoy (1989) introduced linear combinations of regression quantiles and showed that much of the computational and asymptotic theory of $L$-functionals of order statistics for i.i.d.\ models extend to the regression framework. The monograph of Koenker (2005) summarizes these and a number of other aspects of regression quantiles, including nonlinear regression. Other contributions include parametric versions of regression quantiles (Gilchrist, 2000, 2007, Frumento and Bottai, 2016), and the use of $L$-functionals in survival analysis (Frumento and Bottai, 2017). 
\par\medskip
In this paper we review, unify, and extend the use of $L$-functionals for regression models. We consider four classes of $L$-functionals or ratios of $L$-functionals that represent measures of location, scale, skewness, and heavytailedness (kurtosis). Then we introduce order numbers for a large class of $L$-functionals through orthogonal series expansions of quantile functions (Takemura, 1983, Okagbue et al., 2019), and motivate why location, scale, skewness, and heavytailedness functionals have order numbers 1, 2, (3,2), and (4,2) respectively. In this context we describe how a given reference distribution gives rise to a whole collection of $L$-functionals. If Legendre polynomials and a uniform reference distribution is used, the resulting class of $L$-functionals of order $1,2,\ldots$ corresponds to $L$-moments (Sillito, 1969), whereas the ratios of $L$-functionals of order (2,1), (3,2), and (4,2) agree with the $L$-coefficient of variation, the $L$-skewness, and the $L$-kurtosis (Hosking, 1990, 1992, 2006) up to normalizing constants. This Legendre class of $L$-functionals is a natural choice for distributions with bounded support, but Hermite polynomials (with a Gaussian reference distribution) or Laguerre polynomials (with an exponential reference distribution) might be preferable for data whose support is on the real line and on the positive real line respectively, in particular if the distribution of the response is close to the reference distribution. 
\par\medskip
We will apply the framework of $L$-functionals to a wide range of transformed linear regression models, with linear models a special case. We argue that the transformation (or link function) will determine the appropriate type of $L$-functionals to use. We also demonstrate how parametric and nonparametric methods of estimating conditional $L$-functionals can be put into a unified framework, for response variables with or without censoring and truncation.  
\par\medskip
In more detail, the paper is organized as follows. In Section \ref{sec:Lstat} we introduce $L$-functionals for models without covariates, with particular emphasis on functionals that quantify location, scale, skewness or heavytailedness, and how their order numbers can be assessed. Then in Section \ref{Sec:Regr} we generalize the framework of Section \ref{sec:Lstat} to linear and transformed linear regression models, for models with our without censoring or truncation. Section \ref{Sec:Sim} contains numerical examples, in Section \ref{Sec:RealData} we analyze a data set with migration times of birds and introduce a novel version of the coefficient of determination, whereas Section \ref{Sec:Disc} provides a summarizing discussion.               

\section{$L$-functionals without covariates}\label{sec:Lstat}

\subsection{Definition of $L$-functionals}

Let $F_Y(y)=F(y)=P(Y\le y)$ be the unknown distribution function of a random variable $Y$. Suppose we want to infer a certain functional $\theta=T(F)$ of $F$, using a data set $Y_1,\ldots,Y_n$ of independent and identically distributed random variables with $P(Y_i\le y)=F(y)$, for $i=1,\ldots,n$. In this paper we will focus on $L$-functionals, i.e.\ linear combinations 
\beq
\theta = T(F) = \int_0^1 Q(p)dG(p)
\lb{T}
\eeq
of quantiles   
\beq
Q(p) = Q_Y(p) = F^{-1}(p)=\inf\{y;\, F(y)\ge p\} 
\lb{Qp}
\eeq
of $F$, using some weight function $G=G^+-G^-$ that corresponds to a signed measure. When $\theta$ is a measure of location of $F$, the weight function is often a positive measure ($G=G^+$), but this is not the case for measures of scale, skewness and kurtosis (cf. Section 2.3). Following Serfling (1980) we consider measures 
\beq
dG(p) = g(p)dp + \sum_{m=1}^M g_m\de_{\pi_m}(p)
\lb{GDef}
\eeq
on $[0,1]$ that split into one absolutely continuous part $g(p)dp$ and another finite sum of point masses $g_m\de_{\pi_m}$ at $\pi_m$ with weights $g_m$ for $m=1,\ldots,M$, with $0\le \pi_1 < \ldots \pi_M \le 1$. This allows us to work with all the common types of distributions, i.e. continuous, discrete and mixtures, within the same framework. The $L$-functional is robust if extreme quantiles of $F$ are excluded, i.e.\ if the total variation measure $|G|=G^{+}+G^{-}$ satisfies $|G|([0,\pi)\cup (1-\pi,1])=0$ for some sufficiently small $0<\pi\le 0.5$. Then the breakdown point (Rousseeuw and Leroy, 1987) of $T$ is at least $\pi$. 
\par\medskip
A wide class of estimators 
\beq
\hth = T(\hF) = \int_0^1 \hQ(p)dG(p),
\lb{hthLGen}
\eeq
of $\theta$ are obtained by plugging $\hF$, an estimate of $F$, into \re{T}, with $\hQ=\hF^{-1}$ the corresponding estimate of $Q$. This estimate of $F$ could be nonparametric, i.e.\ the empirical distribution function
\beq
\hF(y) = \frac{1}{n}\sum_{i=1}^n \mathbbm{1}(Y_i\le y)
\lb{edf}
\eeq
formed by the sample. Then
\beq
\hth = \sum_{i=1}^n w_i Y_{(i)}
\lb{hthL}
\eeq
is an $L$-statistic, i.e.\ a linear combination of the order statistics $Y_{(1)}\le \ldots \le Y_{(n)}$, with $w_i = \int_{(i-1)/n}^{i/n} dG(p)$ the weight assigned to the $i$th order statistic. It is also possible to insert a parametric estimator 
\beq
\hF(y) = F(y;\hbpsi)
\lb{hFpar}
\eeq
of $F$ into \re{hthLGen}. In this case $F(y)=F(y;\bpsi)$ is fully determined by a finite-dimensional parameter $\psi=(\psi_1,\ldots,\psi_r)$, of which $\hbpsi=(\hpsi_1,\ldots,\hpsi_r)$ is an estimate. If the density function $F^\prime(y;\bpsi)=f(y;\bpsi)$ is tractable, this is typically the maximum likelihood estimator of $\bpsi$. For some distributional families, such as the Generalized Lambda Distributions (Karian and Dudewicz, 2000) or mixtures of quantile functions (Karvanen, 2006, Karvanen and Nuutinen, 2008), the quantile function $Q(p;\bpsi)=F^{-1}(p;\bpsi)$ has a more explicit form. It might then be more tractable to estimate $\bpsi$ by fitting some of the order statistics $Y_{(i)}$ to $Q$ (Gilchrist, 2007) or some empirical $L$-functionals to the corresponding population-based $L$-functionals (Karvanen, 2006). 

\subsection{Asymptotics}

In order to study the large sample behavior of $\hth=\hth_n$ as $n\to\infty$, we introduce $S$, the space of real-valued and integrable functions on $(0,1)$, equipped with a seminorm
\beq
\|Q\|_S = \int_0^1 |Q(p)| |g(p)| dp + \sum_{m=1}^M |g_m| |Q(\pi_m)|.
\lb{QS}
\eeq
where $g(p)$ and $g_1,\ldots,g_M$ refer to the absolutely continuous part and the point masses of the weight measure $G$ in \eqref{GDef} respectively. We may also regard \eqref{QS} as a norm of the restriction of $Q$ to $\mbox{supp}(G)$. If $F_1$ and $F_2$ are two distributions with quantile functions $Q_1$ and $Q_2$, it follows that   
\beq
\begin{split}
|T(F_2)-T(F_1)| &\le  \int_0^1 |Q_2(p)-Q_1(p)|d|G|(p) \\
&= \int_0^1 |Q_2(p)-Q_1(p)||g(p)|dp + \sum_{m=1}^M |g_m||Q_2(\pi_m)-Q_1(\pi_m)|\\
&= \|Q_2-Q_1\|_S.  
\end{split}
\lb{TDiff}
\eeq
This implies in particular that $T(F)$ is a continuous functional with respect to the distance measure $d(F_2,F_1)=\|Q_2-Q_1\|_S$ introduced by $\|\cdot\|_S$. In particular, when $M=0$ and $g(p) \equiv 1$, this distance measure equals the Wasserstein metric of order 1 (Olkin and Pukelsheim, 1982). 
\par\medskip
The asymptotic properties of $\hth_n$ are determined by the weight function $G$ and the large sample behavior of the rescaled quantile process 
\beq
Z_n(p) = \sqrt{n}[\hQ_n(p)-Q(p)], \quad 0<p<1,
\lb{Znp}
\eeq
which is a random element of $S$. We will assume weak convergence (Billingsley, 1999)
\beq
Z_n \Lto Z\mbox{ as }n\to\infty,
\lb{ZnLtoZ}
\eeq
with respect to the topology introduced by the norm in $S$, i.e.\ $E[h(Z_n)]\to E[h(Z)]$ for all bounded and continuous functions $h:S\to\R$. It is further assumed that the limit in \re{ZnLtoZ} is a Gaussian process with mean function $E[Z(p)]=0$ and covariance function $\Cov[Z(p),Z(s)]=R(p,s)$ for $p,s \in \mbox{supp}(G)$. Notice in particular that for discrete weight measures ($g\equiv 0$), only the last term on the right hand side of \re{QS} is present. Then \re{ZnLtoZ} corresponds to weak convergence of finite-dimensional distributions of $Z_n$ towards $Z$ at $\pi_1,\ldots,\pi_m$. This follows from the fact that $Z_n$ and $Z$, restricted to $\mbox{supp}(G) =\{\pi_1,\ldots,\pi_M\}$, represent weak convergence of random vectors of dimension $M$ when $g\equiv 0$. On the other hand, \re{ZnLtoZ} represents functional weak convergence on an infinite index set $p\in\mbox{supp}(G)$ when $g\ne 0$, due to the first term of \re{QS}.
\par\medskip
It is also possible to define 
\beq
\|Q\|_S=\|Q\|_\infty = \sup_{0<p<1} |Q(p)|
\label{QSup}
\eeq
as the supremum of $Q$ on $(0,1)$. It is easy to see that $T(F)$ is continuous with this choice of norm, since  $|T(F_2)-T(F_1)|\le |G|(0,1) \|Q_2-Q_1\|_S$. Note that $\|Q\|_S$ in \re{QSup} does not involve the weight measure $G$. This is advantageous when simultaneous weak convergence of several $L$-functionals, based on different weight functions, is of interest. On the other hand, the advantage of (8) is that this (semi)norm exists for a larger class of functions $Q$. It also gives rise to a weaker topology on $S$, so that (11) requires less. Indeed, when (8) is used, weak convergence of $Z_n$ need only be established on the index set $\mbox{supp}(G)$.
\par\medskip
Since \re{TDiff} implies that $\tilde{T}:S\to\R$, defined by $\tilde{T}(Q)=\int_0^1 Q(p)dG(p)$, is a continuous functional, it follows from the Continuous Mapping Theorem that 
\beq
\sqrt{n}(\hth_n-\theta) = \tilde{T}(Z_n) \Lto \tilde{T}(Z) \sim N(0,\Sigma)
\lb{AsN}
\eeq  
as $n\to\infty$, where the asymptotic covariance matrix of the limiting normal distribution satisfies
\beq 
\Sigma = \int_0^1 \int_0^1 R(p,s) dG(p)dG(s). 
\lb{Si}
\eeq
The expression for the covariance function $R$ depends on which estimator of $F$ that is used in \re{hthLGen}. The nonparametric estimate \re{edf} corresponds to
\beq
R(p,s) = \frac{\min(p,s)-ps}{f(Q(p))f(Q(s))},
\lb{Rpsedf}
\eeq
where $f(y)=F^\prime(y)$ is the density function of $Y$. When \re{Rpsedf} is inserted into \re{Si}, we get a well known expression for the asymptotic variance of $L$-statistics, cf.\ Mosteller (1946), Bennett (1952), Jung (1955), Chernoff et al.\ (1967), and Moore (1968). In the parametric case we have that 
\beq
R(p,s) = \frac{dQ(p;\bpsi)}{d\bpsi} \bV \left(\frac{dQ(s;\bpsi)}{d\bpsi}\right)^T,
\lb{Rpar}
\eeq
provided $\hbpsi=\hbpsi_n$ is asymptotically normal with covariance matrix $\bV$, i.e. 
\beq
\sqrt{n}(\hbpsi_n-\bpsi) \Lto N(\bzero,\bV)
\lb{psiConv}
\eeq
as $n\to\infty$. In particular, $\bV$ equals the inverse of the Fisher information matrix of $\bpsi$, when $\hbpsi_n$ is the ML-estimator of $\bpsi$. 

\subsection{Examples of $L$-functionals}\lb{Sec:ExLFunc}

The choice of weight function $G$ in \re{hthLGen} will determine the type of $L$-functional. In the present paper we will mainly focus on four types of statistical functionals, presented in \ Oja (1981), namely for location, scale, skewness and kurtosis:
\par\medskip
\begin{exa}[Measures of location]
A measure of location ($T=T_{\scr{loc}}$) is equivariant with respect to linear transformations of data, i.e.\
\beq
T_{\scr{loc}}(F_{aY+b}) = aT_{\scr{loc}}(F_Y) + b
\lb{Tloc}
\eeq
for any real valued $a$ and $b$. Since $Q_{aY+b}=aQ_Y+b$, it follows that \re{Tloc} holds for all $L$-functionals with a weight function satisfying 
\beq
\int_0^1 dG(p) = 1.
\lb{GLoc}
\eeq
For location functionals we will also require that $G$ is a positive measure, so that 
\beq
G^- = 0.
\lb{Gloc-}
\eeq
The two regularity conditions \re{GLoc}-\re{Gloc-} imply that $T_{\scr{loc}}$ preserves stochastic ordering of distribution functions. By this we mean that $T_{\scr{loc}}(F_1)\le T_{\scr{loc}}(F_2)$ whenever $F_2$ is stochastically larger than $F_1$, that is, when the quantile functions of the two distributions satisfy $Q_1(p)\le Q_2(p)$ for all $0<p<1$. 
\par\medskip
Examples of location functionals with a discrete weight function ($g\equiv 0$) include quantiles 
\beq
T_{\scr{loc}}(F) = Q(\pi)
\lb{GRQ}
\eeq
for some fixed $0<\pi < 1$ and for distributions $F$ with bounded support the midrange
\beq
T_{\scr{loc}}(F) = \frac{1}{2}\left[Q(0)+Q(1)\right].
\lb{GMidRange}
\eeq
Examples of location functionals with an absolutely continuous weight density include smoothed quantiles
\beq
T_{\scr{loc}}(F) = \int_0^1 \frac{1}{b}K\left(\frac{p-\pi}{b}\right)Q(p)dp,
\lb{GKernel}
\eeq
where $K$ is the smoothing probability density and $0<b<\min(\pi,1-\pi)$ the bandwidth (Parzen, 1979, Sheather and Marron, 1990), and compound expectations (CEs) 
\beq
T_{\scr{loc}}(F) = \frac{1}{\pi_1-\pi_0} \int_{\pi_0}^{\pi_1} Q^{-1}(p)dp
\lb{GCCE}
\eeq
that average quantiles of $F$ between $\pi_0$ and $\pi_1$ for some appropriately chosen $0\le\pi_0<\pi_1\le 1$. Notice in particular that \re{GCCE} is a special case of \re{GKernel} that corresponds to a rectangular kernel $K(x)=1(|x|\le 1)/2$, $\pi=(\pi_1+\pi_2)/2$ and $b=(\pi_1-\pi_0)/2$.  
\par\medskip
A centralized measure of location is one whose weight function is symmetric around $p=0.5$ and satisfies \re{GLoc}. If $F_{Y}$ has a symmetric distribution around $\mu$, a centralized measure of location will equal the center of symmetry ($T_{\scr{loc}}(F_Y)=\mu$). These types of weight functions include the median ($\pi=0.5$ in \re{GRQ}), the midrange \re{GMidRange}, and the trimmed mean ($\pi_0+\pi_1=1$ in \re{GCCE}, see for instance Tukey and McLaughlin (1963), Bickel (1965), and Stigler (1977)). When $\pi_0=0$ and $\pi_1=1$, the trimmed mean simplifies to the expected value 
\beq
T_{\scr{loc}}(F_Y) = \int_0^1 Q_Y(p)dp = E(Y).
\lb{TEY}
\eeq
A location-scale family 
\beq
F(y)=F_0\left(\frac{y-\mu}{\si}\right)
\lb{LocScale}
\eeq
corresponds to a parametric family \re{hFpar} with $\bpsi=(\mu,\si)$ and $F_0$ known. It is then of interest to find the
optimal weight function $G$ of location for the $L$-statistic \re{hthL}. This weight function minimizes the asymptotic variance \re{Si}-\re{Rpsedf} among all $L$-functionals $T_{\scr{loc}}$ that estimate $\mu$ in \re{LocScale}. It is well known (Chernoff et al., 1967, Chapter 8 of Serfling, 1980) that if the density $f_0=F_0^\prime$ is twice differentiable, the asymptotically optimal $L$-estimator of $\mu$ (when $\si$ is known or when $f_0$ is symmetric so that $\sigma$ is Fisher orthogonal to $\mu$) corresponds to an  absolutely continuous weight measure with density
\beq
g(p) = - \frac{(\log f_0)^{\prime\prime}(F_0^{-1}(p))}{I(f_0)},
\lb{gLocOpt}
\eeq
where $I(f_0)=\int f_0^\prime(y)^2/f_0(y)dy$ is the Fisher information. In particular, the expected value \re{TEY} with weight function $g(p)\equiv 1$ is optimal when $F\sim N(\mu,\si^2)$, whereas the midrange \re{GMidRange} is optimal for a uniform distribution, as can be seen by approximating this distribution by a smooth $F_0$.   
\slut\end{exa}
\par\medskip
\begin{exa}[Measures of scale]
A measure of scale ($T=T_{\scr{scale}}$) is non-negative and satisfies
\beq
T_{\scr{scale}}(F_{aY+b}) = |a|T_{\scr{scale}}(F_Y)
\lb{Tscale}
\eeq
for all real-valued $a$ and $b$. In the context of $L$-functionals, \re{Tscale} is satisfied for weight functions that take on positive as well as negative values, in such a way that the conditions
\beq
dG\mbox{ is skew-symmetric around }p = 0.5,
\lb{GScale}
\eeq
and
\beq
\begin{array}{rcl}
G^+ &=& G\mbox{ restricted to }(0.5,1],\\
G^- &=& G\mbox{ restricted to }[0,0.5)
\end{array}
\lb{GScaleb}
\eeq
are fulfilled. In order to see that \re{GScale}-\re{GScaleb} correspond to a measure of scale, it is instructive to insert these conditions into \re{T}. This makes it possible to rewrite the scale functional as  
\beq
\begin{array}{rcl}
T_{\scr{scale}}(F) &=& \int_0^1 (Q(p)-Q(0.5))dG(p)\\
&=& \int_0^1 |Q(p)-Q(0.5)|d|G|(p)\\
&=& \int_{0.5}^1 \left[Q(p)-Q(1-p)\right]dG^+(p).
\end{array}
\lb{TscaleIQRange}
\eeq
Notice in particular that \re{TscaleIQRange} has an intuitive interpretation as a linear combination of interquantile ranges $Q(p)-Q(1-p)$.  
Bickel and Lehmann (1976) and Oja (1981) defined a spread-ordering among distributions, where $F_2$ is said to be at least as spread out as $F_1$ if the difference $Q_2(p)-Q_1(p)$ between the quantile functions of $F_2$ and $F_1$ is non-decreasing. It follows from \re{TscaleIQRange} that scale functionals preserve spread ordering, i.e.\ $T_{\scr{scale}}(F_1)\le T_{\scr{scale}}(F_2)$. Since the unit of $T_{\scr{scale}}$ is somewhat arbitrary, some reference distribution $F_0$ is typically chosen to have scale 1, i.e.  
\beq
T_{\scr{scale}}(F_0) = \int_0^1 Q_0(p)dG(p) = 1.
\lb{GScaleN}
\eeq
If this reference distribution is a standard normal $N(0,1)$, it follows that $T_{\scr{scale}}(F_Y)=\si$ for $Y\sim N(\mu,\si^2)$.
The simplest scale functional that satisfies \re{GScale}-\re{GScaleb} and \re{GScaleN}, is the one for which $G^+$ is a point measure at $\pi$ for some $0.5 < \pi \le 1$. Then 
\beq
T_{\scr{scale}}(F) = \frac{Q(\pi)-Q(1-\pi)}{Q_0(\pi)-Q_0(1-\pi)} 
\lb{QuantRange}
\eeq
equals the standardized interquartile range of $F$ when $\pi=0.75$ and the standardized range when $\pi=1$ and $F$ has bounded support. The standardized Gini's mean difference
\beq
T_{\scr{scale}}(F) = K E|Y_1-Y_2| = \frac{\int_{0.5}^1 \left[Q(p)-Q(1-p)\right](p-0.5)dp}{\int_{0.5}^1 \left[Q_0(p)-Q_0(1-p)\right](p-0.5)dp} 
\lb{TGini}
\eeq
corresponds to having $dG^+(p)=K(p-0.5)1(0.5\le p \le 1)dp$ in \re{TscaleIQRange}, with $K$ chosen so that \re{GScaleN} holds. A third measure of scale 
\beq
T_{\scr{scale}}(F) = \frac{\int_{0.5}^1 \left[Q(p)-Q(1-p)\right]dp}{\int_{0.5}^1 \left[Q_0(p)-Q_0(1-p)\right]dp} 
\lb{TscalegConst}
\eeq
has a constant density $dG^+(p)=K 1(0.5\le p \le 1)dp$ on $(0.5,1)$, with $K>0$ a constant chosen so that \re{GScaleN} holds. 
\par\medskip
It is sometimes of interest to choose the weight function $G$ so that the asymptotic variance $\Si$ in \re{Si} is minimized. It can be shown (Chernoff et al., 1967, Chapter 8 of Serfling, 1980) that for the location-scale family \re{LocScale}, an $L$-functional with absolutely continuous weight function 
\beq
g(p) = K\left[(\log f_0)^\prime\left(F_0^{-1}(p)\right) + F_0^{-1}(p)(\log f_0)^{\prime\prime}\left(F_0^{-1}(p)\right)\right]
\lb{gOptScale}
\eeq
corresponds to an asymptotically optimal estimator of the scale parameter $\si$ whenever the location parameter $\mu$ is known or orthogonal to $\si$, and the density function $f_0$ is twice differentiable, with $K$ chosen so that \re{GScaleN} holds for some appropriately chosen reference distribution. In particular, $K=1$ yields a consistent estimator of the scale parameter $\si$ in \re{LocScale}. For instance, $g(p)=F_0^{-1}(p)$ is optimal for the normal distribution ($F_0\sim N(0,1)$). A robustified, skew-symmetric and trimmed version of this weight function has been studied by Welsh and Morrison (1990). More generally, it follows from \re{gOptScale} that the optimal $g$ is skew-symmetric around 0.5 whenever $f_0$ is symmetric around 0, but not when symmetry of $f_0$ fails. The advantage of having a skew-symmetry requirement \re{GScale} on the weight function is the intuitive interpretation \re{TscaleIQRange} of $T_{\scr{scale}}(F)$ in terms of a linear combination of interquantile ranges. Koenker and Zhou (1994) define scale functionals more generally by dropping the skew-symmetry condition and only requiring that the weight function satisfies $G((0,p])<0$ for all $0<p<1$ and $G([0,1])=0$. An example of such a scale functional is presented in Section \ref{Sec:LOrder}.  
\slut\end{exa}
\par\medskip
\begin{exa}[Measures of skewness] The traditional measure 
\beq
T_{\scr{skew}}(F_Y) = \frac{E\left[(Y-E(Y))^3\right]}
{\left\{E\left[(Y-E(Y))^2\right]\right\}^{3/2}}
\lb{TskewTrad}
\eeq 
of skewness compares the left and right tails of $F_Y$. It is not robust, since the third moment of $Y$ must be finite. Here we will analyze versions of skewness that are more robust than \re{TskewTrad}, and defined as the ratio 
\beq
T_{\scr{skew}}(F)=\frac{T_{\scr{uskew}}(F)}{T_{\scr{scale}}(F)}
\lb{TskewDef}
\eeq
of two $L$-functionals. The numerator of \re{TskewDef} corresponds to an unstandardized measure of skewness, which transforms as
\beq
T_{\scr{uskew}}(F_{aY+b}) = aT_{\scr{uskew}}(F_Y)
\lb{Tskew}
\eeq
under linear mappings, for all real-valued $a$ and $b$. The denominator of \re{TskewDef} is another $L$-functional that measures scale. It follows from \re{Tscale} and \re{Tskew} that the standardized skewness satisfies
\beq
T_{\scr{skew}}(F_{aY+b}) = \sgn(a)T_{\scr{skew}}(F_Y)
\lb{Tstskew}
\eeq
for all non-negative $a$, with $\sgn(a)=a/|a|$. We will assume here that $T_{\scr{skew}}$ is chosen so that
\beq
T_{\scr{skew}}(F_0) = 1
\lb{TstskewF0}
\eeq
holds for some reference distribution $F_0$ whose right tail is heavier than the left tail (e.g.\ an exponential distribution).   
\par\medskip
A class of weight functions whose unstandardized skewness functional transform linearly, as in \re{Tskew}, are those that satisfy 
\beq
\begin{array}{rcl}
\int_0^1 dG(p) &=& 0,\\
dG\mbox{ symmetric around }p &=& 0.5,
\end{array}
\lb{GSkew}
\eeq
and
\beq
\begin{array}{rcl}
\mbox{supp}(G^-) &\subseteq & [\pi,1-\pi],\\
\mbox{supp}(G^+) &\subseteq & [0,\pi]\cup [1-\pi,1],
\end{array}
\lb{GSkewb}
\eeq
for some $0\le \pi\le 0.5$. Indeed, inserting \re{GSkew}-\re{GSkewb} into \re{T}, we find that 
\beq
\begin{array}{rcl}
T_{\scr{uskew}}(F) &=& \int_{0.5}^1 \left[Q(p)+Q(1-p)-2Q(0.5)\right]dG(p)\\
&=& \int_{1-\pi}^1 \left[Q(p)+Q(1-p)-2Q(0.5)\right]dG^+(p)\\
&-& \int_{0.5}^{1-\pi} \left[Q(p)+Q(1-p)-2Q(0.5)\right]dG^-(p).
\end{array}
\lb{Tskew2}
\eeq
In particular, if $F_{Y}$ is symmetric around its center of symmetry $\mu$, if follows that $T_{\scr{uskew}}(F)=0$. Since $Q(p)+Q(1-p)-2Q(0.5)$ measures skewness for each quantile $0.5<p<1$, \re{Tskew2} is a functional that compares skewness of the tails of $F$ with the skewness of the central part of $F$. Of particular interest is the case when $G^-$ is a point measure at 0.5. Then \re{Tskew2} simplifies to 
\beq
T_{\scr{uskew}}(F) = \int_{0.5}^1 \left[Q(p)+Q(1-p)-2Q(0.5)\right]dG^+(p).
\lb{Tskew3}
\eeq 
If the positive part of the weight function is chosen as $dG^+(p) = K\de_{1-\pi}(p)+K\de_\pi(p)$ for some $0<\pi < 0.5$ and $K>0$, and
if the interquantile range in \re{QuantRange} is used for standardization, one obtains the standardized measure 
\beq
T_{\scr{skew}}(F) = K\cdot \frac{Q(1-\pi)+Q(\pi)-2Q(0.5)}{Q(1-\pi)-Q(\pi)}
\lb{GStandSkew}
\eeq
of skewness. It was introduced by Galton (1883) and Bowley (1920) for $\pi=0.25$ and $K=1$, and for arbitrary $\pi$ by Hinkley (1975). Here we will rather define $K$ so that \re{TstskewF0} holds for some appropriate reference distribution $F_0$. An alternative to \re{GStandSkew} is to choose an unstandardized skewness functional \re{Tskew3} for which $G^+$ has a constant density in the numerator of \re{TskewDef}, and then use a scale measure with a constant weight density $g^+$ for quantiles above 0.5, in the denominator of \re{TskewDef}. The corresponding standardized measure 
\beq
T_{\scr{skew}}(F) = K\cdot \frac{\int_{0.5}^1 \left[Q(p)+Q(1-p)-2Q(0.5)\right]dp}{\int_{0.5}^1 \left[Q(1-p)-Q(p)\right]dp}
\lb{GStandSkew2}
\eeq
of skewness with $K=1$ was proposed by Groeneveld and Meeden (1984). Here we will rather choose $K$ in order for \re{TstskewF0} to hold. Notice that the skewness measure \re{GStandSkew2} puts higher weights on the tails of the distribution, compared to \re{GStandSkew}. See also Kim and White (1994) for an overview of different robust measures of skewness. 
\par\medskip
Groeneveld and Meeden (1984) argued that a reasonable skewness measure should satisfy \re{Tstskew} and in addition preserve the partial skewness-ordering\break among distributions, due to van Zwet (1964). By this we mean that $T_{\scr{skew}}(F_1)\le T_{\scr{skew}}(F_2)$ whenever $F_2$ is at least as skewed to the right as $F_1$, i.e.\ if $x\to F_2^{-1}(F_1(x))$ is a convex function. Groeneveld and Meeden verified that \re{GStandSkew} and \re{GStandSkew2} preserve this skewness-ordering among distributions. This and other partial skewness-orderings are discussed by Oja (1981), McGillivray (1986), and Garcia et al.\ (2018). 
\slut\end{exa}
\par\medskip
\begin{exa}[Measures of heavytailedness]
There is no universal agreement whether the (excess) kurtosis
\beq
T(F_Y) = \frac{E[(Y-E(Y))^4]}{\left\{E\left[(Y-E(Y))^2\right]\right\}^2}-3
\lb{Tkurtosis}
\eeq
quantifies peakness versus tails of $F_Y$ or modality versus bimodality. Following Chissom (1970) and Oja (1981), we will regard kurtosis as a measure of heavytailedness. In order to find more robust measures of heavytailedness, we will consider functionals 
\beq 
T_{\scr{heavy}}(F) = \frac{T_{\scr{uheavy}}(F)}{T_{\scr{scale}}(F)}
\lb{TheavyDef}
\eeq
that are defined as the ratio of two $L$-functionals. The functional $T_{\scr{uheavy}}$ in the numerator of \re{TheavyDef} corresponds to an unstandardized measure of heavytailedness, and it transforms as  
\beq
T_{\scr{uheavy}}(F_{aY+b}) = |a|T_{\scr{uheavy}}(F_Y)
\lb{Theavy}
\eeq
under linear mappings, for all real-valued $a$ and $b$, whereas the functional in the denominator of \re{TheavyDef} corresponds to a measure of scale. It follows from \re{Tscale} and \re{Theavy} that $T_{\scr{heavy}}$ is invariant with respect to linear transformations, i.e.\ 
\beq
T_{\scr{heavy}}(F_{aY+b}) = T_{\scr{heavy}}(F_Y).
\lb{Tstheavy}
\eeq
Notice that \re{Theavy} is identical to the corresponding relation \re{Tscale} for scale functionals. But whereas $T_{\scr{scale}}$ is always non-negative, we want $T_{\scr{heavy}}$ to be positive for heavy-tailed distributions and negative for light-tailed distributions. In order to accomplish this we choose the weight function in \re{T}, for $T_{\scr{uheavy}}(F)$, as 
\beq
\begin{array}{rcl}
dG\mbox{ skew-symmetric around }p = 0.5,
\end{array}
\lb{GHT}
\eeq
and
\beq
\begin{array}{rcl}
G^+ &=& G\mbox{ restricted to }[\pi,0.5)\cup (1-\pi,1],\\
G^- &=& G\mbox{ restricted to }[0,\pi)\cup (0.5,1-\pi]
\end{array}
\lb{GHTb}
\eeq
for some $0<\pi<0.5$. In addition, we require that the unstandardized and standardized measures of heavytailedness satisfy
\beq
\begin{array}{rcl}
T_{\scr{uheavy}}(F_0) &=& \int_0^1 Q_0(p)dG(p) = 0,\\
T_{\scr{heavy}}(F_1) &=& \int_0^1 Q_1(p)dG(p) / T_{\scr{scale}}(F_1) = 1
\end{array}
\lb{GHTNormal}
\eeq
for two distributions $F_0$ and $F_1$ with quantile functions $Q_0$ and $Q_1$. When there is no restriction on the range of $Y$, $F_0$ is typically a normal distribution, whereas $F_1$ is a symmetric and moderately light tailed distribution, such as the Laplace distribution. 
\par\medskip
In order to motivate that \re{GHT} leads to a measure of heavytailedness, we insert this equation into \re{T} and notice that
\beq
T_{\scr{uheavy}}(F) = \int_{1-\pi}^1 \left[Q(p)-Q(1-p)\right]dG^+(p) - \int_{0.5}^{1-\pi} \left[Q(p)-Q(1-p)\right]dG^-(p)
\lb{Theavy2}
\eeq
equals the difference between the weighted interquantile differences $Q(p)-Q(1-p)$ in the tails and in the central part of the distribution, respectively. 
Several robust measures of kurtosis fit into our framework. As a first example, Moors (1988) introduced
\beq
T_{\scr{heavy}}(F) = K\left[\frac{Q(0.875)-Q(0.625)+Q(0.375)-Q(0.125)}{Q(0.75)-Q(0.25)} - 1.23\right],
\lb{TMoors}
\eeq
with $K=1$, where the unstandardized kurtosis in the numerator has a weight function $G$ such that the restrictions of $G^+$ and $G^-$ to $(0.5,1)$ have a one point distribution at $0.875$ and a two point distribution at $0.625$ and $0.75$ respectively. The scale measure in the denominator of \re{TMoors}, on the other hand, corresponds to an unstandardized interquartile range. Finally, the constant 1.23 is chosen so that the upper part of \re{GHTNormal} holds when $F_0\sim N(0,1)$. Here we will additionally choose $K>0$ in \re{TMoors} so that the lower part of \re{GHTNormal} is satisfied for some appropriately chosen reference distribution $F_1$. Second, the tail ratio of Gilchrist (2000) can be normalized as
\beq
T_{\scr{heavy}}(F) = K_1\left[ \frac{Q(0.9)-Q(0.1)}{Q(0.75)-Q(0.25)} - K_0\right], 
\lb{TGilchrist}
\eeq
where $K_0$ and $K_1$ are chosen so that \re{GHTNormal} holds. A third class of kurtosis measures
\beq
T_{\scr{heavy}}(F) = K_1\left[\frac{\int_{1-\pi_0}^1 \left[Q(p)-Q(1-p)\right]dp}
{\int_{1-\pi_1}^1 \left[Q(p)-Q(1-p)\right]dp} - K_0\right]
\lb{THogg}
\eeq
was introduced by Hogg (1972,1974) for some conveniently selected $0<\pi_0<\pi_1\le 0.5$, and with $K_0$ chosen so that the upper part of \re{GHTNormal} holds when $F_0\sim N(0,1)$.
Hogg conducted simulations for which $\pi_0=0.05$, $\pi_1=0.5$, and $K_0=2.59$ gave satisfactory results. Whereas Hogg used $K_1=1$, we will rather choose $K_1$ so that the lower part of \re{GHTNormal} is satisfied for some appropriate reference distribution $F_1$. See also Kim and White (1994) for an overview of different robust measures of kurtosis.
\par\medskip
A partial kurtosis-ordering of symmetric distributions (van Zwet, 1964, Oja, 1981) states that $F_2$ has more kurtosis than $F_1$ if $x\to F_2^{-1}(F_1(x))$ is a concave (convex) to the left (right) of the point of symmetry. It can be shown that both \re{TMoors} and \re{THogg} preserve this kurtosis-ordering among symmetric distributions.  
\par\smallskip
Sometimes it is only one tail of $F$ that is of interest. It is possible then to split the weight function of an unstandardized heavytailedness functional as 
$$
dG(p) = \mathbbm{1}(p<0.5)dG(p) + \mathbbm{1}(p>0.5)dG(p) =: dG_{\scr{left}}(p) + dG_{\scr{right}}(p),
$$
where only the low quantiles are included in $G_{\scr{left}}$ in order to study the left tail of $F$, whereas only the high quantiles are used in $G_{\scr{right}}$ to study the right tail of $F$. It follows from \re{GHT} that the weight measure for the right tail satisfies 
\beq
\begin{array}{rcl}
G_{\scr{right}}^+ &=& G\mbox{ restricted to }(1-\pi,1],\\
G_{\scr{right}}^- &=& G\mbox{ restricted to }(0.5,1-\pi].
\end{array}
\lb{GHTRight}
\eeq
The corresponding functional 
\beq
T_{\scr{uheavyr}}(F) = \int_{1-\pi}^1 Q(p)G^+(p) - \int_{0.5}^{1-\pi} Q(p)dG^-(p)
\lb{TheavyRight}
\eeq
can be viewed as a restriction of \re{Theavy2} to an interval $(0.5,1)$ of quantiles, and it is standardized as 
$$
T_{\scr{heavyr}}(F) = \frac{T_{\scr{uheavyr}}(F)}{T_{\scr{scale}}(F)}.
$$  
Suppose, for instance, that $e^Y$ has a heavy tail to the right, in the sense that 
$$
P(e^Y \ge y) = y^{-\alpha}L(y)
$$
for some tail parameter $\alpha>0$, and with a function $L(y)$ that varies slowly as $y\to\infty$, i.e.\ $L(ty)/L(y)\to 1$ as $y\to\infty$ for all $t>0$. Then Hill's estimator of $\alpha$ (Hill 1975, Haeusler and Teugels, 1985) corresponds to an estimator $\hth = T_{\scr{uheavyr}}(\hF)$ in \re{hthL} with weight function 
\beq
\begin{array}{rcl}
g_{\scr{right}}^+(p) &=& \mathbbm{1}(1-\pi< p < 1)/(1-\pi),\\
dG_{\scr{right}}^-(p) &=& \de_{1-\pi}(p).
\end{array}
\lb{GHill}
\eeq
In order for $\hth=\hth_n$ to be a consistent estimator of $\alpha$, it is either required that $L(\cdot)$ is constant for large enough $y$, or that $\pi=\pi_n\to 0$ as $n\to\infty$ at a rate depending on how much $L(\cdot)$ varies for large $y$. 
\slut\end{exa}

\subsection{The order of $L$-functionals}\lb{Sec:LOrder}

In this section we introduce order numbers for a certain subclass of $L$-function\-als. Since each $L$-functional is determined by its weight measure $G$, this amounts to introducing order numbers for a subclass of weight measures. We will state two mandatory conditions on $G$ and a third optional symmetry condition, in order for it to be of order $m\in\{1,2,\ldots\}$. First, there has to exist a disjoint decomposition of $[0,1]$ into $m$ intervals $I_1<I_2<\ldots<I_m$ such that
\beq
\begin{array}{rcl}
G^-(I_k)=0\mbox{ and }G^+(I_k)>0, & \mbox{ if $m-k$ is even},\\
G^-(I_k)>0\mbox{ and }G^+(I_k)=0, & \mbox{ if $m-k$ is odd}. 
\end{array}
\lb{GSign}
\eeq
Second,  
\beq
\int_0^1 dG(p)  = \left\{\begin{array}{rcl}
1, & \mbox{if $G$ has order $m=1$},\\
0, & \mbox{if $G$ has order $m>1$}.
\end{array}\right.
\lb{GInt}
\eeq
Third, 
\beq
dG\mbox{ is symmetric (skew-symmetric) around $p=0.5$ if $m$ is odd (even).}
\lb{GSymm}
\eeq
When all three conditions \re{GSign}-\re{GSymm} hold, we refer to $G$ as a symmetric (skew-symmetric) weight measure of order $m$.
\par\medskip 
It may be verified from Section \ref{Sec:ExLFunc} that location functionals have order 1, scale functionals order 2, unstandardized skewness functionals order 3, and unstandardized heavytailedness functionals order 4. Write 
\beq
T_m(F) = \int_0^1 Q(p)dG_m(p)
\lb{TmDef}
\eeq
for an $L$-functional whose weight function $G=G_m$ is of order $m$, and 
\beq
T_{ml}(F) = \frac{T_m(F)}{T_l(F)}
\lb{TmlDef}
\eeq
for a functional that is the ratio of two $L$-functionals of order $m$ and $l$ respectively. We will refer to $(m,l)$ as the order number of $T_{ml}$. Consequently, skewness functionals have order $(3,2)$ and heavytailedness functionals order $(4,2)$. 
\par\medskip
It is possible to obtain collections $\{T_m(F)\}_{m=1}^\infty$ of $L$-functionals from orthogonal series expansions of the quantile function $Q=F^{-1}$.
Each such collection makes use of a reference distribution $F_0$ with density function $f_0(y)=F_0^\prime(y)$ on its support $[a,b]$, where $-\infty\le a < b \le\infty$. Introduce the scalar product $\langle f,g \rangle =\int_a^b f(y)g(y)f_0(y)dy$ for functions that are square integrable with respect to $f_0(y)dy$, and suppose there exists an orthonormal system of polynomials $\{P_k\}_{k=0}^\infty$ of degrees $k=0,1,2,\ldots$ such that
\beq
\int_a^b P_k(y)P_l(y)f_0(y)dy = 1(k=l).
\lb{PnOrth}
\eeq
Then define the absolutely continuous weight densities ($dG_m(p)=g_m(p)dp$) 
\beq
g_m(p) = P_{m-1}(Q_0(p))
\lb{gmGen}
\eeq
for $m=1,2,\ldots$, where $Q_0=F_0^{-1}$ is the quantile function of $F_0$. It follows from \re{PnOrth}-\re{gmGen} that $\{g_m\}$ forms an orthonormal system of basis functions on $[0,1]$, i.e.\
\beq
\int_0^1 g_m(p)g_l(p)dp = 1(m=l).
\lb{gOrth}
\eeq
Equations \re{TmDef} and \eqref{gOrth} imply that $\{T_m(F)\}_{m=1}^\infty$ are the  coefficients of an orthonormal series expansion of $Q$. Indeed, from Theorem 2.3 of Takemura (1983) we find that if $F$ has a finite second moment
\beq
Q(p) = \sum_{m=1}^\infty T_m(F)g_m(p), 
\lb{QExp}
\eeq
for all continuity points $p$ of $Q(\cdot)$, strictly between 0 and 1. This is to say that $\{T_m(F)\}_{m=1}^\infty$ will quantify all aspects of the quantile function $Q$ of $F$. Takemura (1983) used \re{QExp} in order to define goodness-of-fit tests of a location-scale family \re{LocScale}, where $F_0$ in equation \re{LocScale} is also the reference distribution of the series expansion. If we restrict ourselves to functionals up to order 4, it is clear that $\{T_m(F)\}_{m=1}^4$ carries the same information about $F$ as the four types of functionals of Section \ref{Sec:ExLFunc}, i.e.\ $T_1(F)$, $T_2(F)$, $T_{32}(F)$, and $T_{42}(F)$. 
\par\medskip
In order to verify that \re{TmDef} and \eqref{gmGen} for $m=1,2,\ldots$ define a valid collection of $L$-functionals, the first condition \re{GSign} is equivalent to each polynomial $P_k(y)$ having $k$ distinct zeros with a leading positive coefficient of $y^k$. The second condition \re{GInt} follows by choosing $l=1$ in \re{gOrth}, since $g_1(p)= 1$. The third symmetry condition \re{GSymm} holds whenever $F_0$ is symmetric, and if $P_k(y)$ is an even (odd) function of $y$ when $k$ is even (odd). 

\subsubsection{Symmetric and asymmetric collections of $L$-functionals} \label{sec:symm}
A collection $\{T_m\}_{m=1}^\infty$ of $L$-functionals is symmetric if \re{GSymm} holds for all $T_m$, and otherwise it is asymmetric. For a symmetric collection of functionals $\{T_m\}_{m=1}^\infty$ we introduce a symmetric reference distribution $F_0$ with mean 0 and variance 1 that satisfies
\beq
T_m(F_0) = \left\{\begin{array}{ll}
1, & m=2,\\
0, & m\ne 2.
\end{array}\right.
\lb{F0Def}
\eeq
For a symmetric collection of $L$-functionals, we also introduce one distribution $F_1$ that is moderately skewed to the right, and another symmetric distribution $F_2$, which is more heavytailed than $F_0$, such that
\beq
\begin{array}{rcl}
T_{32}(F_1) &=& 1,\\
T_{42}(F_2) &=& 1.
\end{array}
\lb{F1F2Def}
\eeq
\par\medskip
Let us verify that a symmetric polynomial collection \re{gOrth} of $L$-functionals satisfies  \eqref{F0Def} and \eqref{F1F2Def}. Starting with \eqref{F0Def}, recall that $F_0$ is assumed to be symmetric with expected value 0 and variance 1. Applying \re{gOrth} with $m,l\in\{1,2\}$ we find that $g_1(p)=1$ and $g_2(p)=Q_0(p)$. Then a second application of \re{gOrth} with variable $m$ and $l=2$ implies \re{F0Def}. It is also possible to relax the orthonormality condition \re{gOrth} and adjust $g_3(\cdot)$ and $g_4(\cdot)$ by multiplicative constants in such a way that \re{F1F2Def} holds for some appropriately chosen distributions $F_1$ and $F_2$.  
\par\medskip
When $F$ is a life length distribution we typically choose a non-symmetric reference distribution $F_0$, supported on $[a,b)=[0,\infty)$, 
and use a non-symmetric collection of $L$-functionals. Then we replace condition \re{F0Def} by the milder requirement  
\beq
T_2(F_0) = 1.
\lb{T2NonSymm}
\eeq
Since skewness and kurtosis are somewhat less natural concepts for life lengths, we will not always impose condition \re{F1F2Def} in this context. For instance, we are typically more interested in quantifying how thick the right tail of $F$ is (as in \re{TheavyRight}-\re{GHill}), and this requires combined information from all of $\{T_m(F)\}_{m=1}^4$.

\subsubsection{Examples of polynomial collections of $L$-functionals} \label{sec:poly}
In this subsection we give three collections $\{T_m\}_{m=1}^\infty$ of $L$-functionals based on \re{TmDef} and \re{gmGen}.    
\par\medskip  
\begin{exa}[Legendre collection of $L$-functionals]\lb{LfuncLegendre}
\hspace{10pt} This collection of symmetric $L$-func\-tionals has absolutely continuous weight functions  
\beq
g_m(p)=\sqrt{2m-1}L_{m-1}(2p-1),
\lb{gLegendre}
\eeq
where $L_{m-1}$ is the Legendre polynomial of order $m-1$ on $[-1,1]$.  Thus we have that 
$$
\begin{array}{rcl}
g_1(p) &=& 1,\\
g_2(p) &=& \sqrt{3}(2p-1),\\
g_3(p) &=& \sqrt{5}[3(2p-1)^2-1]/2,\\
g_4(p) &=& \sqrt{7}[5(2p-1)^3-3(2p-1)]/2, 
\end{array}
$$
so that $T_1(F)=E(Y)$ equals the mean \re{TEY} and $T_2(F)$ is proportional to Gini's mean difference \re{TGini}.
\par\medskip
The fact that $g_m$ gives rise to an $L$-functional of order $m$ follows from the general construction in \re{gmGen}, with a uniform reference distribution $F_0\sim U(-\sqrt{3},$ $\sqrt{3})$. Notice in particular that this symmetric distribution has first two moments $E(Y)=0$ and $E(Y^2)=1$, as required above \re{F0Def}. Moreover, since $Q_0(p)=\sqrt{3}(2p-1)=g_2(p)$ it follows from \re{TmDef} and \re{gOrth} that \re{F0Def} holds for all $m$. 
Equation \re{gOrth} follows from well known orthogonality properties of Legendre polynomials, whereas \re{GSign} is referred to as the interlacing property of Legendre polynomials. It is possible to adjust $g_3$ and $g_4$ by multiplicative constants so that \re{F1F2Def} holds for some appropriate reference distributions, for instance the beta distributions $F_1\sim B(2,1)$ and $F_2\sim B(0.5,0.5)$.  
\par\medskip
It turns out that the Legendre system of $L$-functionals is equivalent to the $L$-moments $\ga_m(F)=T_m(F)/\sqrt{2m-1}$ of $F$ or order $1,2,\ldots$, introduced by Sillito (1969). Likewise, $T_{32}(F)$ and $T_{42}(F)$ are equivalent to the $L$-skewness $\tau_3=\ga_3/\ga_2$ and $L$-kurtosis $\tau_4=\la_4/\la_2$ of Hosking (1990). Another closely related concept is the family of probability weighted moments
$$
\mu_{qrs}(F) = E\left[Y^q F(Y)^r(1-F(Y))^s\right]
$$
of Greenwood et al.\ (1979), where $F=F_Y$ is the distribution function of $Y$. In fact, it can be seen that each $\mu_{1rs}(F)$ is a linear combination of Legendre functionals $\{T_m(F)\}_{m=1}^{r+s+1}$ of $F$ up to order $r+s+1$.
\par\medskip
We argue that the Legendre system \re{gLegendre} of $L$-functionals is appropriate for bounded random variables $Y\in [a,b]$, and Hosking (1990) proves that\break $\{T_m(F)\}_{m=1}^\infty$ exist for distributions $F$ with a finite mean. In spite of this, each $T_m$ is non-robust with a breakdown point of 0. It is possible though to define a robustified collection $\{T_m^\pi\}_{m=1}^\infty$ of $L$-functionals, for each $0<\pi < 0.5$, with weight functions
\beq
g_m^\pi(p) = \frac{\sqrt{2m-1}}{1-2\pi}L_{m-1}\left(\frac{2p-1}{1-2\pi}\right)1(\pi < p < 1-\pi)
\lb{gmla}
\eeq
and breakdown point $\pi$. Notice in particular that $g_1$ equals the trimmed mean, i.e.\ $\pi_0=\pi$ and $\pi_1=1-\pi$ in \re{GCCE}. In analogy with \re{QExp}, one finds that $\{T_m^\pi(F)\}_{m=1}^\infty$ provides information about the conditional distribution of $F$ for all quantiles between $\pi$ and $1-\pi$, since 
$$
\frac{Q(p)}{1-2\pi} = \sum_{m=1}^\infty T_m^\pi(F)g_m^\pi(p), 
$$
at all continuity point $\pi<p<1-\pi$ of $Q(\cdot)$. In order for $F_0^\pi$ to serve as a reference \re{F0Def} for $\{T_m^\pi\}_{m=1}^\infty$, it is required that $Y \mid Q_0^\pi(\pi)<Y<Q_0^\pi(1-\pi) \sim U(-\sqrt{3},\sqrt{3})$ whenever $Y\sim F_0^\pi$ and $Q_0^\pi = (F_0^\pi)^{-1}$. More generally, we argue that the robustified Legendre collection of $L$-functionals is appropriate whenever $Y \mid Q(\pi)<Y<Q(1-\pi)$ is bounded within some finite interval $[a,b]$. For instance, if $F_{Y \mid Q(\pi)<Y<Q(1-\pi)}\sim U(a,b)$ is uniform, then $T_1^\pi(F_Y)=(a+b)/2$ and $T_2^\pi(F_Y)=(a-b)/(2\sqrt{3})$ equal the mean and standard deviation of this uniform distribution, whereas $T_m^\pi(F_Y)=0$ for $m\ge 2$. Other ways of robustifing the Legendre system of $L$-functionals have been proposed by Mudholkar and Huston (1988) and Elamir and Seheult (2003). 
\slut\end{exa}
\par\medskip  
\begin{exa}[Hermite collection of $L$-functionals]\lb{Sec:Hermite}
In this example we introduce a collection of symmetric $L$-functionals for which $F_0\sim N(0,1)$ is a reference distribution. The weight densities 
\beq
g_m(p) = \frac{1}{\sqrt{(m-1)!}} H_{m-1}[Q_0(p)]
\lb{gHermite}
\eeq
are defined in terms of the quantile function $Q_0=F_0^{-1}$ and the probabilistic Hermite polynomials 
\beq
H_k(y) = (-1)^k e^{y^2/2} \frac{d^k e^{-y^2/2}}{dy^k}, \quad -\infty < y < \infty,
\lb{Hermite}
\eeq
of order $k=0,1,2,\ldots$. Inserting \re{Hermite} into \re{gHermite}, we find that the first four weight densities have the form
\beq
\begin{array}{rcl}
g_1(p) &=& 1,\\
g_2(p) &=& Q_0(p),\\
g_3(p) &=& [Q_0(p)^2-1]/\sqrt{2},\\
g_4(p) &=& [Q_0(p)^3-3Q_0(p)]/\sqrt{6}.
\end{array}
\lb{gHermite2}
\eeq
Recall from \re{gLocOpt} and \re{gOptScale} that $g_1$ and $g_2$ correspond to optimal $L$-functionals of location and scale for $F_0$. 
It follows from well known orthogonality properties of Hermite polynomials that the weight functions in \re{gHermite} constitute an orthonormal system \re{gOrth} on $[0,1]$, and the series expansion \re{QExp} of $Q$ can be interpreted as a robust Cornish-Fisher expansion (Fisher and Cornish, 1960). In particular, since $g_2(p)=Q_0(p)$, it follows from \re{TmDef} and \re{F0Def} that $F_0$ indeed is a reference distribution for the Hermite class of $L$-functionals. As in Example \ref{LfuncLegendre}, one may multiply $g_3$ and $g_4$ by constants so that \re{F1F2Def} holds for some appropriate reference distributions, for instance a non-central $t$-distribution $F_1$ and central $t$-distribution $F_2$, with appropriate parameters. 
\par\medskip
Robustified versions $T_m^\pi(F)$ of the Hermite functionals are constructed in the same way as in Example \ref{LfuncLegendre}, with weight functions
$$
g^\pi_m(p) = \frac{1}{(1-2\pi)\sqrt{(m-1)!}} H_{m-1}[Q_0(\frac{p}{1-2\pi})]1(\pi < p < 1-\pi).
$$ 
The reference \re{F0Def} of this system is the improper mixture distribution
$$
F_0^\pi \sim \pi\de_{-\infty} + (1-2\pi)N(0,1) + \pi\de_\infty,
$$   
with probabilities $\pi$ at minus and plus infinity. 
\par\medskip
We argue that the (robustified) Hermite system of $L$-functionals is appropriate whenever $Y\in\R$, that is, when there are no upper or lower bound restrictions on $Y$. 
\slut\end{exa}
\par\medskip
\begin{exa}[Laguerre collection of $L$-functionals]\lb{Sec:AsymL} \hspace{10pt}
In this example we consider an asymmetric collection of $L$-functionals, which is of interest when $Y\ge 0$, for instance a lifetime. This class of $L$-functionals will have $F_0\sim\Exp(1)$ as a reference distribution. Let $\mbox{La}_k(y)$ be the Laguerre polynomial of degree $k=0,1,2,\ldots$. These polynomials form an orthonormal system on $[0,\infty)$ with respect to the density measure $dF_0(y)=e^{-y}dy$, i.e.
\beq
\int_0^\infty \mbox{La}_k(y)\mbox{La}_l(y) e^{-y}dy = 1(k=l).
\lb{LanOrth}
\eeq
The weight function of $T_m$ is defined as
\beq
g_m(p) = (-1)^{m-1}\mbox{La}_{m-1}(Q_0(p)) = (-1)^{m-1}\mbox{La}_{m-1}\left(\log \left( (1-p)^{-1}\right)\right),
\lb{gExp}
\eeq
for $m=1,2,\ldots$. The factor $(-1)^{m-1}$ of \re{gExp} ensures that all weight functions $g_m(p)$ will have a leading positive coefficient of $Q_0^{m-1}(p)$, as required by \re{GSign}. It can be seen that the first four weight functions are 
$$
\begin{array}{rcl}
g_1(p) &=& 1,\\
g_2(p) &=& Q_0(p)-1,\\
g_3(p) &=& Q_0^2(p)/2 - 2Q_0(p) + 1,\\
g_4(p) &=& [Q_0^3(p)-9Q_0^2(p)+18Q_0(p)-6]/6.
\end{array}
$$ 
It follows from \re{PnOrth} and \re{gExp} that $\{g_m\}_{m=1}^\infty$ forms an orthonormal system \re{gOrth} of weight functions on $[0,1]$, with   
$$
T_m(F_0) = \int_0^1 Q_0(p)g_m(p)dp = \int_0^1 (g_1(p)+g_2(p))g_m(p)dp
= \left\{\begin{array}{ll}
1, & m=1,2,\\
0, & m\ge 3,
\end{array}\right.
$$
in agreement with \re{T2NonSymm}. Hence, according to this definition $T_{32}(F_0)=0$ although $F_0$ is skewed to the right. The rationale is that most lifetime distributions are skewed in this direction. Using $F_0$ as a yardstick we may therefore interpret $T_{32}(F)>0$ as $F$ being {\sl more} skewed to the right than $F_0$.
\par\medskip
When $Y\ge 0$ it is typically only the heaviness of the right tail of $F=F_Y$ that is of interest. A functional for right-heavytailedness is 
$$
T_{\scr{heavyr}}(F) = \frac{T_{4,\scr{right}}(F)}{T_2(F)},
$$
where $T_{4,\scr{right}}(F)$ has weight function
$$
g_{4,\scr{right}}(p) = g_4(p)1(1-\pi<p<1),
$$
and $0<\pi<1$ is the smallest positive integer satisfying $\int_{1-\pi}^1 Q_0(p)g_4(p)dp=0$. This definitions guarantees that $T_{\scr{heavyr}}(F_0)=0$, so that $F_0$ is a reference for right-heavytailedness.  
\slut\end{exa}

\section{$L$-functionals for regression models}\lb{Sec:Regr}

In this section we regard $Y$ as the outcome variable of a regression model with a vector $\bx=(x_1,\ldots,x_q)^T$ of covariates. Let
\beq 
F_{Y \mid \sbx}(y)=P(Y\le y \mid \bx), \quad -\infty < y < \infty,
\lb{FYx}
\eeq
be the conditional distribution function of $Y$ given $\bx$. By inverting this function we obtain the $p$:th conditional quantile (CQ)  
\beq
Q(p \mid \bx) = F_{Y \mid \sbx}^{-1}(p) = \inf\{y;\, F_{Y \mid \sbx}(y)\ge p\}.
\lb{Qpx}
\eeq
Each functional \re{T} gives rise to a linear combination 
\beq
\theta(\bx) = T(F_{Y \mid \sbx}) = \int_0^1 Q(p \mid \bx)dG(p)
\lb{theta}
\eeq
of CQs. We will refer to \re{theta} as a conditional $L$-functional. The weight functions of Section \ref{sec:Lstat} give rise to different conditional $L$-functionals that correspond to location, scale, unstandardized skewness or unstandardized heavytailedness measures of  $F_{Y \mid \sbx}$. For each $0<p<1$, let $\hQ(p \mid \bx)$ be an estimator of the conditional quantile \re{Qpx}, based on a sample of size $n$. The corresponding estimator  
\beq
\hth(\bx) = \int_0^1 \hQ(p \mid \bx)dG(p)
\lb{hthx}
\eeq
of $\theta(\bx)$ is a conditional $L$-statistic that reduces to \re{hthLGen} for a model without any covariates.  
\par\medskip
In order to study standardized measures of skewness and heavytailedness of $F_{Y \mid \sbx}$, we need to consider ratios of two conditional $L$-functionals $T^1(F_{Y \mid \sbx})$ and $T^2(F_{Y \mid \sbx})$ that involve two different weight functions $G^1$ and $G^2$. We will therefore study quantities
\beq
\theta(\bx) = \frac{\int_0^1 Q(p \mid \bx)dG^1(p)}{\int_0^1 Q(p \mid \bx)dG^2(p)},
\lb{thxRatio0}
\eeq
and their estimators
\beq
\hth(\bx) = \frac{\int_0^1 \hQ(p \mid \bx)dG^1(p)}{\int_0^1 \hQ(p \mid \bx)dG^2(p)}.
\lb{hthxRatio0}
\eeq
For a model with covariates, the definitions of \re{hthx} and \re{hthxRatio0} will depend on whether the response variable $Y$ is censored/truncated or not, and on the type of regression model that is used. In the following subsections we will consider several such models.  

\subsection{Data without censoring/truncation}

In this subsection we assume there is no censoring or truncation, so that the response variable $Y$ is observed. In more detail, suppose that a sample of independent random vectors $(\bx_1,Y_1),\ldots,(\bx_n,Y_n)$ is available, with the same conditional distribution $F_{Y_i \mid \sbx_i}=F_{Y \mid \sbx_i}$ of the response variable as in \re{FYx}. In order to estimate $\theta(\bx)$ in \re{theta} or \re{thxRatio0} we have to make some smoothness assumptions on $\bx\to\theta(\bx)$, which in turn requires smoothness of the conditional quantiles $\bx\to Q(p \mid \bx)$. The most general approach is to estimate $Q(p \mid \bx)$ (and hence also $\theta(\bx)$) by some nonparametric method, for instance local polynomial regression (Chauduri, 1991), smoothing splines (Koenker et al., 1994), regression splines (He and Shi, 1994), piecewise polynomial regression tree methods (Chauduri and Loh, 2002), or a semiparametric linear model with varying coefficients (Kim, 2007). In the following two subsections we will rather concentrate on two fully parametric models for the relation between $\bx$ and $Q(p \mid \bx)$; linear models and transformed linear models.      

\subsubsection{Linear models}\lb{Sec:Lin}

{\bf 3.1.1.1 Models and estimators}
\par\smallskip
When the outcome variable $Y\in (-\infty,\infty)$ is unbounded, it is natural to use a linear model 
\beq
Q(p \mid \bx) = \bx^T\bbe(p),
\lb{QpxLin}
\eeq
so that each conditional quantile \re{Qpx} is a linear function of the covariates, with $\bbe(p)=(\beta_1(p),\ldots,\beta_q(p))^T$ a vector of regression parameters. The simplest special case of \re{QpxLin} is the homoscedastic linear model, with 
\beq
F_{Y \mid \sbx}(y) = F_0(y-\bx^T\bb),
\lb{FHomosc}
\eeq
for some vector $\bb=(b_1,\ldots,b_q)^T$. If the first regression component is an intercept ($\bx=(1,x_2,\ldots,x_q)^T$) and $F_0$ has median 0, it follows that $\bbe(p)=\bb+\left[Q_0(p)-Q_0(0.5)\right]\bfm{e}_1$, where $\bfm{e}_1=(1,0,\ldots,0)^T$ and $Q_0$ is the inverse of $F_0$. In particular, the framework of Section \ref{sec:Lstat} is a special case of \re{FHomosc} with $q=1$ and $b_1$ the median of $F_Y$. It is also possible to incorporate heteroscedastic regression models of type
\beq
F_{Y \mid \sbx}(y) = F_0\left(\frac{y-\bx^T\bb}{\bx^T\bc}\right)
\lb{FHeterosc}
\eeq
into \re{QpxLin}, with $\bc = (c_1,\ldots,c_q)^T$, for all $\bx$ such that $\bx^T\bc > 0$. Again, it follows that $\bbe(p)=\bb+\left[Q_0(p)-Q_0(0.5)\right]\bc$ if $F_0$ has median 0.
\par\medskip
The first $L$-based inference methods of regression focused on estimating $\bb$ for the linear and homoscedastic model \re{FHomosc}, based on some preliminary estimate (Bickel, 1973, Ruppert and Carroll, 1980, Welsh, 1987). Gutenbrunner and Jure\v{c}kov\'{a} (1992) introduced another approach based on regression rankscores. Here we will focus on procedures that estimate the conditional quantile \re{QpxLin} by first finding an estimator of $\bbe(p)$. The most general such estimator is a solution of the minimization problem
\beq
\hbbe(p) = (\hbe_1(p),\ldots,\hbe_q(p))^T = \mbox{arg}\min_{\sbb\in\R^q} \sum_{i=1}^n \rho_p(Y_i-\bx_i^T\bb),
\lb{hbbep}
\eeq
where $\rho_p(y) = [p-1(y<0)]y$ is the so called check function. This estimator is nonparametric in the sense that it makes few assumptions about the conditional distribution of $Y | \bx$, apart from its linear dependency on $\bx$ in \re{QpxLin}. It is usually referred to as a regression quantile and for the model of Section \ref{sec:Lstat}, with no explanatory variables, it simplifies to the sample quantile $Y_{([np])}$. Koenker and Bassett (1978) introduced \re{hbbep} for the homoscedastic regression model \re{FHomosc}, and later it was extended by Koenker and Bassett (1982) and Koenker and Zhao (1994) to the heteroscedastic model \re{FHeterosc}. Efron (1991) proposed a slightly different nonparametric estimate of $\bbe(p)$, by first minimizing an asymmetric squared loss function of the residuals $Y_i-\bx_i^T\bb$.
\par\medskip
Frumento and Bottai (2016) modeled 
\beq
\bbe(p) = \bbe(p;\bpsi) = \sum_{k=1}^r \bpsi_k\beta^k(p)
\lb{bbepPar}
\eeq
parametrically in terms of a $q\times r$ matrix $\bpsi=(\bpsi_1,\ldots,\bpsi_r)$, where $\beta^1(p),\ldots,$ $\beta^r(p)$ are known functions of $p$. This includes, for instance, the homoscedastic and heteroscedastic regression models \re{FHomosc} and \re{FHeterosc}, when $F_0$ is regarded as known. Both of these models have $\beta^1(p)\equiv 1$, $\beta^2(p)=Q_0(p)-Q_0(0.5)$, and $\bpsi_1=\bb$, whereas $\bpsi_2=\bfm{e}_1$ for the homoscedastic and $\bpsi_2=\bc$ for the heteroscedastic model. The model in \re{bbepPar} gives rise to a parametric regression quantile estimator, where
\beq
\hbbe(p) = \bbe(p;\hbpsi) = \sum_{k=1}^r \hbpsi_k\beta^k(p),
\lb{hbbep2}
\eeq
and $\hbpsi_k$ is an estimator of $\bpsi_k$. Frumento and Bottai (2016) used 
\beq
\hbpsi = \arg\min_{\bpsi} \sum_{i=1}^n \int_0^1 \rho_p(Y_i-\bx_i^T\bbe(p;\bpsi)) dp.
\lb{hbpsi}
\eeq 
The conditional $L$-functional \re{theta} of a linear model has a very tractable form. Recall first of all from Example \ref{Sec:Hermite} that the Hermite system \re{gHermite} of weight functions $G=G_m$ is a natural choice of $L$-functionals of order $m=1,2,\ldots$ for a linear model, whenever the range of the outcome variable is unbounded. Moreover, since the conditional quantile \re{QpxLin} is a linear function of $\bx$, this linearity is preserved for conditional $L$-functionals \re{theta}. Indeed, inserting \re{QpxLin} into \re{theta}, we find that  
\beq 
\theta(\bx) = \bx^T \int_0^1 \bbe(p)dG(p) =: \bx^T\bB,
\lb{thetaLin}
\eeq
where $\bB=(B_1,\ldots,B_q)^T$ and $B_j=\int_0^1 \beta_j(p)dG(p)$. 
\par\medskip
In order to estimate the conditional $L$-functional $\theta(\bx)$, the chosen estimator $\hbbe(p)$ (for instance the nonparametric \re{hbbep} or the parametric \re{hbbep2}) is first plugged into \re{QpxLin} and then into \re{hthx}. This gives an estimated conditional quantile $\hQ(p \mid \bx) = \bx^T\hbbe(p)$ and an estimate
\beq
\hth(\bx) = \bx^T\hbB 
\lb{hthLin}
\eeq
of $\theta(\bx)$, where $\hbB = (\hB_1,\ldots,\hB_q)^T = \int_0^1 \hbbe(p)dG(p)$ is a linear combination of all $\hbbe(p)$. In the context of regression quantiles, the estimator $\hbB$ was first proposed by Koenker and Portnoy (1987) for the homoscedastic model \re{FHomosc}, and then extended to heteroscedastic models \re{FHeterosc} by Koenker and Zhao (1994). Garc\'{i}a-Pareja and Bottai (2018) studied \re{thetaLin} for weight functions that correspond to a uniform distribution on $[\pi_0,\pi_1]$ for some $0\le \pi_0 < \pi_1 \le 1$. They referred to $\theta(\bx)$ as a conditional compound expectation, and it generalizes the compound expectation \re{GCCE} of the location model. 
\par\medskip
{\bf 3.1.1.2 Asymptotics}
\par\smallskip
In order to study the large sample asymptotics of $\hth(\bx)=\hth_n(\bx)$, we follow Garc\'{i}a-Pareja and Bottai (2018) and view $\bbe=\{\bbe(p);\ 0<p<1\}$ as an element of the $q$-dimensional product space $S^q$. In our setting the norm  
$$
\|\bbe\|_{S^q} = \sum_{j=1}^q \|\beta_j\|_S
$$
of this space generalizes \re{QS} from $q=1$ to $q\ge 1$. For $q$-dimensional vectors $\bB\in\R^q$ we introduce the $L^1$-norm $|\bB| = \sum_{j=1}^q |B_j|$. Moreover, if $W,W_1,\ldots,W_n,\ldots$ are random variables of a metric space $\cW$, equipped with norm $|\cdot|$, we say that the sequence $W_n$ converges in probability towards $W$ as $n\to\infty$, i.e.\ $W_n \pto W$, if $P(|W_n-W|>\va)\to 0$ for each $\va>0$. Equipped with these preliminaries, we have the following:
\par\medskip
\begin{prop}[Consistency of \re{hthLin}.]\lb{Prop:Cons}
Suppose $\hbbe=\hbbe_n\in S^q$ is a consistent estimator of $\bbe$, i.e.\ $\hbbe_n\pto\bbe$ as $n\to\infty$.
Then $\hbB=\hbB_n\pto \bB$ and the estimated conditional $L$-functional in \eqref{hthLin} is consistent, i.e. $\hth_n(\bx)\pto \theta(\bx)$ as $n\to\infty$.  
\end{prop}
\par\medskip
{\bf Proof.} In analogy with \re{TDiff} we have the inequalities
\beq
\begin{split}
|\hbB_n-\bB| &\le  \sum_{j=1}^q \int_0^1 |\hbe_j(p)-\beta_j(p)|d|G|(p)\\
&= \sum_{j=1}^q \left( \int_0^1 |\hbe_j(p)-\beta_j(p)||g(p)|dp 
+ \sum_{m=1}^M |g_m||\hbe_j(\pi_m)-\beta_j(\pi_m)|\right)\\
&= \| \hbbe_n - \bbe\|_{S^q}
\end{split}
\lb{hbBnCons}
\eeq
By assumption, $P(\|\hbbe_n-\bbe\|_{S^q}>\va)\to 0$ as $n\to\infty$ for each $\va>0$. In conjunction with \re{hbBnCons} it follows that $\hbB_n\pto \bB$. Finally, since $\hth_n(\bx)=h(\hbB_n)$ and $\theta(\bx)=h(\bB)$ for the continuous function $h:\R^q\to\R$, defined by $h(\bB)=\bx^T\bB$, consistency $\hth_n(\bx)\pto \theta(\bx)$ follows by the Continuous Mapping Theorem.
\hfill\slut
\par\medskip
As a next step, in order to establish asymptotic normality of $\hth(\bx)=\hth_n(\bx)$ as $n\to\infty$ we first generalize \re{Znp} from the location model ($q=1$, $\bx_i\equiv 1$) and introduce the rescaled process 
\beq
\bZ_n(p) = \sqrt{n}[\hbbe_n(p)-\bbe(p)], \quad 0 < p < 1
\lb{Gnp}
\eeq
of $\hbbe_n(p)$. Viewing $\bZ_n=\{\bZ_n(p); \, 0<p<1\}$ as a random element of $S^q$, we will assume weak convergence 
\beq
\bZ_n \Lto \bZ \mbox{ as }n\to\infty
\lb{GntoG}
\eeq
with respect to the topology in $S^q$ introduced by $\|\cdot\|_{S^q}$. The limit in \re{GntoG} is assumed to be a Gaussian process $\bZ=\{\bZ(p); \, 0<p<1\}$ with mean $E(\bZ(p))=(0,\ldots,0)^T$ for all $0<p<1$ and covariance function $\bR(p,s)=\Cov(\bZ(p),\bZ(s))$ for all $0< p,s < 1$. 
\par\medskip
In the context of regression quantiles \re{hbbep}, in order to find the covariance function $\bR(p,s)$, it is helpful to rewrite \re{Gnp} as
$$
\bZ_n(p) = \frac{1}{\sqrt{n}}\bD_n^{-1}(p)\sum_{i=1}^n \rho_p^\prime(Y_i-\bx_i^T\bbe(p)) + o_p(1),
$$
where 
\beq
\begin{array}{rcl}
\bD_n(p) &=& \frac{\sum_{i=1}^n f_i(p)\bx_i\bx_i^T}{n} \\
f_i(p) &=& \frac{dF_{Y \mid \sbx_i}(y)}{dy} \bigg|_{y=Q(p \mid \sbx_i)} \\
\rho_p^\prime(y) &=& p - 1(y<0)
\end{array}
\eeq
and $o_p(1)$ is an asymptotically negligible remainder term, i.e.\ $o_p(1)\pto 0$ as $n\to\infty$, uniformly for $p\in \mbox{supp}(G)$. Then assume there exist positive definite matrices $\bA$ and $\bD(p)$ such that 
\beq
\begin{array}{cl}
(i) & \lim_{n\to\infty} \sum_{i=1}^n \bx_i\bx_i^T/n = \bA,\\
(ii) & \lim_{n\to\infty} \bD_n(p) = \bD(p),\\
(iii) & \max_{i=1,\ldots,n} |\bx_i|/\sqrt{n} \to 0.
\end{array}
\lb{RegCond}
\eeq
The regularity conditions in \re{RegCond} imply that the Gaussian limit process has a covariance function
\beq
\bR(p,s) = [\min(p,s)-ps]\bD^{-1}(p)\bA\bD^{-1}(s)
\lb{R}
\eeq
for regression quantiles, see for instance Chapter 4 of Koenker (2005). This covariance function simplifies to
\beq
\bR(p,s) = \frac{\min(p,s)-ps}{f_0(Q_0(p))f_0(Q_0(s))}\bA^{-1}
\lb{RHomosc}
\eeq
for the homoscedastic regression model \re{FHomosc}, cf.\ Koenker and Portnoy (1987). We notice that \re{RHomosc} generalizes the asymptotic covariance function \re{Rpsedf} of the normalized quantile process \re{Znp} for data without covariates ($\bx_i=1$ and $\bA=1$). 
\par\medskip
For the parametric estimator \re{hbbep2}, in order to find the asymptotic covariance function of $\{\hbbe(p); \\ 0<p<1\}$, we first need to generalize \re{psiConv} and establish asymptotic normality of $\hbpsi=\hbpsi_n$. To this end, it is convenient to introduce $\mbox{vec}(\hbpsi)=(\hbpsi_1^T,\ldots,\hbpsi_r^T)^T$, the column vector of length $rq$ in which the columns of $\hbpsi$ are stacked on top of each other. It is shown in Frumento and Bottai (2016), under regularity conditions similar to \re{RegCond}, that 
\beq
\sqrt{n}(\mbox{vec}(\hbpsi_n)-\mbox{vec}(\bpsi)) \Lto N(\bzero,\bV) = N\left(\bzero,\left(\begin{array}{ccc} \bV_{11} & \ldots & \bV_{1r} \\
\vdots & \ddots & \vdots \\ \bV_{r1} & \ldots & \bV_{rr} \end{array}\right)\right),
\lb{psiLConv} 
\eeq
where $\bV$ is a square matrix of order $rq$, and $\bV_{kl}$ is a square matrix of order $q$ that corresponds to the asymptotic covariance matrix between $\hbpsi_k$ and $\hbpsi_l$. The exact form of $\bV$ can be found in Frumento and Bottai (2016). It follows from \re{hbbep2} and \re{psiLConv} that weak convergence \re{GntoG} holds with asymptotic covariance function
\beq
\bR(p,s) = \sum_{k,l=1}^r \beta^k(p)\beta^l(s)\bV_{kl} 
\lb{RpsPar}
\eeq
of the limit process $\bZ$. 
Notice also that this covariance function is a special case of \re{Rpar} when $q=1$ and $\bx_i\equiv 1$, with $dQ(p;\bpsi)/d\bpsi = (\beta^1(p),\ldots,\beta^r(p))$.   
\par\medskip
Equipped with these preliminaries, the following result provides asymptotic normality of $\hbB_n$ and $\hth_n(\bx)$:
\par\medskip
\begin{prop}[Asymptotic normality of \re{hthLin}]\lb{Prop:BnAsN} \hspace{6pt}
Suppose \re{GntoG} holds for some Gaussian li\-mit process $\bZ$ whose covariance function $\bR(p,s)$ is the asymptotic covariance function of $\{\hbbe_n(p);\\ 0<p<1\}$. Suppose further that the function $g$ in \re{GDef} is bounded by $\|g\|_\infty<\infty$. Then the estimator $\hbB=\hbB_n$ of $\bB$ is asymptotically normal as $n\to\infty$, in the sense that 
\beq
\sqrt{n}(\hbB_n-\bB) \Lto N(\bzero,\bSi),
\lb{BnAsN}
\eeq
where the mean vector of the $q$-dimensional limiting normal distribution is $\bzero=(0,\ldots,0)^T$ and the covariance matrix equals  
\beq
\bSi = \int_0^1\int_0^1 \bR(p,s)dG(p)dG(s).
\lb{bSi}
\eeq
Moreover, the estimator $\hth(\bx)=\hth_n(\bx)$ of the conditional $L$-functional in \re{hthLin} is also asymptotically normal, with
\beq
\sqrt{n}[\hth_n(\bx)-\theta(\bx)] \Lto N(0,\bx^T\bSi\bx)
\lb{hthnAsN}
\eeq
as $n\to\infty$. 
\end{prop}
\par\medskip
{\bf Proof.} Introduce the functional $\bh:S^q\to \R^q$ by means of $\bh(\bz)=\!\int_0^1\! \bz(p)dG(p)$. From the proof of Proposition \ref{Prop:Cons}, and the fact that $\|g\|_\infty<\infty$, we know that\break $\bh$ is a continuous functional. Since $\sqrt{n}(\hbB_n-\bB) = \bh(\bZ_n)$ and $\bh(\bZ)\sim N(\bzero,\bSi)$, equation \re{BnAsN} is a consequence of the Continuous Mapping Theorem. Then \re{hthnAsN} follows from \re{BnAsN} by a second application of the Continuous Mapping Theorem since $\sqrt{n}[\hth_n(\bx)-\theta(\bx)]=\sqrt{n}\bx^T(\hbB_n-\bB)=h[\sqrt{n}(\hbB_n-\bB)]$, using the function $h:\R^q\to\R$ defined by $h(\bB)=\bx^T\bB$.  
\hfill\slut
\par\medskip
Propositions \ref{Prop:Cons}-\ref{Prop:BnAsN} are applicable whenever the goal is to estimate conditional $L$-functionals that correspond to measures of location, scale, unstandardized skewness or unstandardized heavytailedness of $F_{Y \mid \sbx}$, defined as in \re{theta}. However, in order to study standardized measures of skewness and heavytailedness of $F_{Y \mid \sbx}$ we consider quantities \re{thxRatio0} that are defined as a ratio of two conditional $L$-functionals $T^1(F_{Y \mid \sbx})$ and $T^2(F_{Y \mid \sbx})$, with different weight functions $G^1$ and $G^2$. For the linear model \re{QpxLin} we find that 
\beq
\theta(\bx) = \frac{\bx^T\bB^1}{\bx^T\bB^2},
\lb{thxRatio}
\eeq
with $\bB^k = \int_0^1 \bbe(p)dG^k(p)$ for $k=1,2$. The corresponding estimator
\beq
\hth(\bx) = \frac{\bx^T\hbB^1}{\bx^T\hbB^2},
\lb{hthxRatio}
\eeq
is defined analogously with $\hbB^k = \int_0^1 \hbbe(p)dG^k(p)$ for $k=1,2$. The following proposition shows that $\hbB^1=\hbB^1_n$ and $\hbB^2=\hbB^2_n$ are jointly asymptotically normal, and as a consequence, that $\hth(\bx)=\hth_n(\bx)$ in \re{hthxRatio} is asymptotically normal as well. 
\par\medskip
\begin{prop}[Consistency and asymptotic normality of \re{hthxRatio}]\lb{Prop:JointBnAsN}
\hspace{5pt}\break Suppose the regular\-ity conditions of Propositions \ref{Prop:Cons} and \ref{Prop:BnAsN} hold. Then $\hbB^1_n$ and $\hbB^2_n$ are consistent and (jointly) asymptotically normal estimators of $\bB^1$ and $\bB^2$ as $n\to\infty$, in the sense that
\beq
\sqrt{n}[(\hbB_n^1,\hbB_n^2)-(\bB^1,\bB^2)] \Lto N\left(\bzero,\left(\begin{array}{cc} \bSi^{11} & \bSi^{12}\\ \bSi^{21} & \bSi^{22}\end{array}\right)\right),
\lb{BnJointAsN}
\eeq
where $\bzero=(0,\ldots,0)^T$ is a $2q$-dimensional vector of zeros,
\beq
\bSi^{kl} = \int_0^1\int_0^1 \bR(p,s)dG^k(p)dG^l(s)
\lb{bSiJoint}
\eeq
for $k,l=1,2$, and $\bR(p,s)$ is the asymptotic covariance function of $\{\hbbe_n(p);\, 0<p<1\}$. Moreover, the estimator $\hth_n(\bx)$ in \re{hthxRatio} of the ratio $\theta(\bx)$ of the two conditional $L$-functionals in \re{thxRatio}, is also asymptotically normal as $n\to\infty$, in the sense that
\beq
\begin{split}
\sqrt{n}[\hth_n(\bx)-\theta(\bx)] &\Lto N\left(0,\Sigma^\dagger\right) \\
\Sigma^\dagger &= \frac{\bx^T\bSi^{11}\bx} {(\bx^T\bB^2)^2}  - \frac{2\bx^T\bB^1\cdot \bx^T\bSi^{12}\bx}{(\bx^T\bB^2)^3} + \frac{(\bx^T\bB^1)^2\cdot \bx^T\bSi^{22}\bx }{(\bx^T\bB^2)^4}.
\end{split}
\lb{hthnAsNRatio}
\eeq
\end{prop}
\par\medskip
{\bf Proof.} Consistency of $\hbB^1_n$ and $\hbB^2_n$ follows as in the proof of Proposition \ref{Prop:Cons}. In order to verify \re{BnJointAsN} we first look at linear combinations 
$$
\alpha\hbB^1_n + \beta \hbB^2_n = \int_0^1 \hbbe(p)(\alpha dG^1(p) + \beta dG^2(p)) =: \int_0^1 \hbbe(p) dG(p)
$$
of $\hbB^1_n$ and $\hbB^2_n$. Then we apply Proposition \ref{Prop:BnAsN} with $G=\alpha dG^1+\beta dG^2$ in order to deduce
$$
\sqrt{n}[\alpha\hbB^1_n + \beta \hbB^2_n - (\alpha\bB^1 + \beta \bB^2)] \Lto N(\bzero,\alpha^2\bSi^{11} + 2\alpha\beta \bSi^{12} + \beta^2\bSi^{22})
$$
as $n\to\infty$ for any real-valued $\alpha$ and $\beta$. Then \re{BnJointAsN} follows from the Cram\'{e}r-Wold device.  
\par\medskip
Next we consider the estimator $\hth_n(\bx)$ of $\theta(\bx)$. Consistency $\hth_n(\bx)\pto\theta(\bx)$ follows as in the proof of Proposition \ref{Prop:Cons}, from the consistency of $(\hbB^1_n,\hbB_n^2)\pto (\bB^1,\bB^2)$ and the Continuous Mapping Theorem. In order to verify asymptotic normality \re{hthnAsNRatio}, it is convenient to denote the numerators and denominators of \re{thxRatio} and \re{hthxRatio} as $N/D$ and $\hN/\hD$ respectively. Notice first that 
$$
\sqrt{n}[(\hN,\hD)-(N,D)] \Lto 
N\left((0,0),\left(\begin{array}{cc} \bx^T\bSi^{11}\bx & \bx^T\bSi^{12}\bx \\ \bx^T\bSi^{21}\bx & \bx^T\bSi^{22}\bx \end{array}\right)\right)
$$
as $n\to\infty$, using first the Cram\'{e}r-Wold device then same argument as in the proof of \re{hthxRatio}. Then \re{hthnAsNRatio} follows after a first order two-dimensional Taylor expansion 
$$
\hth_n(\bx) \approx \theta(\bx) + \frac{1}{D}(\hN-N) - \frac{N}{D^2}(\hD-D)
$$    
of $\hth_n(\bx) = \hN/\hD =: h(\hN,\hD)$ around the point $(N,D)$, noticing that $\theta(\bx)=h(N,D)$. 
\hfill\slut

\subsubsection{Transformed linear models}\lb{Sec:TransLinear}

Suppose the response variable $Y\in [a,b]$ is constrained to lie in an interval with end points $-\infty \le a < b \le \infty$. If at least one of these two end points is finite, the conditional quantile \re{QpxLin} of the linear model may fall outside $[a,b]$ for some covariate vectors $\bx$. In order to avoid this, it is common to introduce a known and strictly increasing link function $h:[a,b]\to\R$ and assume that regression data $(\bx,h(Y))$ for the transformed outcome variable follows the linear model of Section \ref{Sec:Lin}. This is analogous to link functions of Generalized Linear Models (McCullagh and Nelder, 1989), although here we focus on transformations of quantiles rather than of expected values. Since quantiles are preserved by monotone transformations it follows from \re{QpxLin} that the conditional quantiles satisfy
\beq
Q(p \mid \bx) = h^{-1}\left(\bx^T\bbe(p)\right)
\lb{QpxLinLink}
\eeq
for all $0<p<1$. If our objective is to estimate a conditional location, conditional scale, conditional unstandardized skewness or conditional unstandardized kurtosis of $Y$, we use the conditional $L$-functional in \re{theta}. For the transformed linear model \re{QpxLinLink}, this functional  
\beq
\theta(\bx) = \int_0^1 h^{-1}(\bx^T\bbe(p))dG(p)
\lb{thetaLink}
\eeq
is no longer a linear function of $\bx$, as in \re{thetaLin}. In order to estimate $\theta(\bx)$ from data, we proceed as in Section \ref{Sec:Lin}. We first estimate the regression parameter $\bbe(p)$ of the transformed data set by some method. One possibility is to estimate $\bbe(p)$ nonparametrically with regression quantiles 
\beq
\hbbe(p) = \mbox{arg}\min_{\sbb\in\R^q} \sum_{i=1}^n \rho_p(h(Y_i)-\bx_i^T\bb),
\lb{hbbepLink}
\eeq
in analogy with \re{hbbep}. When $p=0.5$, we notice that \re{hbbepLink} is a special case of $L_1$-estimation for nonlinear regression models (Oberhofer, 1982). Alternatively we estimate $\bbe(p)$ parametrically as 
\beq
\hbbe(p) = \arg\min_{\bpsi} \sum_{i=1}^n \int_0^1 \rho_s(h(Y_i)-\bx_i^T\bbe(s;\bpsi)) ds 
\left(\begin{array}{c} \beta^1(p) \\ \vdots \\ \beta^r(p) \end{array}\right),
\lb{hbbep2Link}
\eeq
in analogy with \re{hbbep2}-\re{hbpsi}. As we will see below, for some models it is more convenient to estimate $\bbe(p)$ parametrically by maximum likelihood rather than using \re{hbbep2Link}. 
\par\medskip
By plugging a nonparametric or parametric estimate of $\bbe(p)$ into \re{thetaLink}, we finally obtain an estimator
\beq
\hth(\bx) = \int_0^1 h^{-1}(\bx^T\hbbe(p))dG(p)
\lb{hthLink}
\eeq
of $\theta(\bx)$. For some of the examples below, the link function $h$ is unknown, and then we have to estimate $h(\cdot)$ (parametrically or nonparametrically) as well as $\bbe(p)$. This gives rise to an estimator 
\beq
\hth(\bx) = \int_0^1 \hat{h}^{-1}(\bx^T\hbbe(p))dG(p)
\lb{hthLink2}
\eeq
of $\theta(\bx)$, where $\hat{h}(\cdot)$ is a monotone estimate of $h(\cdot)$. 
\par\medskip
We will now give several examples of link functions for which the corresponding transformed response variables follow a linear model, as described in \re{QpxLinLink}. 
\par\medskip
\begin{exa}[Logit transformations.]\lb{Ex:Logit} \hspace{5pt} Suppose the outcome variables is bounded, so that $-\infty<a<b<\infty$ are both finite. Liu et al.\ (2009) and Bottai et al.\ (2010) used model \re{QpxLinLink} with a logit transformation
$$
h(y) = \mbox{logit} \frac{y-a}{b-a} = \log \frac{y-a}{b-y}
$$  
for such data.   
\slut\end{exa}
\par\medskip
Since $Y$ is bounded in Example \ref{Ex:Logit}, it is possible to use $L$-functionals in \re{hthLink}-\re{hthLink2} with a weight density $dG(p)=g_m(p)dp$ chosen from the Legendre system \re{gLegendre} of Example \ref{LfuncLegendre}. The remaining examples of this section concern outcome variables such as life lengths, where $Y$ is constrained to be positive ($a=0$, $b=\infty$). Recall from Example \ref{Sec:AsymL} that it is possible then choose the weight density $dG(p)=g_m(p)dp$ in \re{hthLink}-\re{hthLink2} from the asymmetric collection \re{gExp} of basis functions, with an exponential reference distribution. The construction in Example \ref{Sec:AsymL} can be modified though so that some other lifetime distribution is used as reference.  
\par\medskip
\begin{exa}[Logarithmic transformations.]\lb{Ex:Log} \hspace{5pt} When $Y$ is non-negative\break ($a=0$, $b=\infty$) it is common to use the logarithmic link 
\beq
h(y) = \log(y).
\lb{hlink}
\eeq 
The Accelerated Failure Time (AFT) model (Kalbfleich and Prentice, 2002) is often used when $Y>0$ is a lifetime, and it corresponds to having a parametric location-scale regression model \re{FHomosc} for the transformed data $\{(\bx_i,\log(Y_i))\}_{i=1}^n$. In more detail, 
\beq
\log\left[Q(p \mid \bx)\right] = \mu + \sum_{k=2}^q x_k b_{k-1} + \sigma Q_0(p) = \bx^T\bbe(p),
\lb{hAFT}
\eeq
where $Q_0$ is the inverse of some reference distribution $F_0$ of the log lifetime. Here $\bx=(1,x_2,\ldots,x_q)^T$ is the covariate vector with an added intercept, $\bbe(p)=(\mu+\si Q_0(p))\bfm{e}_1+(0,\bb^T)^T$, and $\bb=(b_1,\ldots,b_{q-1})^T$ contains the effect parameters of the $q-1$ covariates. When $F_0$ is known, for instance a logistic distribution, equation \re{hAFT} corresponds to a parametric model \re{hbbep2Link} with $r=2$, $\bpsi_1=(\mu,b_1,\ldots,b_{q-1})^T$, $\bpsi_2=(\si,0,\ldots,0)^T$, $\beta^1(p)\equiv 1$ and $\beta^2(p)=Q_0(p)$. 
When computing $\hbbe(p)$ it is possible though to use the maximum likelihood estimates of the $q+1$ nonzero model parameters $\mu,\si,b_1,\ldots,b_{q-1}$ of $\bpsi$, rather than \re{hbbep2Link}. Garc\'{i}a-Pareja et al.\ (2019) used this approach, with a piecewise constant weight function \re{GCCE}, to estimate the conditional compound expectation $\theta(\bx)$ of an AFT model. When $F_0$ is left unspecified it is also possible to estimate $\mu$ and $\bb$ nonparametrically according to \re{hbbepLink} (Ying et al., 1995).         
\slut\end{exa}
\par\medskip
\begin{exa}[Power transformations.]\lb{Ex:Power} Assume as in Example \ref{Ex:Log} that $Y$ is positive, i.e.\ $a=0$ and $b=\infty$. Mu and He (2007) estimated conditional quantiles \re{thetaLink} with class of power link functions (Box and Cox, 1964), i.e.
\beq
h_\ga(y) = \left\{ \begin{array}{ll} 
(y^\ga-1)/\ga, & \ga\ne 0,\\
\log(y), & \ga = 0.
\end{array}\right.
\lb{hga}
\eeq
It was assumed in Mu and He (2007) that not only the regression parameter vector $\bbe(p)$, but also the parameter $\ga=\ga(p)$ of the power transformation \re{hga}, were unknown. For this reason their estimate included a combination of a CUSUM procedure and regression quantiles \re{hbbepLink}. The resulting estimated link function $\hat{h}(y;p)=h_{\hga(p)}(y)$ is then inserted into \re{hthLink2} in order to estimate $\theta(\bx)$.  
\slut\end{exa}
\par\medskip
\begin{exa}[Log cumulative baseline hazard transformations.]\lb{Ex:Cox} \hspace{3pt}
As\break in the previous two examples, consider a positive response variable $Y$, so that $a=0$ and $b=\infty$. The Cox regression model (Cox, 1972) expresses the hazard function 
\beq
\la(y \mid x_2,\ldots,x_q) = \la_{\base}(y)\exp(\sum_{k=2}^q x_k b_{k-1})
\lb{Cox}
\eeq
of the lifetime $Y$ as a product of a baseline hazard $\la_{\base}$ and a term that involves the covariates $x_2,\ldots,x_q$ and a regression vector $\bb=(b_1,\ldots,b_{q-1})^T$. Equivalently, \re{Cox} can be rewritten in terms of the cumulative hazard function as
\beq
\La(y \mid x_2,\ldots,x_q) = \int_0^y \la(z \mid x_2,\ldots,x_q)dz = \La_{\base}(y)\exp(\sum_{k=2}^q x_k b_{k-1}),
\lb{Cox2}
\eeq
where $\Lambda_{\base}(y)=\int_0^y \la_{\base}(z)dz$ is the corresponding cumulative baseline hazard. It is well known (Doksum and Gasko, 1990, Koenker and Geling, 2001, Portnoy, 2003, Garc\'{i}a-Pareja et al., 2019) that \re{Cox2} can be rewritten as a transformed linear model \re{QpxLinLink}, with a link function
\beq
h(y) = \log\left[\Lambda_{\base}(y)\right] \stackrel{y=Q(p \mid \sbx)}{=} \log\left[\log (1-p)^{-1}\right] - \sum_{k=2}^q x_k b_{k-1} = \bx^T\bbe(p)
\lb{hCox}
\eeq
that is the logarithm of the cumulative baseline hazard, $\bx=(1,x_2,\ldots,x_q)^T$ is a covariate vector with added intercept, and the regression vector is $\bbe(p)=\log\left[\log (1-p)^{-1}\right]\bfm{e}_1 - (0,\bb^T)^T$. This is an instance of a parametric model \re{hbbep2Link} for which $r=2$, $\bpsi_1=\bfm{e}_1$ is known, $\bpsi_2 = (0,\bb^T)^T$, $\beta^1(p)=\log\left[\log (1-p)^{-1}\right]$ is a quantile of a Gumbel distribution and $\beta^2(p)\equiv -1$. Garc\'{i}a-Pareja et al.\ (2019) estimate the conditional compound expectation $\theta(\bx)$ of a Cox model based on the piecewise linear weight function \re{GCCE}. 
\par\medskip
For a Cox model it is traditional to estimate $\bb$ directly by partial likelihood rather than using the nonparametric or parametric estimators \re{hbbepLink} and \re{hbbep2Link}. Typically the baseline hazard is estimated nonparametrically (see for instance Kalbfleish and Prentice, 2002). This corresponds to a nonparametric estimate
\beq
\hh(y) = \log [\hLa_{\base}(y)]
\lb{hhyCox}
\eeq
of the link function $h$, which is inserted into \re{hthLink2} in order to estimate $\theta(\bx)$. On the other hand, if the baseline distribution $F_{\base}(y)=1-\exp(-\La_{\base}(y))=1-\exp(-ay^b)$ is Weibull, with $a,b>0$ known, then \re{hCox} is equivalent to an AFT model with a known logarithmic link function \re{hlink}, and with $Q_0(p)$ a quantile of a Gumbel distribution in \re{hAFT}. Gelfand et al.\ (2000) proposed a larger parametric model for the baseline distribution $F_{\base}$; a mixture of Weibull distributions. In our context this corresponds to having a finite number of unknown parameters $\bga$ of the link function \re{hCox}, and a parametric estimate \re{hhyCox} of this link function based on $\hLa_{\base}(y)=\La_{\base}(y;\hbga)$. Finally, Royston and Parmar (2002) assumed a version 
\beq
\log[\La(y \mid x_2,\ldots,x_q)] = s(\log(y);\bga) + \sum_{k=2}^q x_k b_{k-1}
\lb{CoxSpline}
\eeq
of \re{Cox2} where $\log[\La_{\base}(y)]$ is replaced by a cubic spline function $s(\cdot;\bga)$ of $\log(y)$ that is parametriz\-ed by $\bga$. It can be seen that this corresponds to replacing the link function in \re{hCox} by $h(y;\bga)=s(\log(y);\bga)$. This gives rise to an estimated link function 
\beq
\hh(y) = s(\log(y);\hbga),
\lb{hhyRP}
\eeq 
which is inserted into \re{hthLink2} in order to estimate $\theta(\bx)$.    
\slut\end{exa}
\par\medskip
\begin{exa}[Log baseline odds transformations.]\lb{Ex:Logodds} \hspace{5pt}
As in the previous\break example, let $x_2,$ $\ldots,x_q$ represent $q-1$ covariates. Bennett (1983) introduced a model for which the proportional odds of death
satisfies 
\beq
\Gamma(y \mid x_2,\ldots,x_q) = \frac{F_{Y \mid x_2,\ldots,x_q}(y)}{1-F_{Y \mid x_2,\ldots,x_q}(y)} 
= \Ga_{\base}(y)\exp(\sum_{k=2}^q x_k b_{k-1}),
\lb{POdds}
\eeq
for some effect parameters $\bb=(b_1,\ldots,b_{q-1})^T$ and baseline odds function $\Ga_{\base}(y)$. Doksum and Gasko (1990) showed that this model can be rewritten as a transformed linear model \re{QpxLinLink}, with link function  
\beq
h(y) = \log\left[\Gamma_{\base}(y)\right] \stackrel{y=Q(p \mid \sbx)}{=} \log [\frac{p}{1-p}] - \sum_{k=2}^q x_k b_{k-1} = \bx^T\bbe(p).
\lb{hPOdds}
\eeq
The vector $\bx=(1,x_2,\ldots,x_q)^T$ contains an intercept and covariates, whereas $\bbe(p)=\log[p/(1-p)]\bfm{e}_1 - (0,\bb^T)^T$ includes the effect parameters $\bb$ of the covariates and an intercept parameter $\log[p/(1-p)]$ that is a quantile of a standard logistic distribution. This corresponds to a parametric model \re{hbbep2Link} with $r=2$, a known $\bpsi_1=\bfm{e}_1$, and unknown $\bpsi_2=(0,\bb^T)^T$, $\beta^1(p)=\log[p/(1-p)]$, and $\beta^2(p)\equiv -1$.
\par\medskip
If the baseline distribution $$F_{\base}(y)=\Ga_{\base}(y)/(1+\Ga_{\base}(y))=\left[1+\exp(-(\log(y)-a)/b)\right]^{-1}$$ is log-logistic, for some known $a$ and $b>0$, it can be seen that the proportional odds model \re{hPOdds} is equivalent to an AFT model \re{hAFT} with a logarithmic link function, where $Q_0(p)$ is the quantile of a logistic distribution. Royston and Parmar (2002) studied a version 
\beq
\log[\Gamma(y \mid x_2,\ldots,x_q)] = s(\log(y);\bga) + \sum_{k=2}^q x_k b_{k-1}
\lb{POddsSplines}
\eeq
of \re{POdds} where the log baseline odds $\log[\Ga_{\base}(y)]$ is replaced by a cubic spline function $s(\cdot;\bga)$ of $\log(y)$. It can be seen that this corresponds to replacing the link function in \re{hPOdds} by $h(y)=s(\log(y);\bga)$. The corresponding estimate \re{hhyRP} is then inserted into \re{hthLink2} in order to estimate $\theta(\bx)$.           
\slut\end{exa}
\par\medskip
\begin{exa}[Log power transformations.]\lb{Ex:LogPower} \hspace{5pt}
Younes and Lachin (1997)\break considered a class of models which includes proportional hazards and proportional odds as special cases. In more detail, they assumed that the logarithm of a Box-Cox transformation \re{hga} of the survival function $1-F_Y$, conditionally on covariates $x_2,\ldots,x_q$, satisfies
$$
\log \frac{(1-F_{Y \mid x_2,\ldots,x_q}(y))^{-\ga}-1}{\ga} = \log \frac{(1-F_{\base}(y))^{-\ga}-1}{\ga} + \sum_{k=2}^q x_k b_{k-1}
$$
for some $\ga>0$, where $\ga=1$ corresponds to the proportional odds model \re{POdds} and $\ga\to 0$ to the proportional hazards model \re{Cox2}. This is a transformed linear model with link function
\beq
\begin{split}
h(y) &= \log\frac{(1-F_{\base}(y))^{-\ga}-1}{\ga} \\ \{y=Q(p \mid \bx)\} &= \log\frac{(1-p)^{-\ga}-1}{\ga} - \sum_{k=2}^q x_k b_{k-1} \\ &= \bx^T\bbe(p),
\end{split}
\lb{hlogPower}
\eeq
covariate vector $\bx=(1,x_2,\ldots,x_q)^T$, and regression parameter vector $\bbe(p)=\log\{[(1-p)^{-\ga}-1]/\ga\}\bfm{e}_1 - (0,\bb^T)^T$. 
\slut\end{exa}
\par\medskip
In order to study the asymptotic properties of the estimator $\hth(\bx)$ in \re{hthLink} of the conditional $L$-functional $\theta(\bx)$ in \re{thetaLink}, we will assume that the link function $h$ is continuously differentiable, with a strictly positive derivative. It is helpful to approximate \re{hthLink} by a first order Taylor expansion
\beq
\hth(\bx) \approx \theta(\bx) + \bx^T\int_0^1 \left[\hbbe(p)-\bbe(p)\right]dG_{\sbx}(p), 
\lb{hthLinkAppr}
\eeq  
where
\beq
dG_{\sbx}(p) = \frac{1}{h^\prime[h^{-1}(\bx^T\bbe(p))]}\cdot dG(p)
\lb{Gxp}
\eeq
can be thought of as an effective weight function, which determines how much different quantiles contribute to the estimation error of $\theta(\bx)$.   
The following result is a corollary of Propositions \ref{Prop:Cons}-\ref{Prop:BnAsN}:
\par\medskip
\begin{cor}[Consistency and asymptotic normality of \re{hthLink}.]\lb{Cor:thetax} \hspace{4pt}
Sup\-pose that the regularity conditions of Propositions \ref{Prop:Cons}-\ref{Prop:BnAsN} hold for the transformed regression model $(\bx,h(Y))$ and that $h$ is a known link function such that $h^\prime(y)\ge c_{\sbx} > 0$ for all $y\in [a_{\sbx},b_{\sbx}]$, where $(a_{\sbx}+\va,b_{\sbx}-\va)$ includes the set $h^{-1}\left(\{\bx^T\bbe(p);\right.$ $p\left.\in\mbox{supp}(G)\}\right)$ for some $\va>0$. The conditional $L$-statistic $\hth_n(\bx)$ in \re{hthLink} is then a consistent estimator of $\theta(\bx)$, so that $\hth_n(\bx)\pto \theta(\bx)$ as $n\to\infty$. It is asymptotically normal as well, i.e. 
\beq
\sqrt{n}[\hth_n(\bx)-\theta(\bx)] \Lto N(0,\bx^T\bSi_{\sbx}\bx)
\lb{AsNLink}
\eeq
as $n\to\infty$, with a covariance matrix
\beq
\bSi_{\sbx} = \int_0^1\int_0^1 \bR(p,s)dG_{\sbx}(p)dG_{\sbx}(s)
\lb{bSix}
\eeq
that involves the effective weight function $G_{\sbx}$ in \re{Gxp} and $\bR(p,s)$, the asymptotic covariance matrix of $\{\hbbe_n(p);\, 0<p<1\}$. 
\end{cor}
\par\medskip
{\bf Proof.} The result can be derived similarly as in the proofs of Propositions \ref{Prop:Cons}-\ref{Prop:BnAsN}, using the Taylor expansion \re{hthLinkAppr}. Indeed, the regularity conditions on $h^\prime$ imply that the remainder term of this Taylor expansion is asymptotically negligible, both for the consistency and the asymptotic normality proofs.    
\hfill\slut 
\par\medskip
When our objective is to estimate conditional skewness or kurtosis we focus on quantities $\theta(\bx)$, defined as the ratio \re{thxRatio0} of two conditional $L$-functionals $T^1(F_{Y \mid \sbx})$ and $T^2(F_{Y \mid \sbx})$ with different weight functions $G^1$ and $G^2$. For the transformed linear model \re{QpxLinLink} we find that 
\beq
\theta(\bx) = \frac{\int_0^1 h^{-1}(\bx^T\bbe(p))dG^1(p)}{\int_0^1 h^{-1}(\bx^T\bbe(p))dG^2(p)}.
\lb{thetaLinkRatio}
\eeq
We estimate $\theta(\bx)$ by plugging an appropriate estimator of $\bbe(p)$ into \re{thetaLink}, i.e.
\beq
\hth(\bx) = \frac{\int_0^1 h^{-1}(\bx^T\hbbe(p))dG^1(p)}{\int_0^1 h^{-1}(\bx^T\hbbe(p))dG^2(p)}.
\lb{hthLinkRatio}
\eeq
When the link function $h$ is unknown, as in \re{hthLink2} it is possible to define a version of $\hth(\bx)$ where $h$ is replaced by a parametric or nonparametric estimate $\hh$ in the numerator and denominator of \re{hthLinkRatio}. 
\par\medskip 
The following result is essentially a consequence of Proposition \ref{Prop:JointBnAsN} and Corollary \ref{Cor:thetax}:
\par\medskip
\begin{cor}[Consistency and asymptotic normality of \re{hthLinkRatio}.]\lb{Cor:Ratio} \hspace{4pt}
Sup\-pose that the regularity conditions of Proposition \ref{Prop:JointBnAsN} and Corollary \ref{Cor:thetax} hold. The quantity $\hth_n(\bx)$ in \re{hthLinkRatio} is then a consistent estimator of $\theta(\bx)$ in \re{thetaLinkRatio}. It is also asymptotically normal, in the sense that   
\beq
\begin{split}
\sqrt{n}[\hth_n(\bx)-&\theta(\bx)] \Lto N\left(0,\Sigma^\star\right) \\
\Sigma^\star &= \frac{\bx^T\bSi^{11}_{\sbx}\bx}{T^2(F_{Y \mid \sbx})} - \frac{2T^1(F_{Y \mid \sbx})\cdot \bx^T\bSi_{\sbx}^{12}\bx}{(T^2(F_{Y \mid \sbx}))^3} + \frac{(T^1(F_{Y \mid \sbx}))^2\cdot \bx^T\bSi^{22}_{\sbx}\bx }{(T^2(F_{Y \mid \sbx}))^4}
\end{split}
\lb{hthnAsNLinkRatio}
\eeq
as $n\to\infty$, where
\beq
\bSi^{kl}_{\sbx} = \int_0^1\int_0^1 \bR(p,s)dG^k_{\sbx}(p)dG^l_{\sbx}(s),
\lb{Sikl}
\eeq
$G^1_{\sbx}$ and $G^2_{\sbx}$ are defined as in \re{Gxp}, with $G^1$ and $G^2$ in place of $G$, and $\bR(p,s)$ is the asymptotic covariance matrix of $\{\hbbe_n(p);\, 0<p<1\}$.  
\end{cor}
\par\medskip
The consistency and asymptotic normality of $\hth_n(\bx)$ in Corollaries \ref{Cor:thetax}-\ref{Cor:Ratio} is a unified result for nonparametric or parametric estimates \re{hbbepLink}-\re{hbbep2Link} of $\bbe=\{\bbe(p);\, 0<p<1\}$. It is only required that the link function $h$ is known, and that $\hbbe=\{\hbbe(p);\, 0<p<1\}$ is a consistent and asymptotically normal estimator of $\bbe$, whose covariance function $\bR(p,s)$ appears in \re{bSix} and \re{Sikl}. In the nonparametric case, $\bR(p,s)$ is defined as in \re{R}, provided the density in the expression for $\bD_n(p)$ is changed to $f_i(p) = \left. dF_{h(Y) \mid \sbx_i}(v)/dv\right|_{v=h(Q(p \mid \sbx_i))}$. In the parametric case, $\bR(p,s)$ is defined as in \re{RpsPar}. See also Newey and McFadden (1994), Ying et al.\ (1995), Chen et al.\ (2003), Mu and He (2007), and references therein, for a discussion on when the estimator \re{hthLink2} of $\theta(\bx)$ with estimated link function is asymptotically equivalent to the corresponding estimator \re{hthLink} where $h$ is known.

\subsection{Censored and truncated data}

Assume there exists a collection $\{(\bx_i,Y_i)\}_{i=1}^n$ of i.i.d.\ random vectors, and that $F_{Y \mid \sbx}$ follows the transformed linear model \re{QpxLinLink}. The objective is to estimate the $L$-functional \re{thetaLink} or the ratio \re{thetaLinkRatio} of $L$-functionals, when some data is lost due to censoring or truncation. This boils down to finding an estimator $\hbbe(p)$ of $\bbe(p)$ for all $0<p<1$, and then plugging this estimator into \re{hthLink}, \re{hthLink2} or \re{hthLinkRatio}. Notice that Corollaries \ref{Cor:thetax}-\ref{Cor:Ratio} apply to censored and truncated data as well, if consistency and asymptotic normality is established for $\{\hbbe(p);\ 0<p<1\}$, with some limiting covariance function $\{\bR(p,s); \, 0<p,s<1\}$.   

\subsubsection{Censoring}

When a distorted version $\tY$ of the outcome variable $Y$ is observed, it is often the case that a censoring variable $C$ causes this distortion. The two most common types of censoring are 
$$
\tY = \left\{\begin{array}{ll}
\max(Y,C), & \mbox{left-censoring},\\
\min(Y,C), & \mbox{right-censoring}.
\end{array}\right.
$$
Without loss of generality we restrict ourselves to right-censoring. To this end, assume that $\{(\bx_i,Y_i,C_i)\}_{i=1}^n$ are i.i.d.\ copies of covariate vectors, response and censoring variables, such that $Y_i$ and $C_i$ are conditionally independent given $\bx_i$, and with $\tY_i=\min(Y_i,C_i)$ the right-censored version of $Y_i$. 
\par\medskip
Let us first assume that all censoring variables $C_i$ are observed, whether $Y_i$ is censored or not. This corresponds to having a data set consisting of the observations $(\bx_i,\tY_i,C_i)$ for $i=1,\ldots,n$. In this case the conditional quantile
\beq
\tQ(p \mid \bx,C) = F_{\tY \mid \sbx}^{-1}(p) = \min(\bx^T\bbe(p),C)
\lb{tQ}
\eeq
of censored data is an explicit function of $\bbe(p)$. This gives rise to the estimate 
$$
\hbbe(p) = \mbox{arg}\min_{\sbb} \sum_{i=1}^n \rho_p[h(\tY_i)-\min(\bx_i^T\bb,C_i)]
$$
of $\bbe(p)$ due to Powell (1986). 
\par\medskip
For the general right-censoring problem only $\{(\bx_i,\tY_i,\De_i)\}_{i=1}^n$ is observed,\break where $\De_i=1(Y_i=\tY_i)$ indicates whether observation $i$ has been censored or not. The relation between the quantile functions of censored and uncensored response variables is then somewhat more complicated than \re{tQ}. Because of conditional independence of $Y$ and $C$ given $\bx$, it follows that  
\beq
F_{\tY \mid \sbx}(y) = 1 - (1-F_{Y \mid \sbx}(y))(1-F_{C \mid \sbx}(y)),
\lb{FtYFC}
\eeq
and consequently the quantile function $\tQ(\cdot \mid \bx)$ of the right censored response variable is related to the quantile function $Q(\cdot \mid \bx)$ of the uncensored response as    
\beq
\tQ(\pi(\bx,p) \mid \bx) = Q(p \mid \bx),
\lb{tQQ}
\eeq
where $\pi(\bx,p) = 1-(1-p)(1-F_{C \mid \sbx}(Q(p \mid \bx)))$ tells which quantile of the censored observation an uncensored $p$-quantile relates to. If $F_{C \mid \sbx}$ would be known for all $\bx$, then based on \re{tQQ} and the fact that quantiles are preserved under the monotone transformation $h$, the estimator
\beq
\hbbe(p) = \mbox{arg}\min_{\sbb} \sum_{i=1}^n \rho_{\pi(\sbx_i,p)}(h(\tY_i)-\bx_i^T\bb)
\lb{hbbepCens}
\eeq
of $\bbe(p)$ due to Lindgren (1997) could be used. When the censoring distribution is unknown, a nonparametric Kaplan-Meier estimator $\hF_{C \mid \sbx}$ was employed by Lindgren (1997) in order to replace $\pi(\bx_i,p)$ in \re{hbbep} by an estimator $\hpi(\bx_i,p)$. 
\par\medskip
Another consequence of \re{FtYFC} is that the random variable
$$
\frac{1(h(\tY)-\bx^T\bbe(p)\ge 0)}{1-F_{C \mid \sbx}(\bx^T\bbe(p))} - (1-p)
$$
has zero expectation. 
This motivated Ying et al.\ (1995) and Leng and Tong (2013) to propose and study the score-based estimator
\beq
\hbbe(p) = \mbox{arg}\min_{\sbb} \left| \sum_{i=1}^n \bx_i \left[\frac{1-\tom_i(\bb)}{1-\hF_{C \mid \sbx_i}(\bx_i^T\bb)} - (1-p)\right] \right|  
\lb{hbbepCens2}
\eeq
of $\bbe(p)$, with $\tom_i(\bb)=1(h(\tY_i)\le \bx_i^T\bb)$. Leng and Tong (2013) proved that this estimator is asymptotically normal with a $\sqrt{n}$-rate of convergence. As a drawback, \re{hbbepCens2} does not account for which observations that have been censored ($\De_i=0$) or not. In order to use this information, Wang and Wang (2009) introduced
\beq
\hbbe(p) = \mbox{arg}\min_{\sbb} \left| \sum_{i=1}^n \bx_i \left[p-\tom_i(\bb) + \frac{(1-\De_i)\tom_i(\bb)(1-p)}{1-\hF_{Y \mid \sbx_i}(\tY_i)}\right] \right|,
\lb{hbbepCens3}
\eeq
where $\hF_{Y \mid \sbx}$ is a local Kaplan-Meier estimate of $F_{Y \mid \sbx}$. This estimator can be motivated by noticing that each term of \re{hbbepCens3} has zero expectation when $\bb=\bbe(p)$, and $\hF_{Y \mid \sbx_i}$ is replaced by the true but unknown $F_{Y \mid \sbx_i}$. Other estimators of $\bbe(p)$, for censored observations, have been proposed by Yang (1999), Portnoy (2003), Neocleous et al.\ (2004), and Peng and Huang (2008).   

\subsubsection{Truncation}

Truncation means that some observations $(\bx,Y)$ are lost, depending on the value of some other truncation random variable $L$. The two most common types of truncation are that 
$$
(\bx,Y) \mbox{ is lost when }\left\{\begin{array}{ll}
Y<L, & \mbox{left-truncation},\\
Y>L, & \mbox{right-truncation}.
\end{array}\right.
$$
We will assume that left-truncation occurs together with right-censoring, so that the observed data set is $\{(\bx_i,\tY_i,\De_i,L_i);\, i=1,\ldots,n, Y_i>L_i\}$.
Frumento and Bottai (2017) generalized \re{hbbepCens3} and presented an estimator 
\beq
\begin{split}
\hbbe(p) = \mbox{arg}\min_{\sbb} \left| \sum_{i=1}^n \bx_i \left[\om_i(\bb) -\frac{\om_i(\bb)(1-p)}{1-\hF_{Y \mid \sbx_i}(L_i)}\right.\right. \hspace{74pt} \\ 
\hspace{74pt} \left.\left. - \tom_i(\bb) +  \frac{(1-\De_i)\tom_i(\bb)(1-p)}{1-\hF_{Y \mid \sbx_i}(\tY_i)}\right] \right|,
\end{split}
\lb{hbbepTrunc}
\eeq
of $\bbe(p)$, where $\om_i(\bb)=1(h(L_i)\le \bx_i^T\bb)$. Frumento and Bottai (2017) gave conditions under which this estimator of $\bbe(p)$ is consistent and asymptotically normal. 

\section{Numerical examples}\lb{Sec:Sim}


In this section we will analyze numerical properties of $L$-functionals without covariates, as described in Section \ref{sec:Lstat}. 
In more detail we consider the four (standardized) $L$-functionals $T_1(F)$, $T_2(F)$, $T_{32}(F)$, and $T_{42}(F)$ of a target distribution $F$. First, in Section \ref{Sec:QuantAppr}, we quantify how well these functionals approximate selected, well known distributions' quantile functions using the polynomial series of Examples \ref{LfuncLegendre}-\ref{Sec:AsymL}. For the Legendre system \re{gLegendre}, this has previously been done by Hosking (1990, 1992), Karvanen (2006, 2008) and Elamir and Seheult (2003). We will find that typically, matching the ``essential support'' (or form) of a linearly standardized version of the target distribution $F$ (i.e.\ a region harbouring most of the probability mass of a standardized version of $F$) with the support (or form) of the polynomial system's reference distribution $F_0$, produces the most accurate approximations. In Section \ref{Sec:Tgraphs}, we present graphs showing the change in the $L$-functionals $T_2$, $T_{32}$ and $T_{42}$ for Beta distributions $F\sim B(\psi_1,\psi_2)$, when the two shape parameters $\psi_1$ and $\psi_2$ are varied in such a way that the expected value $T_1(F)=\psi_1/(\psi_1+\psi_2)$ in \re{TEY} is held constant. All computations were performed in \textsf{R} (R Core Team, 2021).

\subsection{Quantile function approximation errors}\lb{Sec:QuantAppr}

In this subsection we will investigate the suitability of the three orthogonal polynomial series of Examples \ref{LfuncLegendre}-\ref{Sec:AsymL} for approximating the quantile function $Q(p)=F^{-1}(p)$, using $m_0$ terms. For instance, $m_0=4$ terms is the approximation using only location, scale, skewness and kurtosis. This will be done for distributions $F$ with various types of essential support (that is, a set that supports close to 1 of the probability mass of $F$). The three types of essential support are interval support $[a,b]$, half-infinite support of life-time distributions $[0, \infty)$ and doubly infinite support on the real line $\R$. As a rule of thumb we hypothesize that the polynomial series whose reference distribution $F_0$ has a support of the same type as the essential support of the distribution $F$ whose quantile function we try to approximate, will be the most suitable. Thus, the Legendre polynomials should suit distributions with bounded essential support, Laguerre polynomials should be appropriate for distributions with half infinite essential support, whereas Hermite polynomials are preferable for distributions with unbounded essential support, to the left and right.

The approximated quantile function is denoted $Q_{\scr{appr}}(p) = \sum_{m=1}^{m_0} T_m(F)g_m(p)$. 
Note that, because of the orthogonality property \re{gOrth}-\re{QExp}, the Integrated Squared Error (ISE) of the approximation error $Q-Q_{\scr{appr}}$, satisfies  
\beq
\begin{split}
 \int_0^1 \left( Q(p) - Q_{\scr{appr}}(p) \right)^2 \rmd p &= \int_0^1 \left( \sum_{m=1}^\infty T_m(F)g_m(p) - \sum_{m=1}^{m_0} T_m(F)g_m(p) \right)^2 \rmd p \\
 &= \int_0^1 \left( \sum_{m=m_0+1}^\infty T_m(F)g_m(p)\right)^2 \rmd p,\\
&= \sum_{m=m_0+1}^\infty T_m^2(F),
\end{split}
\lb{QAppr}
\eeq
meaning that the approximation error amounts to what cannot be summarized about the distribution from the first $m_0$ $L$-moments. To compute the ISE, we evaluate  
\begin{align} \label{RISE}
\int_0^1 \left( Q(p) - \sum_{m=1}^{m_0} \left( \int_0^1 Q(\rho) g_m(\rho) \rmd \rho \right) g_m(p) \right)^2 \rmd p
\end{align}
numerically. This makes it possible to compute how large fraction 
\beq
\Delta_{m_0} = 1 - \frac{\mbox{ISE}}{ \int_0^1 Q^2(p) \rmd p }
\lb{Delta}
\eeq
of the variation of $Y$ that is explained by the first $m_0$ $L$-moments. In Table \ref{approxerror} we have computed $\Delta_4$ for selected distributions, using the abovementioned three polynomial series.

Sometimes the numerical procedures fail to evaluate either integral of \eqref{RISE} when $F$ has heavy tails. Most commonly, $m=4$ causes problems for the Hermite polynomials, as can be seen from \re{gHermite}, \re{gHermite2}, and \eqref{RISE}, and less commonly for the higher order $L$-functionals using the Laguerre polynomials (cf. \re{gExp} and \eqref{RISE}). Shortening the inner and outer integration intervals of \eqref{RISE} to $(\eps, 1-\eps)$ would mitigate these problems. Therefore, in order to reduce the numerical problems associated with the approximation error \eqref{RISE} we construct a new series $g_m^\eps$ of orthonormal polynomials to replace $g_m$ in $Q_{\scr{appr}}$. We do this by using the Gram-Schmidt process to orthonormalize the base functions $g_1(p),\ldots,g_4(p)$ with an inner product
\begin{align}
 \int_\eps^{1-\eps} g_k(p)g_l(p)\rmd p, \quad k,l \in \{1,\ldots,4\}.
\end{align}
Thus, we can approximate \eqref{RISE} using $g_m^\eps(p)$ for any chosen $\eps$ and distribution. In particular, $\eps = 0$ corresponds to using $g_m^0=g_m$ in $Q_{\scr{appr}}$. Notice that $g_m^\eps$ is conceptually different from $g_m^\pi$ used in the robustified  functionals in Examples \ref{LfuncLegendre}, \ref{Sec:Hermite} and \ref{Sec:AsymL}. Since the numerical issues occur when $p$ approaches 0 and 1, we cannot mitigate it with $g_m^\pi$, since it still utilizes the whole support of the chosen polynomial series.

\begin{table}
\centering
\footnotesize
\def\arraystretch{1.3}
 \begin{tabular}{|c||c|c|c|c|c|c|}
\hline
& \multicolumn{6}{|c|}{Polynomial system} \\
\cline{2-7}
& \multicolumn{2}{|c|}{Legendre} & \multicolumn{2}{|c|}{Hermite} & \multicolumn{2}{|c|}{Laguerre} \\ \cline{2-7}
Distribution $F$ & $\Delta_4$ & $\epsilon$ & $\Delta_4$ & $\epsilon$ & $\Delta_4$ & $\epsilon$ \\
\hline
$\U(0,1)$ & {\bf 100\%} & 0 & 99.87\% & 0 & 99.61\% & 0 \\
$\beta(0.1,0.1)$ & {\bf 98.23\%} & 0 & 93.27\% & 0 & 93.96\% & 0 \\
$\text{N}(0,1)$ & 98.84\% & 0 & {\bf 100\%} & 0 & 94.31\% & 0 \\
$t_{10}$ & 97.36\% & 0 & {\bf 99.998\%} & $10^{-5}$ & 92.19\% & 0 \\
$\Exp(1)$ & 96.88\% & 0 & 99.99\% & $10^{-5}$ & {\bf100\%} & 0 \\
$\Exp(1)$ & 96.87\% & $10^{-5}$ & 99.99\% & $10^{-5}$ & {\bf100\%} & $10^{-5}$ \\
$\Exp(10)$ & 96.88\% & 0 & 99.99\% & $10^{-5}$ & {\bf 100\%} & 0 \\
$\Ga(10,1)$ & 99.84\% & 0 & {\bf 100\%} & 0 &  99.80\% & $10^{-7}$ \\
$\Wei(3,1)$ & 99.91\% & 0 & {\bf 99.999\%} & $10^{-5}$ & 99.57\% & $10^{-6}$ \\
$\Wei(1/2, 1)$ & 73.78\% & 0 & 97.95\% & 0 & {\bf100\%} &0 \\
\hline
 \end{tabular}
\caption{Values of how large a fraction ($\Delta_4$) of the quantile functions $Q$ that is explained by the first 4 terms of a polynomial series expansion, for selected distributions $F$, according to \eqref{Delta}. The first column denotes the distribution, the next three specify which polynomial series was used to approximate the quantile function, and the subcolumns specify the values of $\Delta_4$ and $\eps$ for each approximation. The best approximation is highlighted. $\Exp(\psi)$ refers to an exponential distribution with scale parameter $\psi$, so that $E_F(Y)=\psi$, $\Gamma(\psi_1,\psi_2)$ is a gamma distribution with shape parameter $\psi_1$ and scale parameter $\psi_2$, so that $E_F(Y)=\psi_1\psi_2$ and $\Gamma(1,\psi)=\Exp(\psi)$. $\Wei(\psi_1,\psi_2)$ corresponds to a Weibull distribution with shape parameter $\psi_1$ and scale parameter $\psi_2$, so that $\Wei(1,\psi)=\Exp(\psi)$. Notice how it is not the true support of the distribution that reveals which polynomial series that best approximates the quantile function, but rather the essential support. For instance, a $\Wei(3,1)$ distribution has support on $[0,\infty)$, but the Hermite polynomials give the best approximation since the left {\sl and} right tails are light, similarly to a normal distribution.}
 \label{approxerror}
\end{table}

\begin{figure} 
\begin{subfigure}{.5\textwidth}
  \centering
  \includegraphics[width=\linewidth]{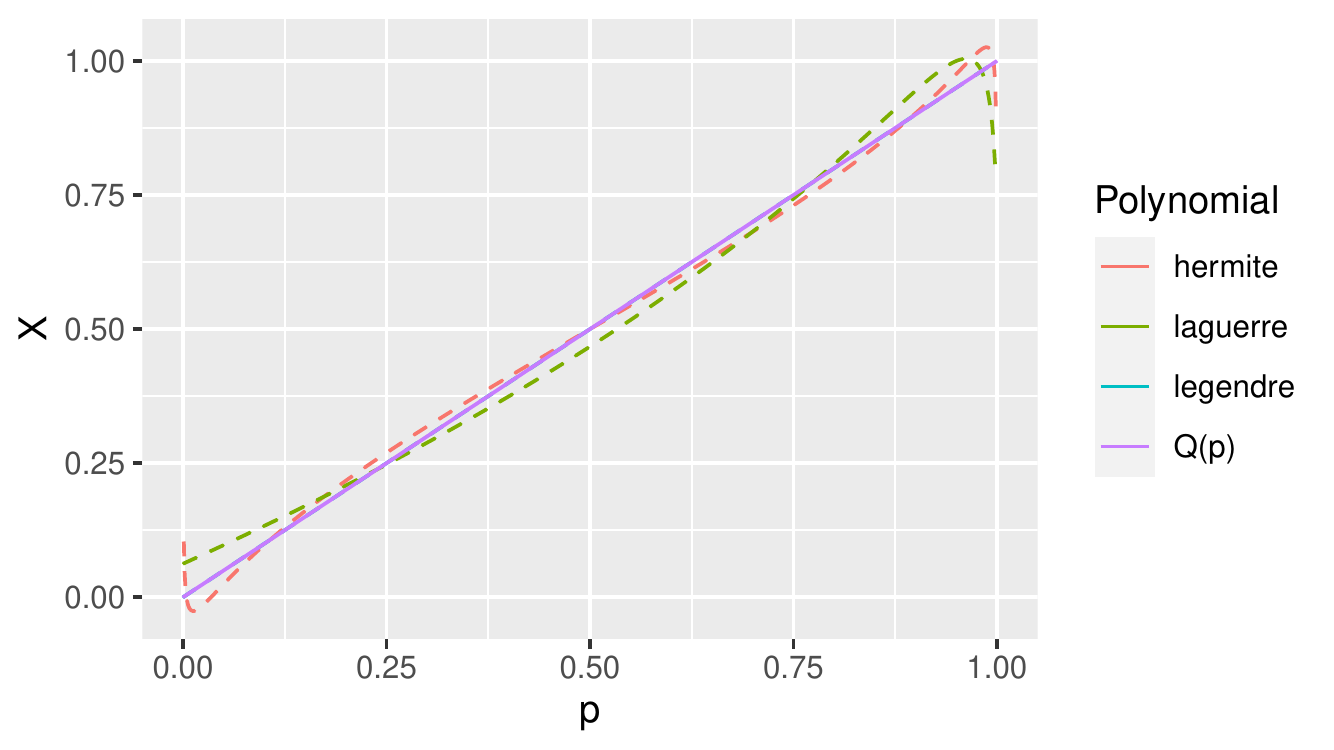}
  \caption{$X \sim \U(0,1)$ with $\eps = 0$.}
  \label{fig:sfig1}
\end{subfigure}%
\begin{subfigure}{.5\textwidth}
  \centering
   \includegraphics[width=\linewidth]{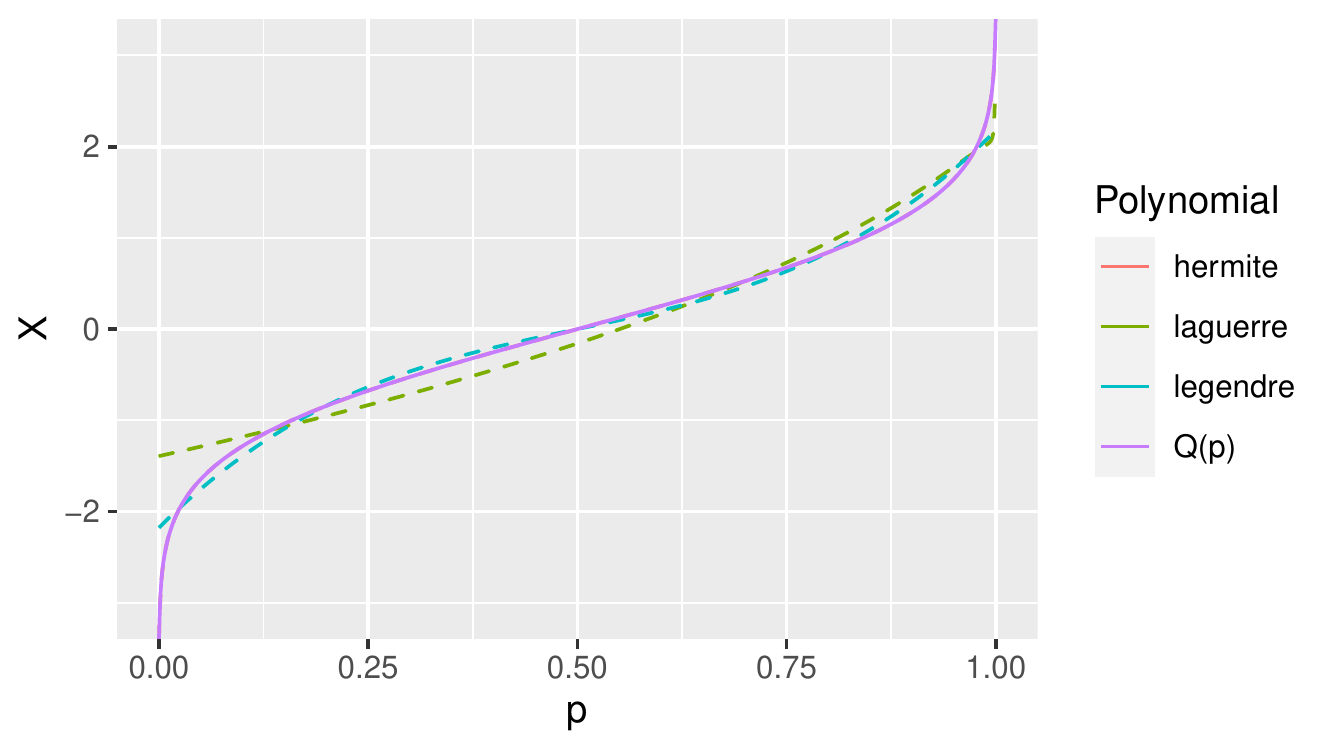}
  \caption{$X \sim\text{N}(0,1)$ with $\eps = 0$.}
  \label{fig:sfig2}
\end{subfigure}
\begin{subfigure}{.5\textwidth}
  \centering
   \includegraphics[width=\linewidth]{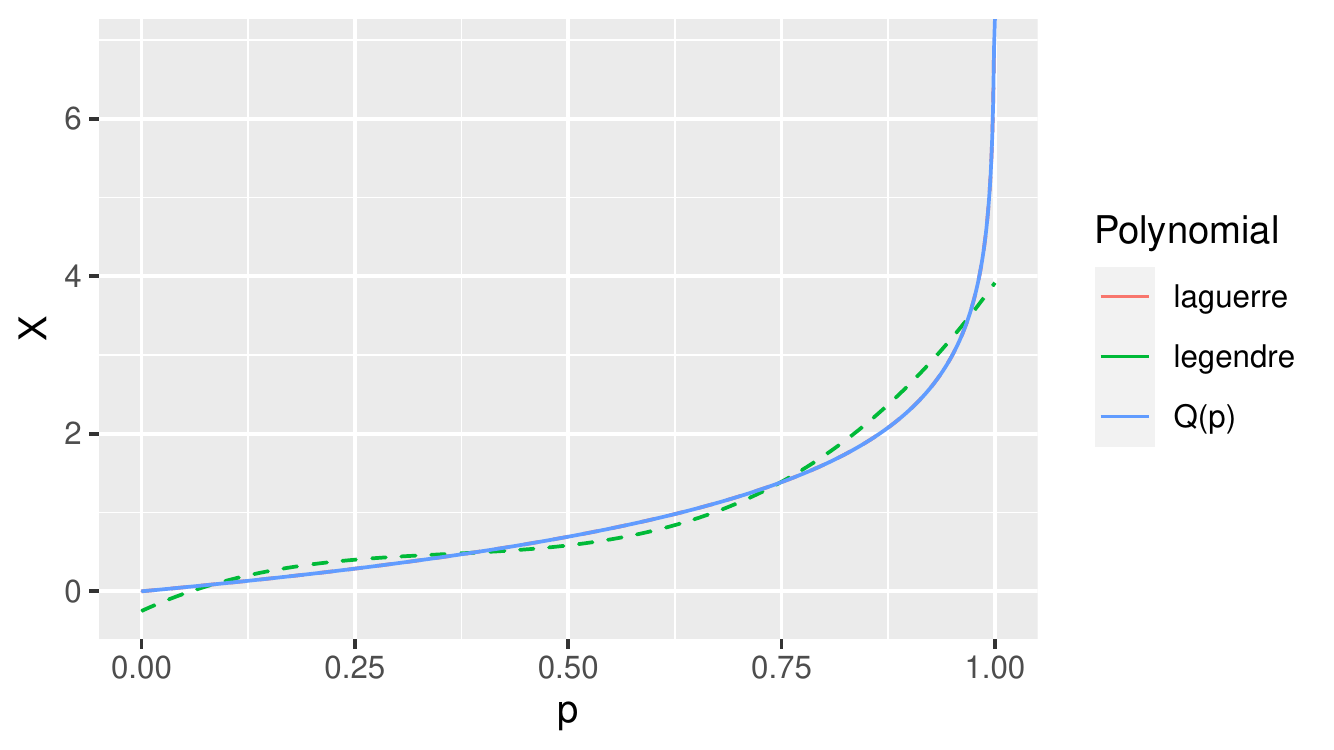}
  \caption{$X \sim\Exp(1)$ with $\eps = 0$.}
  \label{fig:sfig3}
\end{subfigure}%
\begin{subfigure}{.5\textwidth}
  \centering
   \includegraphics[width=\linewidth]{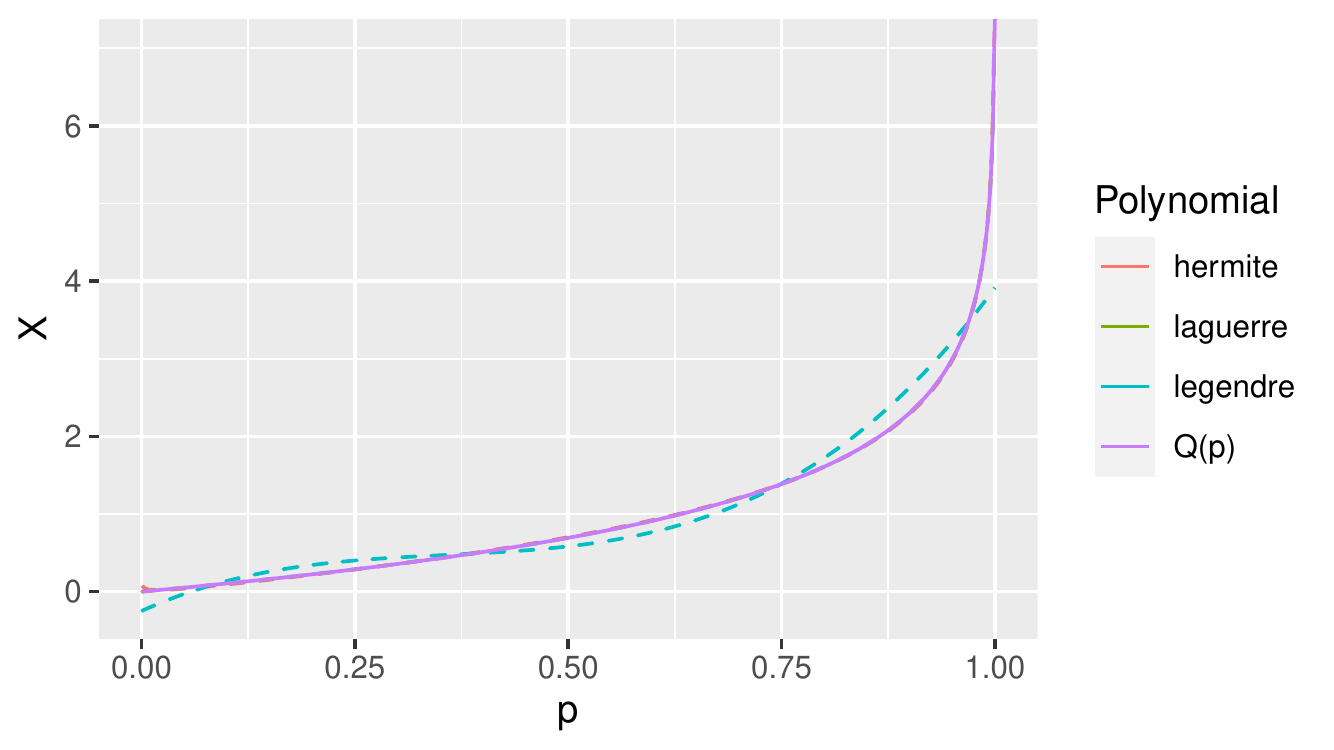}
  \caption{$X \sim\Exp(1)$ with $\eps = 10^{-5}$.}
  \label{fig:sfig4}
\end{subfigure}
\caption{Quantile function $Q(p)$, and the associated approximations $Q_{\scr{appr}}(p)$, for any of the polynomial systems of Section \ref{Sec:LOrder} that is computable, for each chosen $\eps$. The most suitable polynomial series for approximating $Q$, measured by $\Delta_4$, is the Legendre system for (a), the Hermite system for (b) and the Laguerre system for (c) and (d). For (a) to (c), the reference distribution $Q_0$ of the best performing polynomial series equals $Q$, see Examples \ref{LfuncLegendre}-\ref{Sec:AsymL}. For (d), notice how the quantile function $Q$ three times crosses the approximating $Q_{\scr{appr}}$ when using Laguerre polynomials.}
\label{approxlines}
\end{figure}

In Figure \ref{approxlines} we showcase four distributional approximations, where (a)-(c) are approximations using $\eps = 0$ and (d) is the same approximation as in (c), but using $\eps = 10^{-5}$. A comparison between c) and d) reveals the effects of extreme quantiles on $Q_{\scr{appr}}$.

\clearpage

\subsection{Plotting the change in scale, shape and location for beta distributions with fixed location} \lb{Sec:Tgraphs}

\begin{figure} 
\begin{subfigure}{.5\textwidth}
  \centering
  \includegraphics[width=\linewidth]{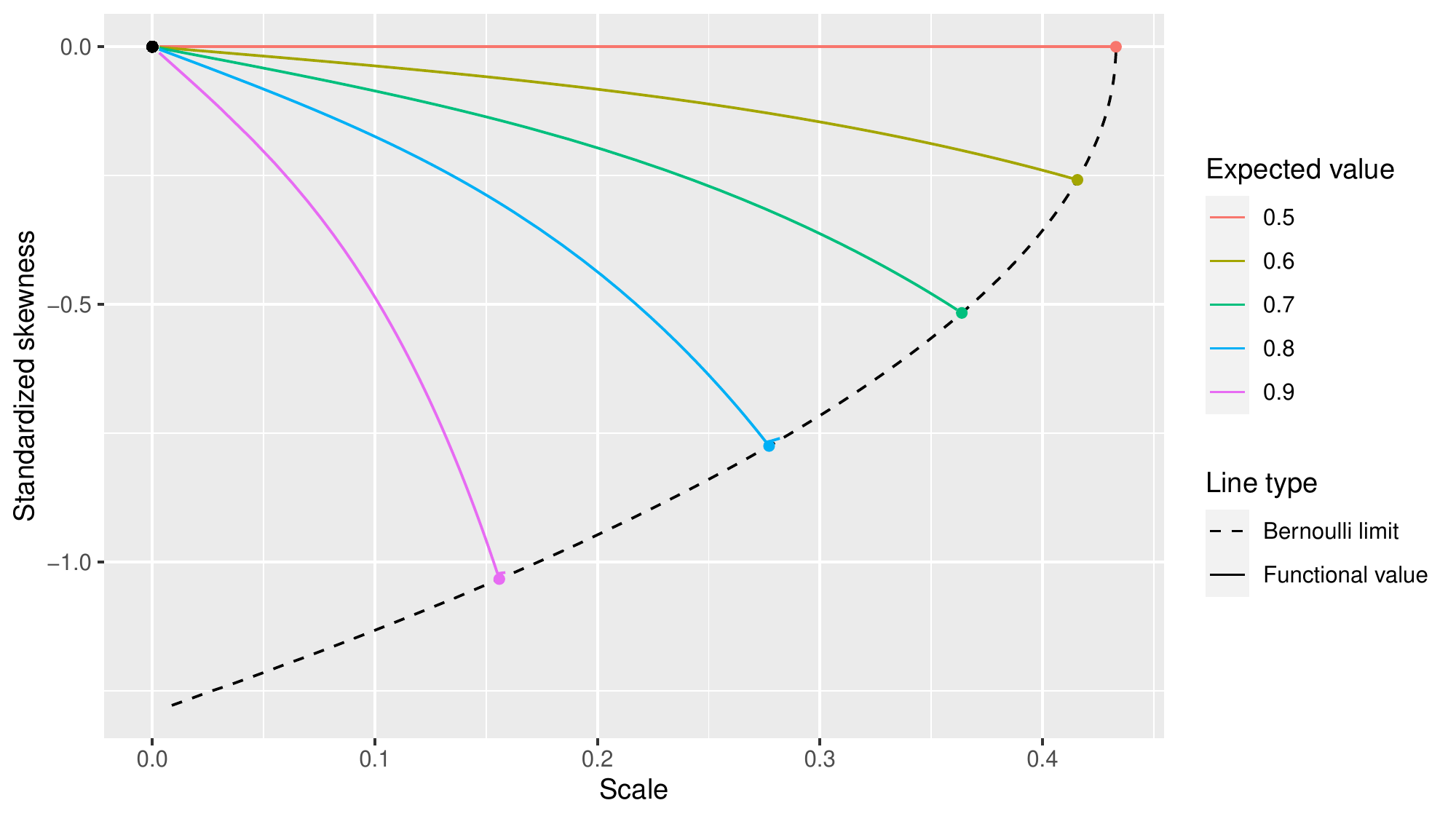}
  \caption{Standardized skewness.}
  \label{fig:m3legbeta}
\end{subfigure}%
\begin{subfigure}{.5\textwidth}
  \centering
   \includegraphics[width=\linewidth]{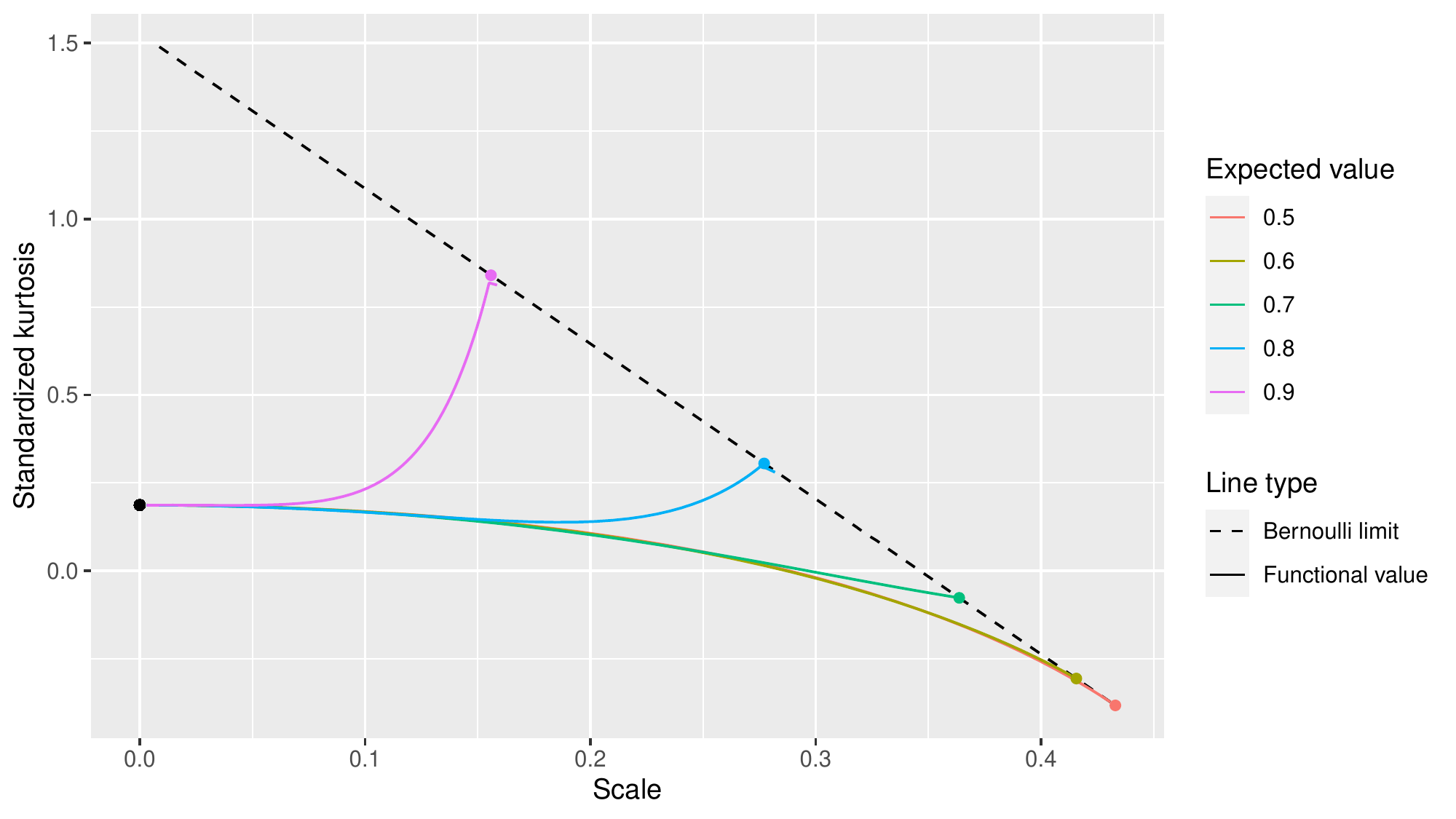}
  \caption{Standardized kurtosis.}
  \label{fig:m4legbeta}
\end{subfigure}
\caption{Plots of standardized skewness $T_{32}(F)$ (a) and standardized kurtosis $T_{42}(F)$ (b) versus scale $T_2(F)$ for various beta distributions $F\sim B(\psi_1,\psi_2)$ using the Legendre polynomials and $\eps = 0$. The black dots correspond to values of $T_{32}(\text{N}(\cdot,\cdot))$ and $T_{42}(\text{N}(\cdot,\cdot))$. The colored dots illustrate the values of $T_{32}(\Be(T_1(F)))$ and $T_{42}(\Be(T_1(F)))$ for each of the five fixed values of $T_1(F)$. Notice that close to the Bernoulli limit, because of the numerical approximations involved in computing the $L$-functionals, both $T_{32}$ and $T_{42}$ show slight deviations from the ideal value for $T_1(F) \in \{0.8,0.9\}$.}
\label{legbeta}
\end{figure}

\begin{figure} 
\begin{subfigure}{.5\textwidth}
  \centering
  \includegraphics[width=\linewidth]{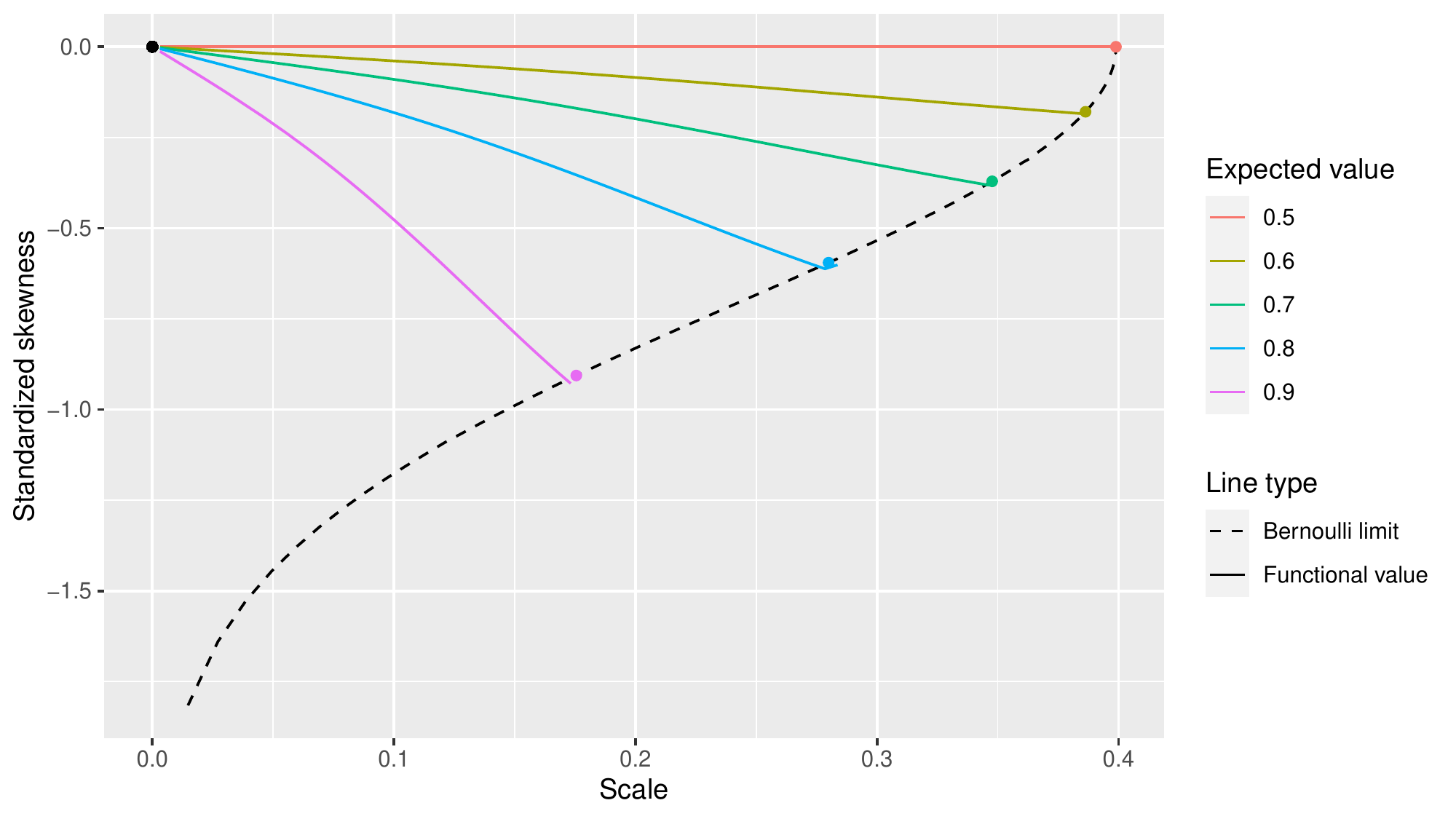}
  \caption{Standardized skewness.}
  \label{fig:m3herbeta}
\end{subfigure}%
\begin{subfigure}{.5\textwidth}
  \centering
   \includegraphics[width=\linewidth]{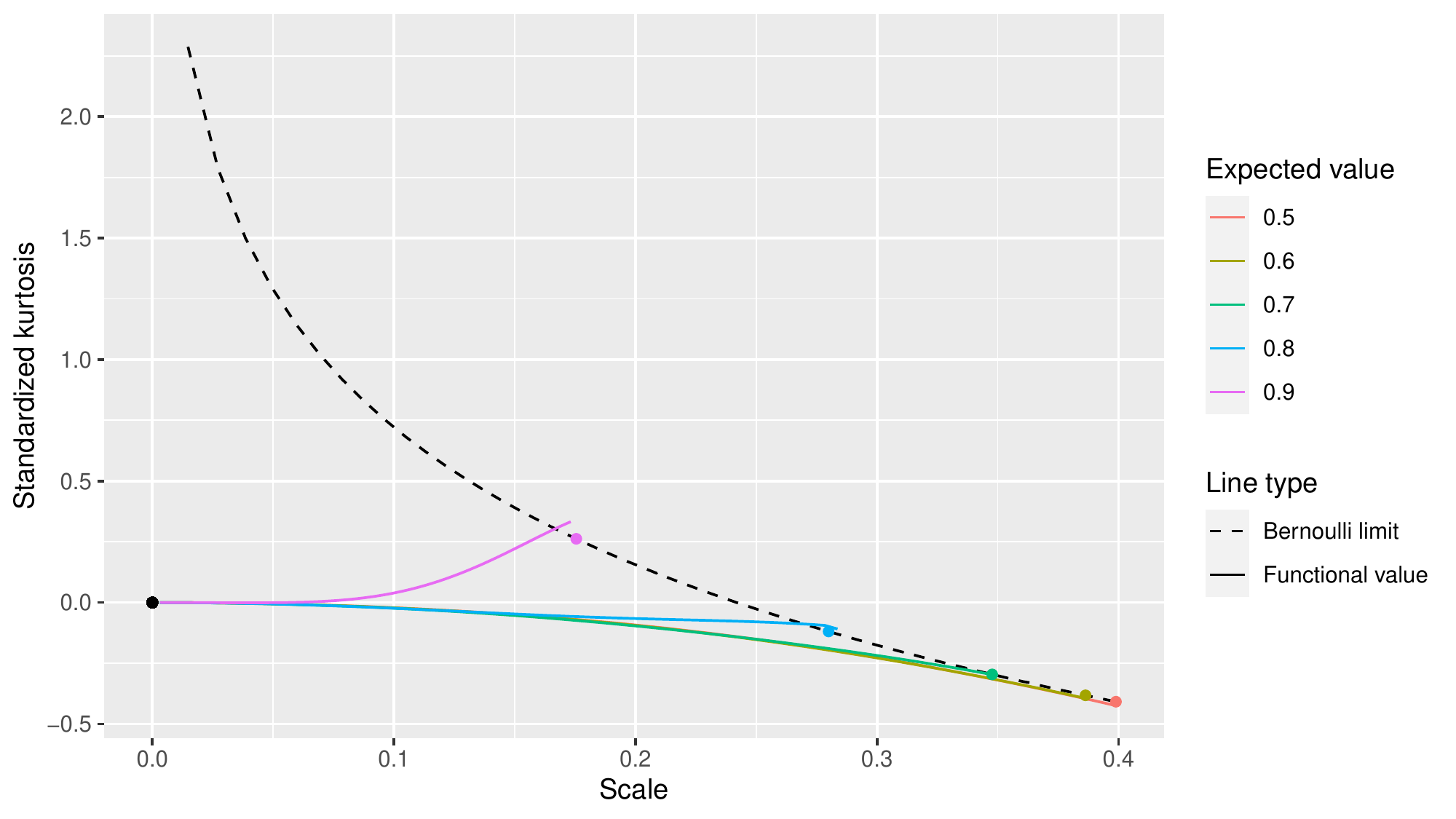}
  \caption{Standardized kurtosis.}
  \label{fig:m4herbeta}
\end{subfigure}
\caption{Plots of standardized skewness $T_{32}(F)$ (a) and standardized kurtosis $T_{42}(F)$ (b) versus scale $T_2(F)$ for various beta distributions $F\sim B(\psi_1,\psi_2)$ using the Hermite polynomials and $\eps = 10^{-3}$. The black dots correspond to values of $T_{32}(\text{N}(\cdot,\cdot))$ and $T_{42}(\text{N}(\cdot,\cdot))$. The colored dots illustrate the values of $T_{32}(\Be(T_1(F)))$ and $T_{42}(\Be(T_1(F)))$ for each of the five fixed values of $T_1(F)$. Notice that the deviations from the expected Bernoulli limits increase as the expected value increases. This is due to an increasing probability mass in the extreme tails, outside of $(\eps,1-\eps)$, which is missed by our adapted orthogonal weight functions $g_m^\eps$, to a higher extent the larger $T_1(F)$ is. }
\label{herbeta}
\end{figure}

\begin{figure} 
\begin{subfigure}{.5\textwidth}
  \centering
  \includegraphics[width=\linewidth]{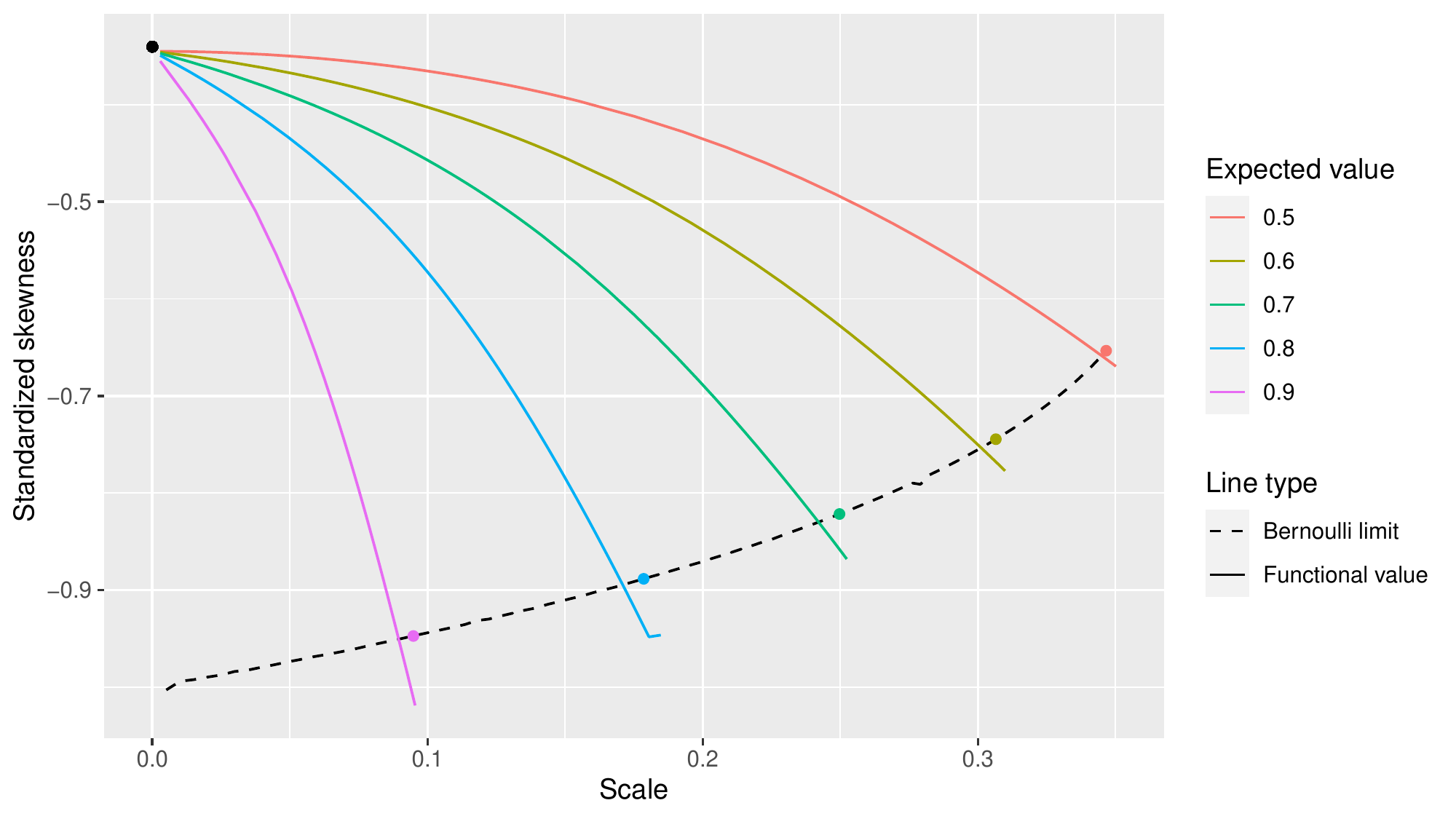}
  \caption{Standardized skewness.}
  \label{fig:m3lagbeta}
\end{subfigure}%
\begin{subfigure}{.5\textwidth}
  \centering
   \includegraphics[width=\linewidth]{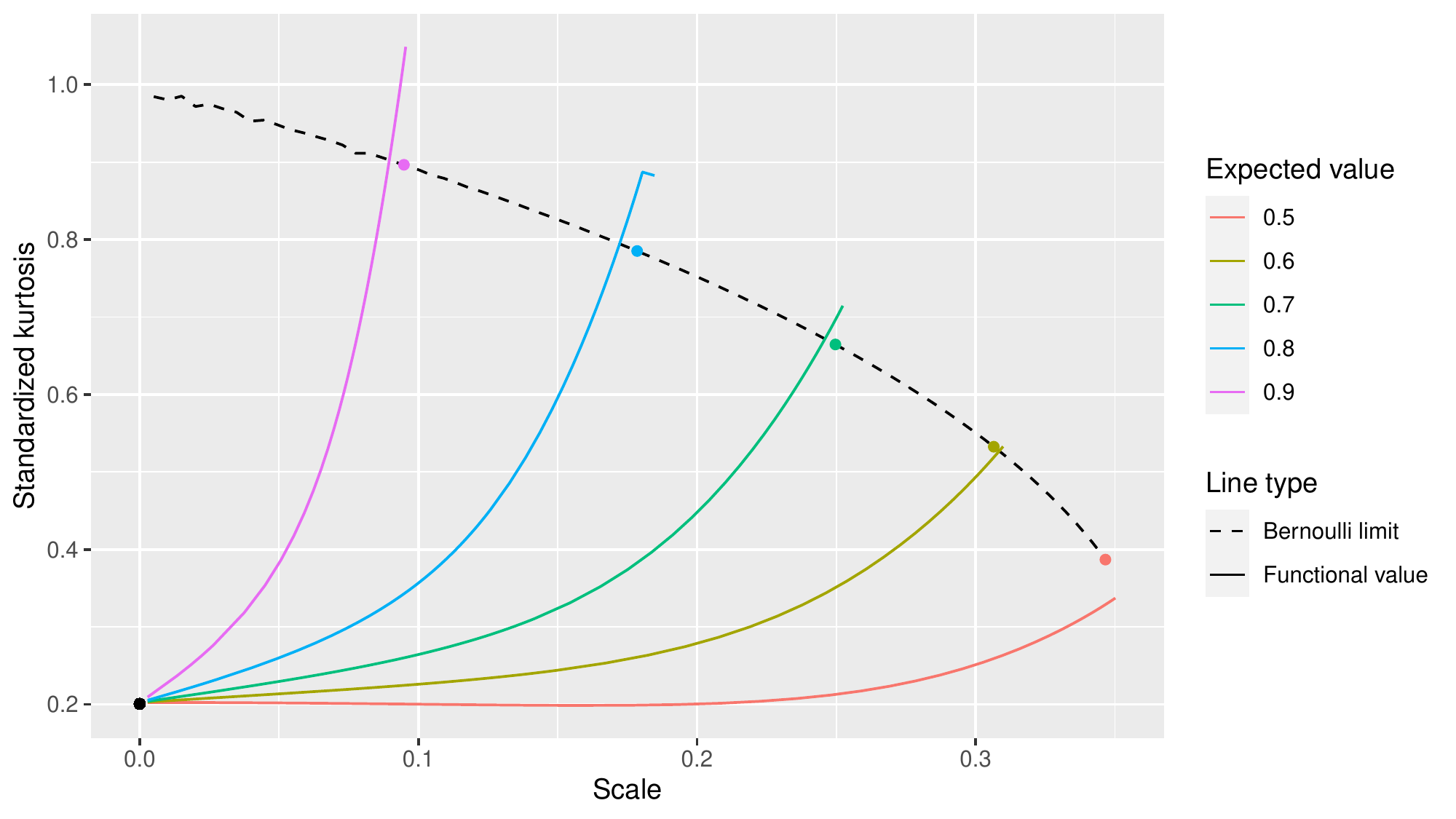}
  \caption{Standardized kurtosis.}
  \label{fig:m4lagbeta}
\end{subfigure}
\caption{Plots of standardized skewness $T_{32}(F)$ (a) and standardized kurtosis $T_{42}(F)$ (b) versus scale $T_2(F)$ for various beta distributions $F\sim B(\psi_1,\psi_2)$ using the Laguerre polynomials and $\eps = 10^{-3}$. The black dots correspond to values of $T_{32}(\text{N}(\cdot,\cdot))$ and $T_{42}(\text{N}(\cdot,\cdot))$. The colored dots illustrate values of $T_{32}(\Be(T_1(F)))$ and $T_{42}(\Be(T_1(F)))$ for each of the five fixed values of $T_1(F)$. The deviations of the colored curves from the expected Bernoulli limits of $T_{32}$ and $T_{42}$ are larger compared to the Hermite polynomials of Figure \ref{herbeta}, due to the asymmetry of the Laguerre polynomials.}
\label{lagbeta}
\end{figure}

Using the canonical parametrization  $B(\psi_1,\psi_2)$  of the beta distribution we generated a grid of parameter values for five fixed choices $(0.5, 0.6, 0.7, 0.8, 0.9)$ of the expected value $T_1(F)=E_F(Y)=\psi_1/(\psi_1+\psi_2)$. For each value pair $(\psi_1, \psi_2)$ we computed the three (standardized) $L$-moments $T_2(F)$, $T_{32}(F)$, and $T_{42}(F)$ according to \eqref{TmDef} and \re{TmlDef}, using all three polynomial series of Section \ref{sec:poly}. 
As $T_2(F)\to 0$, it can be seen from properties of the beta distribution that the standardized version of $F$ will converge to a normal distribution, and consequently $T_{32}(F)\to T_{32}(\text{N}(\cdot,\cdot))$. Since the Hermite and Legendre polynomials generate symmetric collections of $L$-functionals, it follows from \eqref{F0Def} of Section \ref{sec:symm} that $T_{32}(\text{N}(\cdot,\cdot)) = 0$, whereas $T_{32}(\text{N}(\cdot,\cdot)) = -0.340$ for Laguerre polynomials. Similarly, $T_{42}\to T_{42}(\text{N}(\cdot,\cdot))$, where $T_{42}(\text{N}(\cdot,\cdot))=0$ for Hermite polynomials, whereas $T_{42}(\text{N}(\cdot,\cdot))=0.187$ for Legendre polynomials and $T_{42}(\text{N}(\cdot,\cdot))=0.201$ for Laguerre polynomials.
As $T_2(F)$ increases, $F$ will converge to a Bernoulli $\mbox{Be}(T_1(F))$-distribution, so that in the limit $F(\{0\})=1-T_1(F)$ and $F(\{1\})=T_1(F)$. In Figures \ref{legbeta}-\ref{lagbeta} we have plotted $T_{32}(F)$ and $T_{42}(F)$ as functions of $T_2(F)$ for each expected value and polynomial series. It can be seen, for instance, that for the Hermite and Legendre systems, standardized skewness is zero (negative) for a beta distribution with $T_1(F)=0.5$ ($T_1(F)>0.5$). On the other hand, for the Laguerre system, skewness is negative for all beta distributions with $T_1(F)\ge 0.5$. The reason is that skewness of the Laguerre system is quantified in relation to the asymmetric exponential reference distribution.

\section{Bird migration timing analysis}\lb{Sec:RealData}


As an application of the regression methodology of Section \ref{Sec:Regr}, we consider bird migration timing analysis, often called phenological analysis. We will use $L$-functionals and quantile regression to analyze changes over time in location, scale, skewness and kurtosis of bird migration timing distributions. We will also present the varying effects of covariates across the quantiles, and highlight reflections on how the covariates influence the values of the $L$-functionals. Our approach enables analysis of a larger number of distributional aspects than available with canonical methods, such as the one in Lehikoinen et.al.\ (2019).

\subsection{The data}
The data set we will analyze is collected by the Falsterbo Bird Observatory and concerns the species Common Redstart (\textit{Phoenicurus phoenicurus}). This species was selected since it is a bird where the plumages of juveniles, adults, females and males are all distinct. Thus, reliable information on sex and age is available for use as covariates in the analysis.
Birds were captured and ringed by the bird observatory under similar schemes each year during the period 1980-2019, although we will use data from the years 1982-2019, since the identification of juvenile females was not performed the first two years. We will only use data on newly ringed birds. Recaptures between years are rare, but ideally these data should be incorporated as well. 
Due to the extent of the ringing effort each year, it is safe to assume that the whole migration period was covered by the annual sampling window. The covariates and response variable of the data set are summarized in Table \ref{covariatevalues}.

\begin{table}
 \centering
 \begin{tabular}{|c|c|c|}
  \hline
   Variable & Type & Values \\
  \hline
  \hline
  Julian day & Response & Any integer between 80 and 162. \\
  Age & Covariate & \textit{Juvenile} and \textit{adult}.  \\
  Sex & Covariate & \textit{Female} and \textit{male}.  \\
  Year & Covariate & $1982,\ldots,2019$, centered around 2001. \\
 \hline
 \end{tabular}
\caption{Variables used in the phenological data set. We fit the model with \textit{Julian day} as response, and the other variables as covariates. Of primary interest is to study how the arrival distribution varies with \textit{year}, while controlling for the effects of \textit{sex} and \textit{age}.}
\label{covariatevalues}
\end{table}



\subsection{Exploratory visualizations}

\begin{figure}
\includegraphics[width=\linewidth]{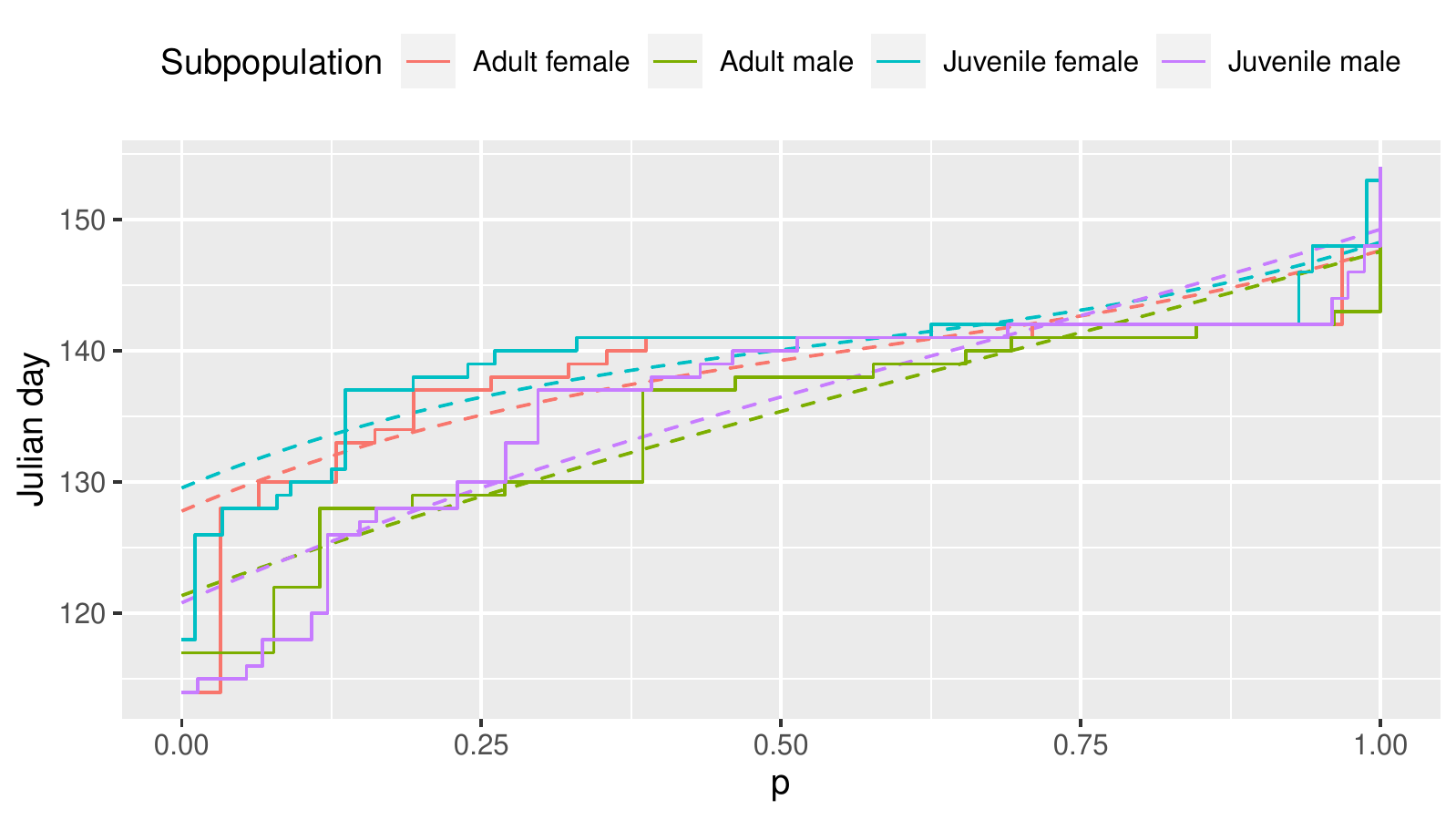}
\caption{The conditional empirical quantile function $\hQ(p|\bx)$ is plotted for each combination of age and sex, for year 2010 (solid lines), as well as the corresponding approximation $\hQ_{\scr{appr}}(p|\bx)$ based on $m=4$ terms (dashed lines), using the Legendre polynomial series.}
\label{redstartsubsets}
\end{figure}

Selecting all the observations from 2010, and splitting the data set into four subsets - one for each combination of age and sex - allows us to compute the first four $L$-moments for each subset this particular year. In Figure \ref{redstartsubsets}, the empirical conditional quantile function 
\beq
 \hQ(p; \bx) = \hF^{-1}(p;\bx) = \inf\{y;\, \hF(y|\bx)\ge p\}
\lb{Qpxhat}
\eeq
for each subpopulation $\bx$ of year 2010, age, and sex is presented, along with the Legendre-based approximation $\hQ_{\scr{appr}}(\cdot|\bx)$ of each empirical conditional quantile function. The function in \re{Qpxhat} is the inverse of the conditional empirical response distribution 
\beq
\hF(y|\bx) = \frac{\sum_{i=1}^n 1(Y_i\le y,\bx_i=\bx)}{\sum_{i=1}^n 1(\bx_i=\bx)},
\lb{edfx}
\eeq
for each covariate vector $\bx$ that appears in the data set. Although the plot in Figure \ref{redstartsubsets} contains data from just one of the years under study, it is possible to trace differences in location, and to some extent differences in scale and skewness, between the subpopulations. 

\subsection{Quantile regression}
Fitting quantile regression models to the Common Redstart data makes it possible to find out whether covariates have different magnitudes of effect for different quantiles. Moreover, we can estimate the first four conditional $L$-functionals from \re{hthx}, and the  standardized conditional $L$-functionals from \re{hthxRatio0}. In particular, we will study how these estimated conditional $L$-functionals change with covariates. 

\subsubsection{Model setup}
We fitted models using two approaches: an identity link and a logit link. As presented in Section 3.1.2, quantiles are preserved under monotone transformations. The response variable is \textit{Julian day} and thus the response values will be located within an interval $[a,b]$, where the sampling window has end points $a=80$ and $b=162$. As presented in Example 8 the \textit{logit}-link function is a monotone transformation. In accordance with Section \ref{Sec:TransLinear}, we fit the model
\beq
 \logit\left( \frac{Q\left( p\mid \bx\right)-a}{b-a}\right) = \bx^T  \boldsymbol{\beta}(p),
\lb{BirdEst1}
\eeq
with $q=4$ parameters in $\bbe(p)$ (corresponding to an intercept and an effect parameter for each one of the three covariates of Table \ref{covariatevalues}). Then we find $\hat{\boldsymbol{\beta}}(p)$ nonparametrically through \re{hbbepLink}, i.e.
\beq
 \hat{\boldsymbol{\beta}}(p) = \arg \min_{\sbfm{b} \in \R^q} \sumiton\rho_p\left(\logit\left( \frac{Y_i-a}{b-a}\right) - \bx_i^T \bb\right).
\lb{BirdEst2}
\eeq
For comparison with \re{BirdEst1}-\re{BirdEst2}, we will also fit a linear quantile regression model \re{QpxLin}. All implementations were made in \textsf{R} (R Core Team, 2021). We chose to implement the objective function ourselves and used \texttt{optim} to find estimates, but it is also possible to use software packages readily available for many types of quantile regression, such as \texttt{qgam} (Fasiolo  et al., 2017) and \texttt{quantreg} (Koenker, 2020).

\subsubsection{Model selection}
When fitting the quantile regression model we may choose any link function that is a monotone transformation. Moreover, when having fitted the model we may choose between the three sets of weight functions $g_m(p)$ presented in Examples \ref{LfuncLegendre}-\ref{Sec:AsymL} and then compute the associated four conditional (standardized) $L$-functionals $\hth(\bx)$. An ideal model selection method should yield the combination of link function and weight functions that gives the best overall approximation of $\hQ(\cdot\mid\bx)$ for all $\bx$. It should also be useful for selecting which covariates to use in the model. To this end, we modify \eqref{RISE} to create
\begin{align} \label{modelRISE}
R^2_{m_0} = 1 - \frac{\int_{\cal X} \int_0^1 \left(\hat{Q}(p\mid\bx)-\hat{Q}_{\text{appr}}(p\mid\bx)\right)^2 \rmd p \rmd \hF_{\sbX}(\bx)}{\int_0^1 (\hat{Q}(p)-E_{\hat{F}}(Y))^2 \rmd p},
\end{align}
where $\hat{F}_{\sbX}$ gives equal weight to all covariates vectors that appear in the data set (regardless of the number of observations with this covariate vector),
\begin{align}
\hQ(\cdot\mid\bx) = \hF^{-1}(\cdot \mid \bx) \\ 
\hQ(p) = \hat{F}^{-1}(p) =  \inf \{y; \hF(y) \ge p\}
\end{align}
and 
\begin{align}
\hF(y) = \int_{\cal X} \hF(y\mid\bx) \rmd \hF_{\sbX}(\bx).
\end{align}
Thus, $E_{\hF}(Y) = T_1(\hF) = T_1(\hF)g_1(p)$, and note that this is not the sample mean, but rather a weighted sample mean, since we weight the unique covariate vectors $\bx$ (not the observations) equally. For the particular analysis of the Common Redstart data, this has the effect of weighting data from each year equally, meaning we do not let years with a larger number of registered birds have a larger influence. Since the number of birds can vary a lot between years, this is deemed advantageous for the purpose of the analysis.

As with \eqref{RISE}, the $m_0$ index specifies the degree of $\hat{Q}_{\text{appr}}$. With $m_0 = 1$, this reduces to a weighted version of the classical $R^2$, and for any $m_0\ge 1$ we get a measure of how much variation in the response is captured by the first $m_0$ conditional $L$-moments. Similarly as in Hössjer (2008), any $R^2_{m_0}$ with $m_0\ge 2$ includes variation in the response not explained be the covariates. In particular, if $R^2_{m_0}$ increases significantly when $m_0$ gets larger, this indicates that the corresponding weight functions $g_{m_0}(p)$ capture an essential part of the variation in the conditional response distributions.  

We fitted models for the identity link and logit link functions and computed $\hth(\bx)$ for all three polynomial weight function sets. When computing the integrals in \re{hthx} and \re{hthxRatio0} we interpolated the grid of $\beta_j(p)$-estimates using the cubic splines of Forsythe (1977). The resulting $R^2_{m_0}$ values are presented in Table \ref{modselect}. Observe that $R_{m_0}^2$ for $m_0>1$ not necessarily is an increasing function of the number of covariates included in the model. This is partly due to the fact that all covariate vectors are weighted equally in the definition of $R_{m_0}^2$, regardless of their number of observations.

\begin{table}
 \centering
 \footnotesize
 \def\arraystretch{1.3}
 \begin{tabular}{|c|c|c|c|c|c|c|}
  \hline
  Covariates & Link function & $g_m$& $R^2_1$ & $R^2_2$ & $R^2_3$ & $R^2_4$\\
  \hline
  \hline
\multirow{6}{*}{Age, sex, year} &&Legendre & 15.29 \% & 56.99\%  & 56.94\% & 58.13\%  \\
&Identity &Hermite  & 15.29\% & 57.55\% & 57.60\% & 57.77\% \\
&&Laguerre & 15.29\% & 50.22\%  & 55.63\% & 56.52\% \\
\cline{2-7}
&&Legendre & 15.30\% & 56.90\% & 56.84\% & 57.94\% \\
&Logit &Hermite  & 15.30\% & 57.34\% & 57.31\% & 57.37\% \\
&&Laguerre & 15.30\% & 49.71\% & 55.27\% & 56.15\% \\
 \hline \hline
\multirow{6}{*}{Year} &&Legendre & 2.89 \% & 67.55\%  & 67.61\% & 70.91\%  \\
&Identity &Hermite  & 2.89\% & 70.88\% & 71.01\% & 71.05\% \\
&&Laguerre & 2.89\% & 60.61\%  & 66.53\% & 67.59\% \\
\cline{2-7}
&&Legendre & 2.95\% & 67.58\% & 67.63\% & 70.93\% \\
&Logit &Hermite  & 2.95\% & 70.85\% & 70.98\% & 71.02\% \\
&&Laguerre & 2.95\% & 60.58\% & 66.53\% & 67.57\% \\
 \hline
 \end{tabular}
\caption{$R^2_{m_0}$-values for two choices of covariates, the identity and logit link functions and the three polynomial systems of weight functions defined in Examples \ref{LfuncLegendre}-\ref{Sec:AsymL}. (A weighted version of) the classical coefficient of determination, $R^2_1$, reflects that more of the variation in the response is explained when \textit{age} and \textit{sex} are added as covariates, compared to having only year as covariate. Notice though that $R^2_{m_0}$ for $m_0 >1$ is higher for a model with year as the only covariate compared to using all three covariates. The reason is that the year only model has an increased number of observations in each covariate specific subset, which leads to a smoother $\hat{Q}(p|\bx)$ with smaller jumps at each discontinuity point. This, in turn, leads to a closer fit of $\hat{Q}_{\text{appr}}(p|\bx)$.}
\label{modselect}
\end{table}

\subsubsection{Covariate and $L$-functional estimates}
\begin{figure} 
\begin{subfigure}{.5\textwidth}
  \centering
  \includegraphics[width=\linewidth]{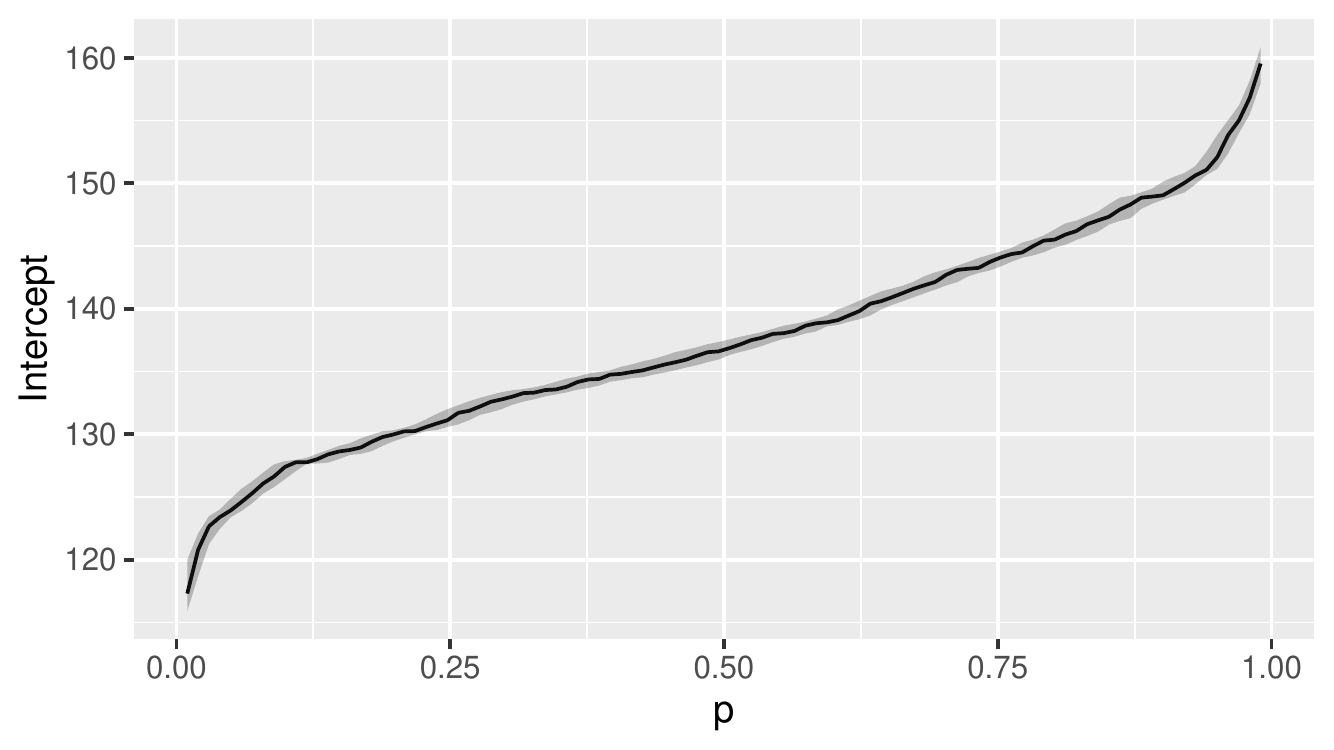}
  \caption{The intercept.}
  \label{fig:intercept}
\end{subfigure}%
\begin{subfigure}{.5\textwidth}
  \centering
   \includegraphics[width=\linewidth]{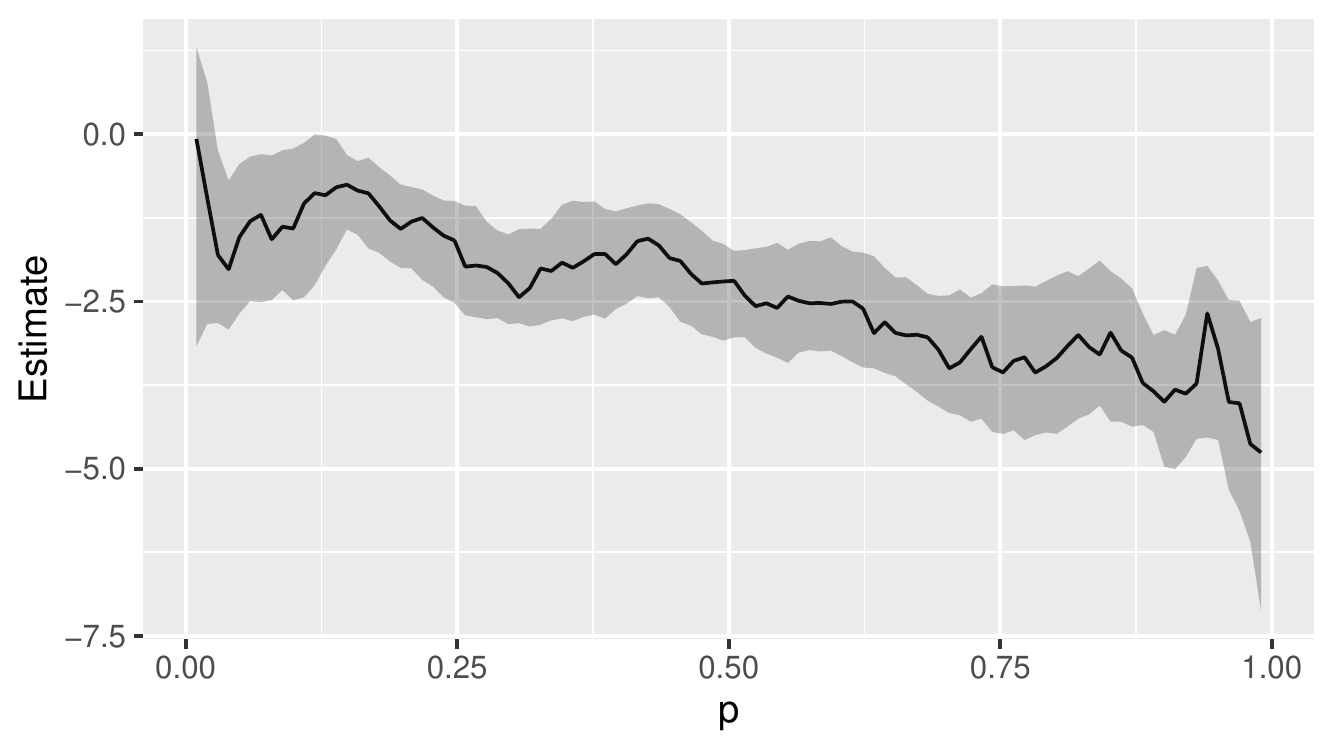}
  \caption{Effect of the bird being adult.}
  \label{fig:adult}
\end{subfigure}
\begin{subfigure}{.5\textwidth}
  \centering
   \includegraphics[width=\linewidth]{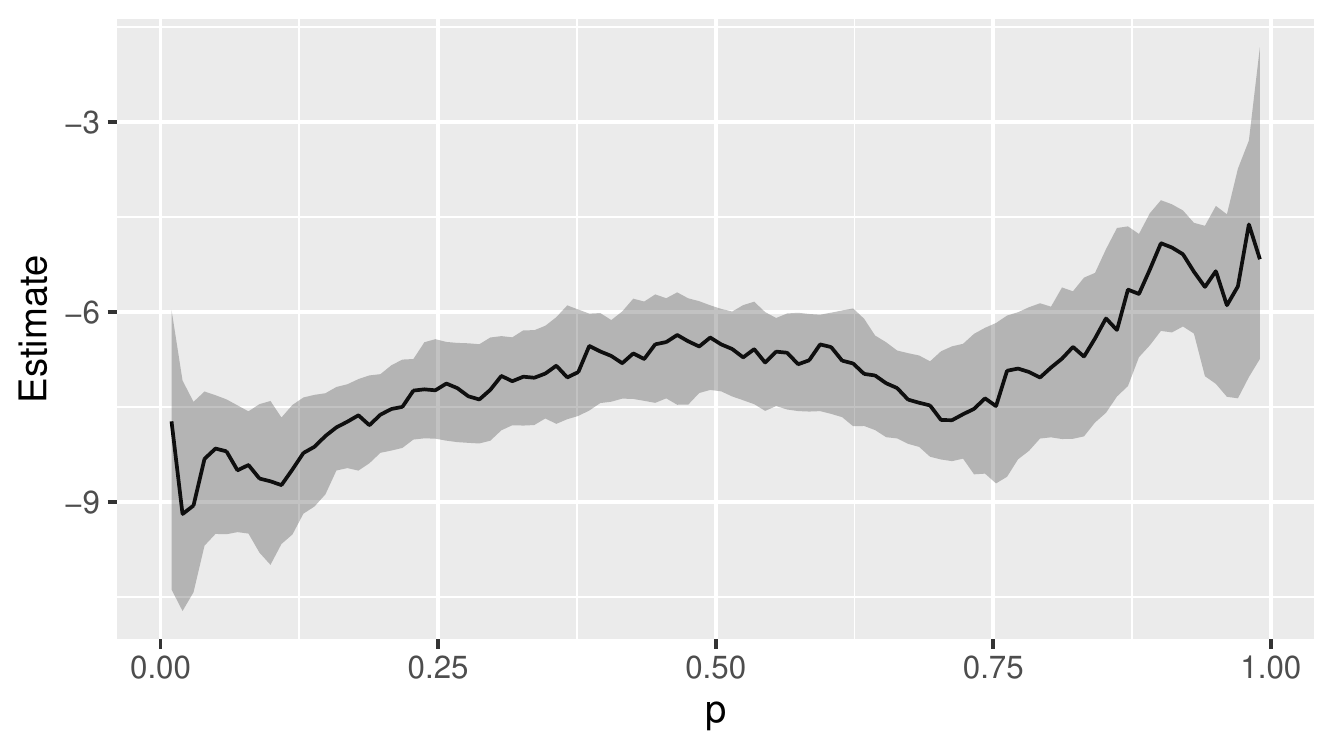}
  \caption{Effect of the bird being male.}
  \label{fig:male}
\end{subfigure}%
\begin{subfigure}{.5\textwidth}
  \centering
   \includegraphics[width=\linewidth]{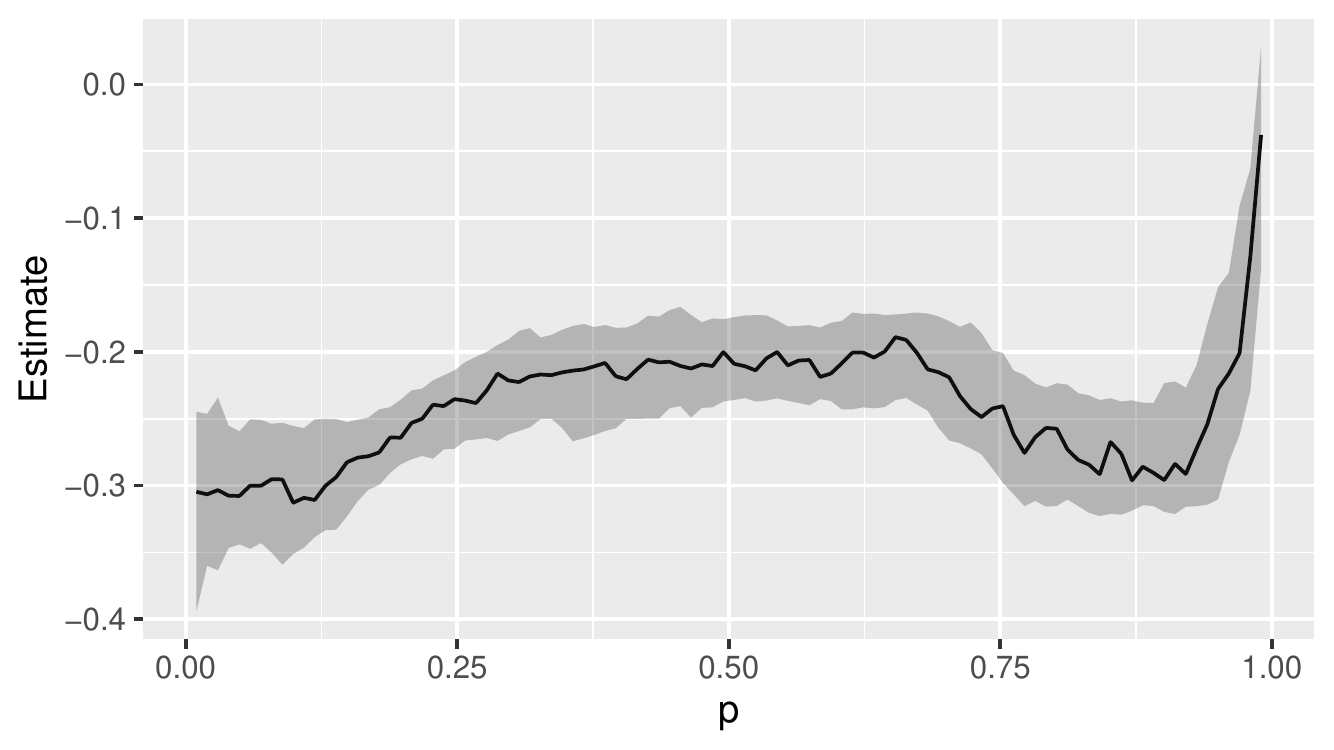}
  \caption{Effect of year.}
  \label{fig:year}
\end{subfigure}
\caption{Plots of parameter estimates for the regression model with the best fit, using the identity link. The black lines represent the $q=4$ estimates $\hbe_1(p),\hbe_2(p),\hbe_3(p),\hbe_4(p)$, where $p$ takes values on a grid of 100 equispaced points in $(0,1)$. The grey ribbon of each subplot consists of approximate $95\%$ pointwise bootstrapped confidence interval for $\{\beta_k(p); \, 0<p<1\}$. Note that in (b), the confidence interval is very broad for the lowest quantiles, which is due to there being extremely few juvenile birds among the early arrivers. In (d), the effect of year quickly approaches 0 in the upper quantiles. This might be due to an increased chance of catching stray birds instead of migrating birds at that time of the year. Lastly, smoothing $\hat{\boldsymbol{\beta}}(p)$ seems like a good idea, given the wigglyness of the estimates.}
\label{betas}
\end{figure}

\begin{figure} 
\begin{subfigure}{.5\linewidth}
  \centering
  \includegraphics[width=\linewidth]{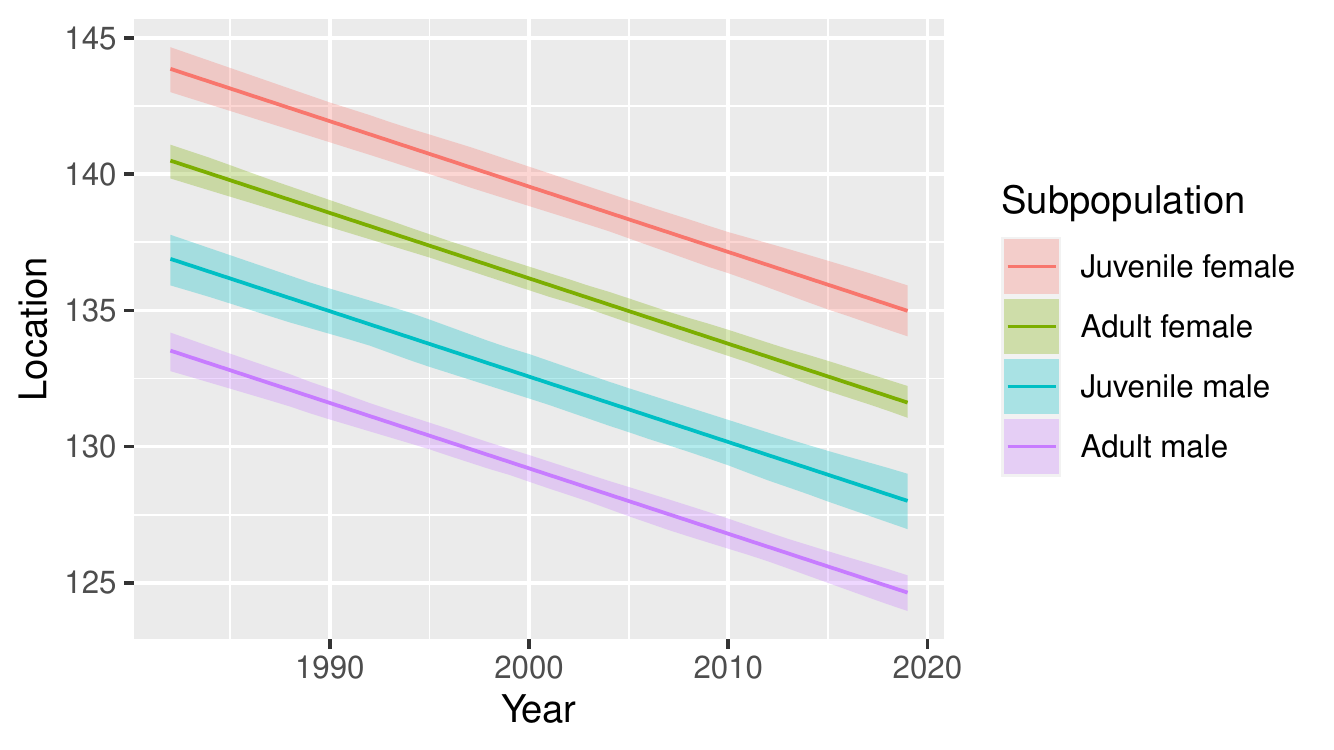}
  \caption{$T_1(\hF_{Y|\sbfm{x}})$}
  \label{fig:T1}
\end{subfigure}%
\begin{subfigure}{.5\linewidth}
  \centering
   \includegraphics[width=\linewidth]{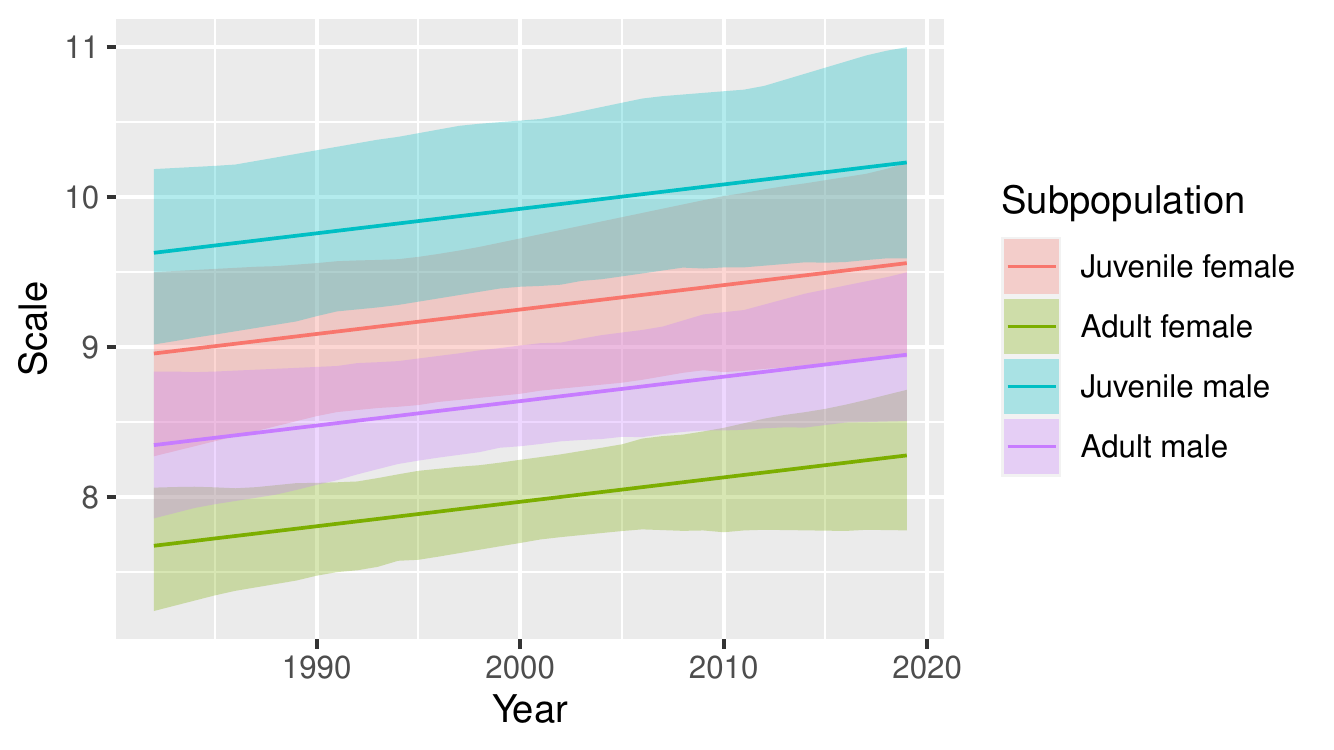}
  \caption{$T_2(\hF_{Y|\sbfm{x}})$}
  \label{fig:T2}
\end{subfigure}
\begin{subfigure}{.5\linewidth}
  \centering
   \includegraphics[width=\linewidth]{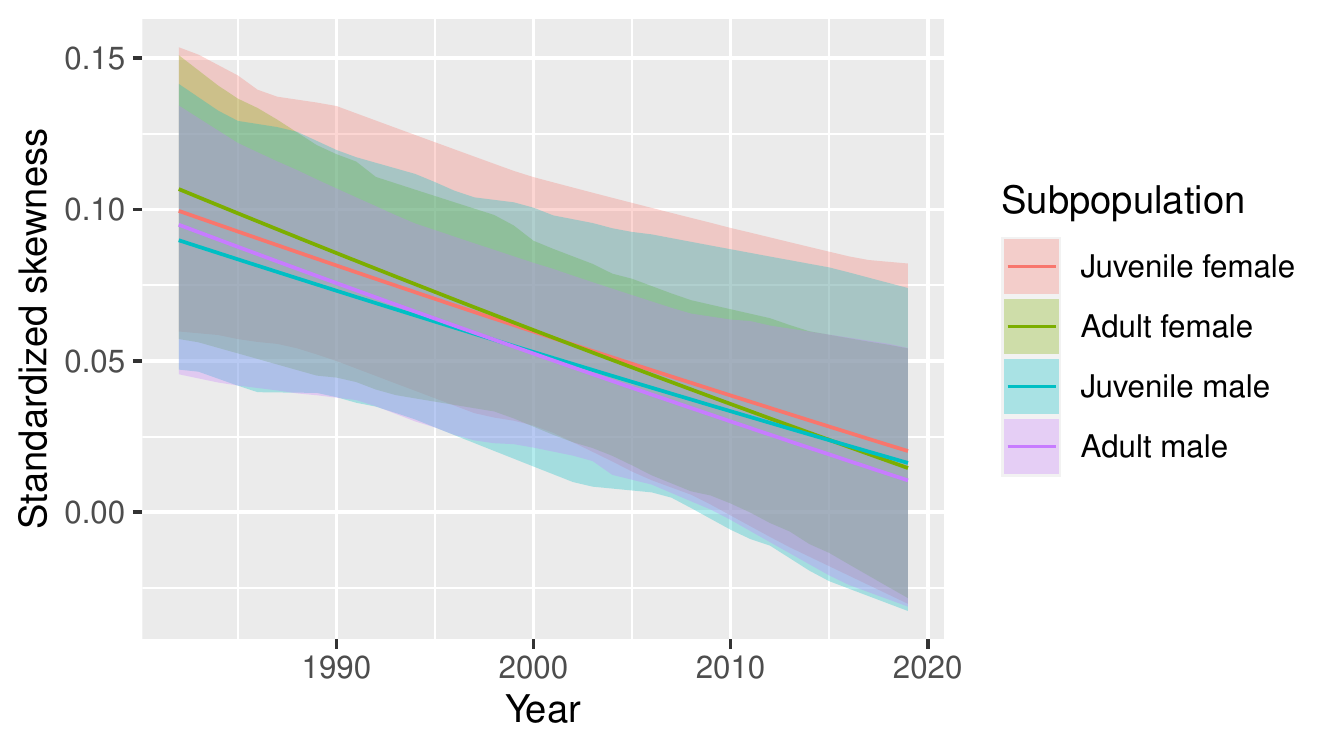}
  \caption{$T_{32}(\hF_{Y|\sbfm{x}})$}
  \label{fig:T32}
\end{subfigure}%
\begin{subfigure}{.5\linewidth}
  \centering
   \includegraphics[width=\linewidth]{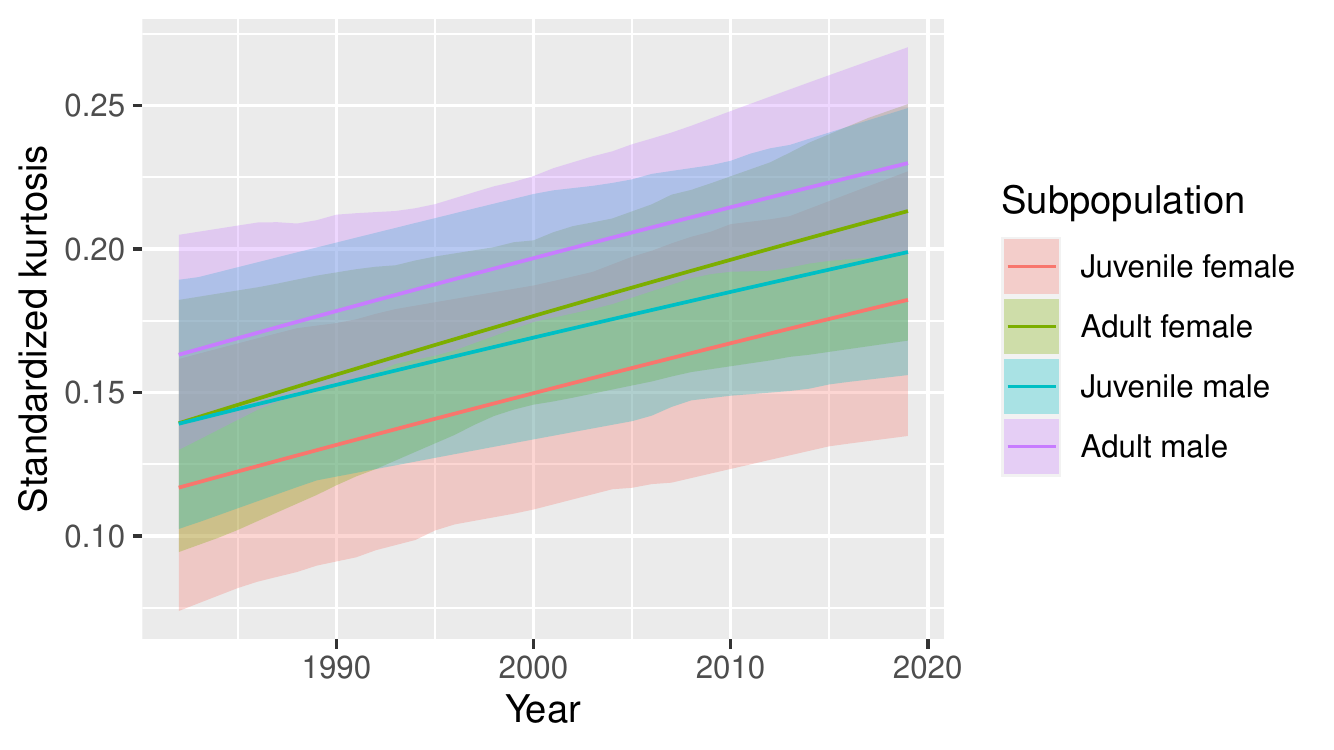}
  \caption{$T_{42}(\hF_{Y|\sbfm{x}})$}
  \label{fig:T42}
\end{subfigure}
\caption{Estimates of the four conditional $L$-functionals of location, scale, standardized skewness and standardized kurtosis for the bird migration data, for each combination $\bx$ of year and the binary covariates age and sex. These conditional $L$-functionals are estimated from a linear quantile regression model, cf.\ \re{hthLin} and \re{hthxRatio}, using the Legendre weight functions $g_m(p)$. A subplot illustrates a specific conditional $L$-functional. It contains four curves, each one of which corresponds to a fixed combination of age and sex, whereas year varies along the horizontal axis. For the location plot (a), all polynomial series would give the same result (the conditional mean), since the first polynomial in each series is 1. For the scale plot (b) the lines correspond to Gini's mean difference times a proportionality constant (cf. Example \ref{LfuncLegendre}). The standardized skewness (c) and standardized kurtosis (d) are measured relative to a uniform distribution, since we use weight functions based on Legendre polynomials. The important information obtained from the figure is qualitative; how the various $L$-functionals change over time, rather than their actual values. From (c) we notice that the arrival distribution shifts from right-skewed towards symmetric over the study period, whereas from (d) we observe that the kurtosis increases. The conclusion would be that the left tail slowly becomes as thick as the right tail, meaning that early arrivers become more frequent over time.}
\label{Ts}
\end{figure}

In this section we will first illustrate how the $q=4$ parameter estimates of \re{BirdEst2} change as a function of the quantile $p$ (cf.\ Figure \ref{betas}) over a grid of 100 values in $(0,1)$ for the case of an identity link function. Next we plot the change in the location, scale, standardized skewness and standardized kurtosis over time, for each subpopulation, i.e. each combination of age and sex. These plots are shown in Figure \ref{Ts}. All estimates include approximate 95\% bootstrapped pointwise confidence intervals, represented by ribbons. The negative estimates of the effect of year in Figure \ref{fig:year}, and the negative slopes over time of the expected value of conditional response distribution, in Figure \ref{fig:T1}, both demonstrate that birds in recent years tend to arrive earlier.

We also use the resampling results to present tables over the Mahalanobis distance between different pairs of covariates' $L$-functionals. For this we need some additional notation. For any subset $I\subset \{1,2,32,42\}$ of order numbers of location, scale, standardized skewness and standardized kurtosis, we introduce the collection
$$
\bth_I(\bx) = (T_m(F_{Y|\sbfm{x}}); \, m\in I)^\top
$$
of conditional (standardized) $L$-moments for covariate vector $\bx$. The corresponding vector of estimated conditional $L$-moments is denoted 
\beq
\hbth_I(\bx) = (T_m(\hF_{Y|\sbfm{x}}); \, m\in I)^\top.
\lb{thI}
\eeq
Recall that each $\bx$ corresponds to a subpopulation (a combination of age group and sex) at a specific time point. Assume that the data set is resampled $B$ times, and let $F^{\ast b}_{Y|\sbfm{x}}$ be the $b$:th resampled conditional response distribution ($b=1,\ldots,B$) for covariate vector $\bx$. The corresponding resampled vector of conditional $L$-moments is   
\beq
\bth_I^{\ast b}(\bx) =  (T_m(F^{\ast b}_{Y|\sbfm{x}}); \, m\in I)^\top.
\lb{thIast}
\eeq 
These $B$ vectors will be scattered around \re{thI}, with an estimated covariance matrix
$$
\hat{\bSi}_I(\bx) = \frac{1}{B}\sum_{b=1}^B (\bth_I^{\ast b}(\bx)-\hbth_I(\bx)) (\bth_I^{\ast b}(\bx)-\hbth_I(\bx))^\top.
$$
In particular, the Mahalanobis distance between \re{thIast} and the center point \re{thI} of the distribution is 
\beq \label{Mahalanobis}
M_I^b(\bx) = \sqrt{(\bth_I^{\ast b}(\bx)-\hbth_I(\bx))^T \hat{\bSi}_I(\bx)^{-1}(\bth_I^{\ast b}(\bx)-\hbth_I(\bx))}.
\eeq
This distance can be interpreted as how many standard deviations away from the center point $\hbth_I(\bx)$ the point $\bth_I^{\ast b}(\bx)$ is. We would however like to measure the distance between the center of the resampled clouds for different pairs $\bx_1$ and $\bx_2$ of covariate vectors, rather than each resampled point's distance to its center. The Mahalanobis distance between $\bx_1$ and $\bx_2$, for the collection $I$ of conditional $L$-moments, is 
\beq
M_I(\bx_1,\bx_2) =\! \sqrt{ (\hbth_I(\bx_1)-\hbth_I(\bx))^T [0.5(\hat{\bSi}_I(\bx_1)+\hat{\bSi}_I(\bx_2))]^{-1}  (\hbth_I(\bx_1)-\hbth_I(\bx)) }.
\lb{ModMah}
\eeq
In Tables \ref{mahal4d} and \ref{mahal4dy} we present values of the modified Mahalanobis distances \re{ModMah} for different combinations of $I$, $\bx_1$, and $\bx_2$. From these tables it can be seen that the location functional is most important for distinguishing the arrival distributions between subpopulations (Table \ref{mahal4d}) and years (Table \ref{mahal4dy}). Although the scale, standardized skewness and standardized heavytailedness functionals are less important, they still help to discriminate even more between the arrival time distributions of these groups.


\begin{table}[ht]
\centering
\footnotesize
\def\arraystretch{1.3}
\begin{tabular}{|r|rrrr|}
\hline
& \multicolumn{4}{|c|}{Subpopulation comparison} \\
\cline{2-5}
Order & Juvenile female & Juvenile male & Juvenile female & Adult female \\
\multicolumn{1}{|r|}{number set $I$} & v. Adult female & v. Adult male & v. Juvenile male & v. Adult male \\
\hline
$\{1\}$ & 12.41 & 10.46 & 19.27 & 32.44 \\ %
$\{2\}$ &  7.14 & 6.21 & 2.97 & 4.36 \\ %
 $\{32\}$ & 0.03 & 0.03 & 0.33 & 0.53 \\ %
 $\{42\}$ & 1.71 & 1.77 & 1.13 & 1.44 \\ %
 $\{1,2\}$ & 14.83 & 11.81 & 19.43 & 33.22 \\ %
 $\{1,32\}$ & 12.59 & 10.47 & 19.53 & 32.44 \\ %
 $\{1,42\}$ & 12.49 & 10.47 & 19.36 & 33.19 \\ %
 $\{2,32\}$ &  7.29 & 6.21 & 2.97 & 4.37 \\ %
 $\{2, 42\}$ & 7.39 & 6.31 & 3.25 & 4.60 \\ %
 $\{32, 42\}$ & 1.75 & 1.77 & 1.14 & 1.54 \\ %
  $\{1,32,42\}$ & 15.36 & 11.82 & 19.65 & 33.23 \\ %
  $\{1, 2, 42\}$ & 14.89 & 11.81 & 19.51 & 33.98 \\ %
  $\{1, 32, 42\}$ & 12.75 & 10.48 & 19.69 & 33.19 \\ %
  $\{2, 32, 42\}$ & 7.47 & 6.31 & 3.26 & 4.61 \\ 
  $\{1,2,32,42\}$ & 15.56 & 11.82 & 19.79 & 33.98 \\
   \hline
\end{tabular}
\caption{The Mahalanobis distance $M_I(\bx_1,\bx_2)$ between different pairs of subpopulations for the year 2000. These distances are computed for different collections $I$ of conditional $L$-moments. In accordance with Figure \ref{Ts} the subpopulations mainly differ in location, although some of them have quite large values of $M_I$ for scale as well. It can be seen that skewness differs the least and kurtosis the second least between subpopulations. }
\label{mahal4d}
\end{table}

\begin{table}[ht]
\centering 
\footnotesize
\def\arraystretch{1.3}
\begin{tabular}{|r|rrrr|}
\hline
& \multicolumn{4}{|c|}{Subpopulation} \\
\cline{2-5}
\multicolumn{1}{|c|}{$I$} & Juvenile female & Juvenile male & Adult female & Adult male \\ 
\hline
$\{1\}$ & 23.66 & 19.46 & 29.81 & 27.11 \\ %
$\{2\}$ & 2.41 & 2.05 & 2.80 & 2.46 \\ %
$\{32\}$ & 3.20 & 3.01 & 4.13 & 3.57 \\ %
$\{42\}$ & 3.27 & 2.86 & 3.73 & 3.17 \\ %
$\{1, 2\}$ & 23.67 & 19.78 & 30.06 & 28.76 \\ %
$\{1, 32\}$ & 24.91 & 20.04 & 30.32 & 27.17 \\ %
$\{1, 42\}$ & 23.67 & 19.48 & 30.27 & 27.17 \\ %
$\{2, 32\}$ & 3.71 & 3.61 & 4.73 & 4.73 \\ %
$\{2, 42\}$ & 4.25 & 3.89 & 4.86 & 4.11 \\ %
$\{32,42\}$ & 4.00 & 4.31 & 5.84 & 5.33 \\ %
$\{1, 2, 32\}$ & 24.92 & 20.34 & 30.49 & 28.94 \\ %
$\{1, 2, 42\}$ & 23.67 & 19.78 & 30.47 & 28.84 \\ %
$\{1, 32, 42\}$ & 25.15 & 20.05 & 30.71 & 27.25 \\ %
$\{2, 32, 42\}$ & 4.61 & 5.05 & 6.39 & 6.35 \\ %
$\{1,2,32,42\}$ & 25.19 & 20.34 & 30.84 & 29.08 \\ %
\hline
\end{tabular}
\caption{The Mahalanobis distance $M_I(\bx_1,\bx_2)$, for each subpopulation, between the years 1982 and 2019. The distances are computed for different collections $I$ of conditional $L$-functionals. As for Table \ref{mahal4d} the main differences are in location. However, standardized skewness and standardized kurtosis both have values of $M_I$ slightly below 4. If we regard 4 standard deviations as significant, these values are on the borderline of demonstrating significant differences in skewness and kurtosis between 1982 and 2019. Notice also that the Mahalanobis distances of $I = \{2,32,42\}$ are larger than 4 for all four subpopulations, indicating that there are other changes other than those in location in the arrival distribution between 1982 and 2019.}
\label{mahal4dy}
\end{table}

\section{Discussion}\lb{Sec:Disc}

In this paper we developed a general theory of $L$-functionals of the response variable distribution of regression models. Based on orthogonal series expansions of the quantile functions of these distributions we generalized the concept of $L$-moments and identified collections $\{T_m(F)\}_{m=1}^4$ of $L$-functionals that correspond to measures of location, scale, unstandardized skewness and unstandardized heavytailedness of the response. Different collections of $L$-functionals were introduced, depending on whether the domain of the response variable is bounded, or unbounded in one or two directions. 
\par\medskip
A number of generalizations of our work is possible. The first extension is to study more systematically the standardized $L$-functionals of skewness, heavytailedness, ... $T_{m2}(F)=T_m(F)/T_2(F)$ when $m\ge 3$. In particular, it is of interest to know which types of stochastic orderings between response distributions (Oja, 1981) these functionals preserve. Second, Takemura (1983) defined orthogonal series expansions of quantile functions for arbitrary reference distributions $F_0$ for which $T_m(F_0)$ vanish when $m\ge 3$. In Examples \ref{LfuncLegendre}-\ref{Sec:AsymL} we considered expansions for uniform, normal, and exponential distributions. It is also possible to define a system \re{gmGen} of $L$-functionals for Weibull, log-logistic, and other reference distributions $F_0$ that are of interest in survival analysis. Third, a general way of robustifying a collection $\{T_m(F)\}$ of $L$-functionals is to introduce a weight function $w$ that downweights contributions from the lower and upper tails of $F$. The weight density 
$$
g_m(p) = w(p)P_{m-1}(Q_0(p))
$$
of order $m$ is a generalization of \re{gmGen}, using a variable rather than a constant weight function $w(p)\equiv 1$. For instance, the trimmed $L$-moments of Elamir and Seheult (2003) correspond to a weight function $w(p)=(1-p)^t p^t$ for some positive integer $t$. Notice however that in general the orthogonality property \re{gOrth} is lost for such a system of weight functions. Fourth, recursive estimation of quantiles (Stephanou et al., 2017) could be extended to online estimation of $L$-functionals. Fifth, regression quantiles have been used for time series (Cai and Zu, 2008, White et al., 2008) in order to estimate quantiles and robust measures of skewness/kurtosis of predictive distributions. It is of interest to analyze (ratios of) $L$-functionals of such predictive distributions. Sixth, measures of location, scale, skewness and heavytailedness of multivariate distributions can be defined in terms of multivariate $L$-statistics (Liu, 1990, Liu et al. 1999, Zuo et al., 2004, and Dang et al., 2009). In this context it is of interest to study different polynomial systems of $L$-functionals and their order numbers.

\section*{References}

Bennett, C.A. (1952). Asymptotic properties of ideal linear estimators. Ph.D.\ dissertation, University of Michigan. 
\par\smallskip
Bennett, S. (1983). Analysis of survival data by the proportional odds model. {\sl Statistics in Medicine} {\bf 2}, 273-277. 
\par\smallskip
Bickel, P.J. (1965). On some robust estimators of location. {\sl Annals of Mathematical Statistics} {\bf 36}, 847-858. 
\par\smallskip
Bickel, P.J. (1973). On some analogues of linear combinations of order statistics in the linear model. {\sl Annals of Statistics} {\bf 1}, 597-616. 
\par\smallskip
Bickel, P.J.\ and Lehmann, E.L. (1975). Descriptive statistics for nonparametric models II: Location. {\sl Annals of Statistics} {\bf 3}, 1045-1069. 
\par\smallskip
Bickel, P.J.\ and Lehmann, E.L. (1976). Descriptive statistics for nonparametric models III: Dispersion. {\sl Annals of Statistics} {\bf 4}, 1139-1158. 
\par\smallskip
Billingsley, P. (1999). {\sl Convergence of Probability Measures}, 2nd ed. Wiley, New York. 
\par\smallskip
Bottai, M., Cai, B.\ and McKeown, R.E. (2010). Logistic quantile regression for bounded outcomes. {\sl Statistical Medicine} {\bf 29}, 309-317. 
\par\smallskip
Bowley, A.L. (1920). {\sl Elements of Statistics}, Schribner's, New York. 
\par\smallskip
Box, G.E.P.\ and Cox, D.R. (1964). An analysis of transformations. {\sl Journal of the American Statistical Association} {\bf 26}, 211-252. 
\par\smallskip
Cai, Z.\ and Xu, X. (2008). Nonparametric quantile estimations for dynamic smooth coefficient models. {\sl Journal of the American Statistical Association} {\bf 103}, 1595-1608. This reference was removed since it deals time series data. 
\par\smallskip
Chauduri, P. (1991). Global nonparametric estimation of conditional quantile functions and their derivatives. {\sl Journal of Multivariate Analysis} {\bf 39}, 246-269. 
\par\smallskip
Chauduri, P.\ and Loh, W.-L. (2002). Nonparametric estimation of conditional quantiles using quantile regression trees. {\sl Bernoulli} {\bf 8}, 561-576. 
\par\smallskip
Chen, X., Linton, O.\ and van Kellegom, I. (2003). Estimation of semiparametric functions when the criterion function is not smooth. {\sl Econometrica} {\bf 71}, 1591-1608. 
\par\smallskip
Chernoff, H., Gastwirth, J.L.\ and Jones, M.V. (1967). Asymptotic distribution of linear combinations of functions of order statistics with applications to estimation. {\sl Annals of Mathematical Statistics} {\bf 38}, 52-72. 
\par\smallskip
Chissom, B.S. (1970). Interpretation of the kurtosis statistics. {\sl American Statistician} {\bf 24}(4), 19-23. 
\par\smallskip
Cox, D.R. (1972). Regression models and life tables. {\sl Journal of the Royal Statistical Society, Series B}, {\bf 34}, 187-220.
\par\smallskip
Dang, X., Serfling, R., and Zhou, W. (2009). Influence functions of some depth functions, and application to depth-weighted L-statistics. {\sl Journal of Nonparametric Statistics}, {\bf 21}(1), 49-66.
\par\smallskip
Doksum, K.\ and Gasko, M. (1990). On a correspondence between models in binary regression and survival analysis. {\sl International Statistical Review} {\bf 58}, 243-252. 
\par\smallskip
Elamir, E.A.\ and Seheult, A.H. (2003). Trimmed L-moments. {\sl Computational Statistics and Data Analysis} {\bf 43}, 299-314. 
\par\smallskip
Efron, B. (1991). Regression percentiles using asymmetric squared error loss. {\sl Statistica Sinica} {\bf 1}(1), 93-125. 
\par\smallskip
Efron, B.\ and Tibshirani, R. (1993). {\sl An introduction to the bootstrap}, Chapman and Hall, Boca-Raton, FL. 
\par\smallskip
Fasiolo M., Goude Y., Nedellec R. and Wood S. N. (2017). Fast calibrated
  additive quantile regression. URL: https://arxiv.org/abs/1707.03307
\par\smallskip  
Fisher, R.A.\ and Cornish, E.A. (1960). The percentile points of distributions having known cumulants. {\sl Technometrics} {\bf 2}, 209-225.  
\par\smallskip
Forsythe, G. E. (1977). Computer methods for mathematical computations. Prentice-Hall series in automatic computation, 259.
\par\smallskip
Frumento, P.\ and Bottai, M. (2016). Parametric modeling of quantile regression coefficient models. {\sl Biometrics} {\bf 72}, 74-84. 
\par\smallskip
Frumento, P.\ and Bottai, M. (2017). An estimating equation for censored and truncated quantile regression. {\sl Computational Statistics and Data Analysis} {\bf 113}, 53-63.  
\par\smallskip
Galton, F. (1883). {\sl Enquiry to human faculty and its development}. London, MacMillan.  
\par\smallskip
Garcia, V.J., Martel-Escobar, M.\ and V\'{a}zquez-Polo, F.J. (2018). A note on ordering probability distributions by skewness. {\sl Symmetry} {\bf 10}(7), 286. 
\par\smallskip
Garc\'{i}a-Pareja, C.\ and Bottai, M. (2018). On mean decomposition for summarizing conditional distributions. {\sl Statistics} 2018:7:e208. 
\par\smallskip
Garc\'{i}a-Pareja, C., Santacatterian, M., Ekstr\"{o}m, A.M.\ and Bottai, M. (2019). In {\sl Topics in mathematical statistics for medical applications: summary measures and exact simulation of diffusions.} PhD Thesis, Institute of Environmental Medicine, Karolinska Institutet, Stockholm.  
\par\smallskip
Gelfand, A.E., Ghosh, S.K., Christiansen, C., Somerai, S.B.\ and McLaughlin, T.J. (2000). Proportional hazards models: a latent competing risk approach. {\sl Applied Statistics} {\bf 49}, 385-397. 
\par\smallskip
Gilchrist, W. (2000). {\sl Statistical Modelling with Quantile Functions}. Chapman and Hall. 
\par\smallskip
Gilchrist, W. (2007). Modelling and fitting quantile distributions and regressions. {\sl American Journal of Mathematical and Management Sciences} {\bf 27}, 401-439. 
\par\smallskip
Greenwood, J.A., Landwehr, J.M., Matalas, N.C.\ and Wallis, J.R. (1979) Probability weighted moments: definition and relation to parameters of several distributions expressable in inverse form. {\sl Wat.\ Resour.\ Res.} {\bf 15}, 1049-1054.
\par\smallskip
Groeneveld, R.A.\ and Meeden, G. (1984). Measuring skewness and kurtosis. {\sl The Statistician} {\bf 33}, 391-399. 
\par\smallskip
Gutenbrunner, C.\ and Jure\v{c}kov\'{a}, J. (1992). Regression quantile and regression rank score processes in the linear model and derived statistics. {\sl Annals of Statistics} {\bf 20}, 305-330.  
\par\smallskip
Haeusler, E.\ and Teugels, J.L. (1985). On asymptotic normality of Hill's estimator for the exponent of regular variation. {\sl Annals of Statistics} {\bf 13}(2), 743-756.
\par\smallskip
He, X.\ and Shi, P. (1994). Convergence rate of $B$-spline estimators of nonparametric quantile functions. {\sl Journal of Nonparametric Statistics} {\bf 33}, 299-308. 
\par\smallskip
Hill, B.M. (1975). A simple general approach to inference about the tail of a distribution. {\sl Annals of Statistics} {\bf 3}, 1163-1174. 
\par\smallskip
Hinkley, D.V. (1975). On power transformations to symmetry. {\sl Biometrika} {\bf 62}, 101-111.
\par\smallskip
Hogg, R.V. (1972). More light on the kurtosis and related statistics. {\sl Journal of the American Statistical Association} {\bf 67}, 422-424.
\par\smallskip
Hogg, R.V. (1974). Adaptive robust procedures: A partial review of some suggestions for future applications and theory. {\sl Journal of the American Statistical Association} {\bf 67}, 422-424.
\par\smallskip
Hosking, J.R.M. (1990). L-moments: analysis and estimation of distributions using linear combinations or order statistics. {\sl Journal of the Royal Statistical Society Ser.\ B} {\bf 52}(1), 105-124 
\par\smallskip
Hosking, J.R.M. (1992). Moments or $L$ moments? An example comparing two measures of distributional shape. {\sl The American Statistician} {\bf 46}(3), 186-189. 
\par\smallskip
Hosking, J.R.M. (2006). On the characterization of distributions by their L-moments. {\sl Journal of Statistical Inference and Planning} {\bf 136}(1), 193-198.
\par\smallskip
Hössjer, O. (2008). On the coefficient of determination for mixed regression models.
{\sl Journal of Statistical Planning and Inference} {\bf 138}, 3022-3038. 
\par\smallskip
Jung, J. (1955). On linear estimates defined by continuous weight function. {\sl Arkiv f{\"u}r Matematik} {\bd Bd 3}, 199-209. 
\par\smallskip
Kalbfleish, J.D.\ and Prentice, R.L. (2002). {\sl Statistical Analysis of Failure Time Data} (2nd ed.), Wiley, Noboken, NJ. 
\par\smallskip
Karian, Z.A.\ and Dudewicz, E.J. (2000). {\sl Fitting Statistical Distributions: The Generalized Lambda Distribution and Generalized Bootstrap Methods}. CRC Press, Boca Raton, Florida.  
\par\smallskip
Karvanen, J. (2006). Estimation of quantile mixtures via L-moments and trimmed L-moments. {\sl Computational Statistics and Data Analysis} {\bf 51}(2), 947-959.
\par\smallskip
Karvanen, J.\ and Nuutinen, A. (2008). Characterizing the generalized lambda distribution by L-moments. {\sl Computational Statistics and Data Analysis} {\bf 52}, 1971-1983. 
\par\smallskip
Kim, M.-O. (2007). Quantile regression with varying coefficients. {\sl Annals of Statistics} {\bf 35}(1), 92-108. 
\par\smallskip
Kim, T.-H.\ and White, H. (2004). On more robust estimation of skewness and kurtosis. {\sl Finance Research Letters} {\bf 1}, 56-73. 
\par\smallskip
Koenker, R. (2005). {\sl Quantile Regression}. Cambridge University Press, Cambridge.
\par\smallskip 
Koenker, R. (2020). \texttt{quantreg}: Quantile Regression. R package version 5.73. https://CRAN.R-project.org/package=quantreg
\par\smallskip 
Koenker, R.\ and Bassett Jr, G. (1978). Regression quantiles. {\sl Econometrica} 46(1), 33-50.
\par\smallskip 
Koenker, R.\ and Bassett Jr, G. (1982). Robust tests for heteroscedasticity based on regression quantiles. {\sl Econometrica} {\bf 50}, 43-61. 
\par\smallskip
Koenker, R.\ and Geling, O. (2001). Reappraising medfly longevity: A quantile regression survival analysis. {\sl Journal of the American Statistical Association} {\bf 96}, 458-468. 
\par\smallskip
Koenker, R.\ and Hallock, K. (2001). Quantile regression: An introduction. {\sl Journal of Economic Perspectives} {\bf 15}, 143-156.
\par\smallskip
Koenker, R., Ng, N.\ and Portnoy, S. (1994). Quantile smoothing splines. {\sl Biometrika} {\bf 81}, 673-680.
\par\smallskip
Koenker, R.\ and Portnoy, S. (1987). $L$-estimation for linear models. {\sl Journal of the American Statistical Association} {\bf 82}, 85-1857. 
\par\smallskip
Koenker, R.\ and Zhao, Q. (1994). $L$-estimation for linear heteroscedastic models. {\sl Journal of Nonparametric Statistics} {\bf 3}, 223-235. 
\par\smallskip
Lehikoinen, A., Lind\'{e}n, A., Karlsson, M., Andersson, A., Crewe, T. L., Dunn, E. H., Gregory, G., Karlsson, L., Kristiansen, V., Mackenzie, S. \ and others (2019). Phenology of the avian spring migratory passage in Europe and North America: Asymmetric advancement in time and increase in duration. {\sl Ecological Indicators} {\bf 101}, 985-999.
\par\smallskip
Leng, C.\ and Tong, X. (2013). A quantile regression estimator for censored data. {\sl Bernoulli} {\bf 19}(1), 344-361. 
\par\smallskip
Lind\'{e}n, A., Meller, K. \ and Knape, J (2017). An empirical comparison of models for the phenology of bird migration. {\sl Journal of Avian Biology} {\bf 48}, 255-265.
\par\smallskip
Lindgren, A. (1997). Quantile regression with censored data using generalized $L_1$ minimization. {\sl Computational Statistics and Data Analysis} {\bf 23}, 509-524. 
\par\smallskip
Liu, R. Y. (1990). On a notion of data depth based on random simplices. {\sl The Annals of Statistics}, 405-414.
\par\smallskip
Liu, Y.\ and Bottai, M. (2009). Mixed-effects models for conditional quantiles with longitudinal data. {\sl International Journal of Biostatistics} {\bf 5}(1), Article 28.  
\par\smallskip
Liu, R. Y., Parelius, J. M., and Singh, K. (1999). Multivariate analysis by data depth: descriptive statistics, graphics and inference (with discussion and a rejoinder by liu and singh). {\sl The annals of statistics}, {\bf 27}(3), 783-858.
\par\smallskip
McCullagh, P.\ and Nelder, J.A. (1989). {\sl Generalized Linear Models}, Chapman and Hall/CRC, New York. 
\par\smallskip
McGillivray, H.L. (1986). Skewness and asymmetry: Measures and orderings. {\sl Annals of Statistics} {\bf 14}, 994-1011. 
\par\smallskip
Moore, D.S. (1968). {\sl Annals of Mathematical Statistics} {\bf 39}(1), 263-265. 
\par\smallskip
Moors, J.J.A. (1988). Q quantile alternative for kurtosis. {\sl The Statistician} {\bf 37}, 25-32.
\par\smallskip
Mosteller, F. (1946). On some useful ``inefficient'' statistics. {\sl Annals of Mathematical Statistics} {\bf 17}, 377-408. 
\par\smallskip
Mu, Y.\ and He, X. (2007). Power transformation toward a linear regression quantile. {\sl Journal of the American Statistical Association} {\bf 102}, 269-279.
\par\smallskip
Mudholkar, G.S.\ and Hutson, A.D. (1998). $LQ$-moments: analogs of $L$-mo\-ments. {\sl Journal of Statistical Planning and Inference} {\bf 71},
191-208.
\par\smallskip
Neocleous, T., Vanden Branden, K.\ and Portnoy, S. (2006). Correction to Portnoy (2003). {\sl Journal of the American Statistical Association} {\bf 101}, 860-861. 
\par\smallskip
Newey, W.K.\ and McFadden, D.I. (1994). Large sample estimation and hypothesis testing In: Engle, R.F., McFadden, D.I.\ (eds.) {\sl Handbook of Econometics} {\bf 4}, North-Holland, Amsterdam, pp.\ 2211-2245.  
\par\smallskip
Oberhofer, W. (1982). The consistency of nonlinear regression minimizing the $L_1$ norm. {\sl Annals of Statistics} {\bf 10}, 316-319. 
\par\smallskip
Oja, H. (1981). On location, scale, skewness and kurtosis of univariate distributions. {\sl Scandinavian Journal of Statistics} {\bf 8}, 154-168. 
\par\smallskip
Okagbue, H.I., Adamu, M.O.\ and Anake, T.A. (2019). Quantile mechanics: Issues arising from critical review. {\sl International Journal of Advanced and Applied Sciences} {\bf 6}(1), 9-23. 
\par\smallskip
Olkin, I.\ and Pukelsheim, F. (1982). The distance between two random vectors with given dispersion matrices. {\sl Linear Algebra Appl.} {\bf 48}, 257-263. 
\par\smallskip
Parzen, E. (1979). Nonparametric statistical modelling (with commmets). {\sl Journal of the American Statistical Association} {\bf 74}, 105-131.
\par\smallskip
Peng, L.\ and Huang, Y. (2008). Survival analysis with quantile regression models. {\sl Journal of the American Statistical Association} {\bf103}, 637-649.
\par\smallskip
Portnoy, S. (2003). Censored quantile regression. {\sl Journal of the American Statistical Association} {\bf 98}, 1001-1012. 
\par\smallskip
Powell, J.L. (1986). Censored regression quantiles. {\sl Journal of Econometrics} {\bf 32}, 143-155.
\par\smallskip
R Core Team (2021). R: A language and environment for statistical
  computing. {\sl R Foundation for Statistical Computing}, Vienna, Austria. \\
  https://www.R-project.org/.
\par\smallskip
Rousseeuw, P.J.\ and Leroy, A. (1987). {\sl Robust regression and outlier detection}, Wiley, New York. 
\par\smallskip
Royston, P.\ and Parmar, M.K.B. (2002). Flexible parametric propor\-tional-hazards an proportional-odds models for censored survival data, with application to prognostic modelling and estimation of treatment effects. {\sl Statistics in Medicine} {\bf 21}, 2175-2197. 
\par\smallskip
Ruppert, D.\ and Carroll, R. (1980). Trimmed least squares estimation in the linear model. {\sl Journal of the American Statistical Association} {\bf 75}, 828-838. 
\par\smallskip
Serfling, R. (1980). {\sl Asymptotic Theory of Statistics}. John Wiley, New York. 
\par\smallskip
Sheather, S.J\ and Marron. J.S. (1990). Kernel quantile estimators. {\sl Journal of the American Statistical Association} {\bf 85}, 410-416. 
\par\smallskip
Sillito, G. (1969). Derivation of approximants to the inverse distribution function of a continuous univariate population from the order statistics of a
sample. {\sl Biometrika} {\bf 56}(3), 641-650. 
\par\smallskip
Stephanou, M., Varughese, M.\ and Mcdonald, I. (2017). Sequential quantiles via Hermite series density estimation. arXiv:15017.05073v2. 
\par\smallskip
Stigler, S. (1977). Do robust estimators work with real data? {\sl Annals of Statistics} {\bf 5}, 1055-1098. 
\par\smallskip
Takamura, A. (1983). Orthogonal expansion of quantile function and components of the shapiro-statistic, Report No.\ TR-8, Stanford University, Department of Statistics, California, USA. 
\par\smallskip
Tukey, J.W.\ and McLaughlin, G. (1963). Less vulnerable confidence and significance procedures for location based on a single sample. {\sl Sankya Ser.\ A} {\bf 35}, 331-352.
\par\smallskip
Wang, H.J.\ and Wang, L. (2009). Locally weighted censored quantile regression. {\sl Journal of the American Statistical Association} {\bf 104}, 1117-1128.
\par\smallskip
van Zwet, W.R. (1964). {\sl Convex Transformations of Random Variables}, Math.\ Centrum, Amsterdam. 
\par\smallskip
Welsh, A.H. (1987). The trimmed mean in the linear model. {\sl Annals of Statistics} {\bf 15}(1), 20-36.
\par\smallskip
Welsh, A.H.\ and Morrison, H.L. (1990). Robust L-estimation of scale with an application to astronomy. {\sl Journal of the American Statistical Association} {\bf 85}, 729-743.  
\par\smallskip
White, H., Tae-Hwan, K.\ and Manganelli, S. (2008). Modeling autoregressive conditional skewness and kurtosis with multi-quantile CAViaR. European Central Bank, Working paper 957, November 2008. 
\par\smallskip
Yang, S. (1999). Censored median regression using weighted empirical survival and hazard functions. {\sl Journal of the American Statistical Association} {\bf 94}, 137-145.
\par\smallskip
Ying, Z., Jung, S.H. and Wei, L.J. (1995). Survival analysis with median regression models. {\sl Journal of the American Association} {\bf 90}, 178-184.
\par\smallskip
Younes, N.\ and Lachin, J. (1997). Link-based models for survival data with interval and continuous time censoring. {\sl Biometrics} {\bf 53}, 1199-1211.
\par\smallskip
Zuo, Y., Cui, H., and Young, D. (2004). Influence function and maximum bias of projection depth based estimators. {\sl The Annals of Statistics}, {\bf 32}(1), 189-218.
\end{document}